\documentclass[preprint,3p]{elsarticle}

\usepackage{amssymb}
 \usepackage{amsthm}
\usepackage{latexsym}
\usepackage{amsmath}
\usepackage{amssymb}
\usepackage{diagbox}
\usepackage{slashbox}
\usepackage{color}
\usepackage{subfigure}
\usepackage{graphicx}
\usepackage{epstopdf}
\usepackage{arydshln}

\newcommand{\udots}{\mathinner{\mskip1mu\raise1pt\vbox{\kern7pt\hbox{.}}
\mskip2mu\raise4pt\hbox{.}\mskip2mu\raise7pt\hbox{.}\mskip1mu}}
\usepackage{amsmath}

\newtheorem{theorem}{Theorem}

\newtheorem{proposition}{Proposition}
\begin{document}

\begin{frontmatter}



\title{ Multivariate doubly truncated moments for generalized skew-elliptical distributions with application to multivariate tail conditional risk measures
}

\author{Baishuai  Zuo}
\author{Chuancun Yin\corref{cor1}}
\cortext[cor1]{Corresponding author.}
 \ead{ccyin@qfnu.edu.cn}

\address{School of Statistics and Data Science, Qufu Normal University, Qufu, Shandong 273165, P. R. China}

\begin{abstract}
 In this paper, we focus on multivariate doubly truncated first two moments
of generalized skew-elliptical (GSE) distributions and derive explicit expressions for them.
  It includes many useful distributions, for examples, generalized skew-normal (GSN), generalized skew-Laplace (GSLa), generalized skew-logistic (GSLo) and generalized skew student-$t$ (GSSt) distributions,  all
as special cases. We also give formulas of multivariate doubly truncated expectation and covariance for GSE distributions. As applications, we show the results of multivariate tail conditional expectation (MTCE) and multivariate tail covariance (MTCov) for GSE distributions.
\end{abstract}


\end{frontmatter}

\baselineskip=20pt

\section{Introduction}
Truncated moments' expansions are applied to many areas including the design of experiment (see Thompson, 1976), robust estimation (see Cuesta-Albertos et al., 2008), outlier detections (see, e.g., Riani et al., 2009; Cerioli, 2010), robust regression (see Torti et al., 2012), robust detection (see Cerioli et al. (2014)), statistical estates' estimation (see Shi et al., 2014), risk averse selection (see Hanasusanto et al., 2015), entropy computation and application (see, e.g., Milev et al., 2012; Zellinger and Moser, 2021). Therefore, the related research of truncated moment are developed to  different distributions by many scholars.

 Since Tallis (1961) has given an explicit formula for the first two moments of a lower truncated multivariate standard normal distribution by moment generating function, Amemiya (1974) and Lee (1979) used
Tallis' results (see Tallis, 1961) for extending Tobin's model (see Tobin, 1958) to the multivariate regression and simultaneous equations models when the dependent variables are truncated normal.
Lien (1985) provided the expressions for the moments of lower truncated bivariate log-normal distributions. Manjunath and Wilhelm (2012) computed the first and second moments for the rectangularly double truncated multivariate normal density, and extended the derivation of Tallis (1961) to general $\mu$, $\Sigma$ and for doubly truncation.
Arismendi (2013) derived formulae for the higher order tail moments of the lower truncated multivariate standard normal, Student's $t$, lognormal and a finite-mixture of multivariate normal distributions.
Arismendi and Broda (2017) deeply derived multivariate elliptical lower truncated moment generating function, and first second-order moments of quadratic forms of the multivariate normal, Student and generalized hyperbolic distributions. Moreover, Ho et al. (2012) presented
general formulae for computing the first two moments of the  truncated
multivariate $t$ distribution under
the doubly truncation. Recently, Kan and Robotti (2017) provided expressions of the moments for folded and doubly truncated multivariate normal distribution. Galarza et al. (2021) and Morales et al. (2022) generalized moments of folded and doubly truncated to multivariate Student-$t$ and extended skew-normal distributions, respectively. Also, Ogasawara (2021) derived a non-recursive formula for various moments of the multivariate normal distribution with sectional truncation, and introduced the importance of truncated moments in biological field, such as animals or plants breeding programs (see Herrend\"{o}rfer and Tuchscherer, 1996) and medical treatments with risk variables as blood pressures and pulses, where low and high values of the variables are of primary concern. Galarza et al. (2022) further computed doubly truncated moments for the selection elliptical class of distributions, and established sufficient and necessary
conditions for the existence of these truncated moments.

From a practical viewpoint, skewed distribution model is more used than non skewed distribution model because of data sets possessing large skewness and/or kurtosis measures (for instant, in economic
and financial data sets). Based on this reason, Roozegar et al. (2020) derived explicit expressions of the first two moments for doubly truncated multivariate normal mean-variance mixture distributions.

Inspired by those works, we derived multivariate doubly truncated first two moments
for GSE distributions, and provided expressions of multivariate doubly truncated expectation and covariance for those distribution. Some important cases of those distributions, including GSN, GSSt, GSLo and GSLa distributions, also were presented. As applications, formulas of MTCE and MTCov for GSE distributions were derived by our results established.

The remainder of the paper proceeds as follows. In Section 2 we introduce preliminaries, including the GSE distributions and notation. Section 3 focuses on multivariate doubly truncated moments for GSE distributions, proposing
formulas for their multivariate doubly truncated first two moments. In Section 4 we show some special cases.
As applications, expressions of multivariate tail conditional risk measures for GSE are given in Section 5.
 Numerical illustration is shown in Section 6. Finally, the paper closes
with the concluding remarks.

\section{Preliminaries}
We start introducing the definition of GSE~distributions as follows.\\
\noindent$\mathbf{2.1.~Generalized ~skew}$-$\mathbf{elliptical ~distributions}$

 Let $\mathbf{X}\sim E_{n}(\boldsymbol{\mu},\boldsymbol{\Sigma},g_{n})$ (if it exists) be an $n$-dimensional elliptical random vector with location vector $\boldsymbol{\mu}$, scale matrix $\mathbf{\Sigma}$ and density generator $g_{n}(u)$, $u\geq0$.  Its probability density function (pdf) takes the form (see, for instance, Landsman and Valdez, 2003)
 \begin{align}\label{(2)}
f_{\boldsymbol{X}}(\boldsymbol{x}):=\frac{c_{n}}{\sqrt{|\boldsymbol{\Sigma}|}}g_{n}\left\{\frac{1}{2}(\boldsymbol{x}-\boldsymbol{\mu})^{T}\mathbf{\Sigma}^{-1}(\boldsymbol{x}-\boldsymbol{\mu})\right\},~\boldsymbol{x}\in\mathbb{R}^{n},
\end{align}
where
\begin{align*}
c_{n}=\frac{\Gamma(n/2)}{(2\pi)^{n/2}}\left[\int_{0}^{\infty}s^{n/2-1}g_{n}(s)\mathrm{d}s\right]^{-1}
\end{align*}
is normalizing constant.
The density generator $g_{n}(u)$ satisfies
\begin{align}\label{(3)}
\int_{0}^{\infty}s^{n/2-1}g_{n}(s)\mathrm{d}s<\infty.
\end{align}
We will call the random
vector $\mathbf{Y}\sim GSE_{n}\left(\boldsymbol{\mu},~\boldsymbol{\Sigma},~g_{n},~H\right)$ is an $n$-dimensional generalized skew-elliptical (GSE) random vector if its pdf (exists) takes the form (see McNeil et al., 2005; Adcock et al., 2019)
\begin{align}\label{(1)}
f_{\boldsymbol{Y}}(\boldsymbol{y})=\frac{2c_{n}}{\sqrt{|\boldsymbol{\Sigma}|}}g_{n}\left\{\frac{1}{2}(\boldsymbol{y}-\boldsymbol{\mu})^{T}\mathbf{\Sigma}^{-1}(\boldsymbol{y}-\boldsymbol{\mu})\right\}H\left(\mathbf{\Sigma}^{-\frac{1}{2}}(\boldsymbol{y}-\boldsymbol{\mu})\right),~\boldsymbol{y}\in\mathbb{R}^{n},
\end{align}
where
$H$ is the skewing function, and it satisfies $H(\boldsymbol{-t})=1-H(\boldsymbol{t})$ and $0\leq H(\boldsymbol{t})\leq1$ for $\boldsymbol{t}\in\mathbb{R}^{n}$.
Note that $H:\mathbb{R}^{n}\rightarrow\mathbb{R}$, we can also define skewing function $J:\mathbb{R}\rightarrow\mathbb{R}$ through $H(\boldsymbol{t})=J(\boldsymbol{\gamma}^T\boldsymbol{t})$ for $\boldsymbol{t}\in\mathbb{R}^{n}$.

In Shushi (2016), it was proved that the characteristic function of any $n\times1$ random vector $\mathbf{Y}\sim GSE_{n}\left(\boldsymbol{\mu},~\boldsymbol{\Sigma},~g_{n},~H\right)$ takes the following form
$$c(\boldsymbol{t})=2\mathrm{e}^{i\boldsymbol{t}\mu}\psi_{n}\left(\frac{1}{2}\boldsymbol{t}^{T}\Sigma\boldsymbol{t}\right)k_{n}(\boldsymbol{t}),\boldsymbol{t}\in \mathbb{R}^{n},$$
where $\psi_{n}(t):[0,\infty)\rightarrow\mathbb{R}$ is the characteristic generator of elliptical random vector $\mathbf{X}$ (see Fang et al., 1990),
$k_{n}$ is a function, and it satisfies $k_{n}(-\boldsymbol{t})=1-k_{n}(\boldsymbol{t})$.

We define cumulative generator $\overline{G}_{n}(u)$ and $\overline{\mathcal{G}}_{n}(u)$ as follows (see Landsman et al., 2018):
\begin{align*}
\overline{G}_{n}(u)=\int_{u}^{\infty}{g}_{n}(v)\mathrm{d}v
\end{align*}
and
\begin{align*}
\overline{\mathcal{G}}_{n}(u)=\int_{u}^{\infty}{G}_{n}(v)\mathrm{d}v,
\end{align*}
and their normalizing constants are, respectively, written as (see Zuo et al., 2021):
\begin{align*}
c_{n}^{\ast}=\frac{\Gamma(n/2)}{(2\pi)^{n/2}}\left[\int_{0}^{\infty}s^{n/2-1}\overline{G}_{n}(s)\mathrm{d}s\right]^{-1}
\end{align*}
and
\begin{align*}
c_{n}^{\ast\ast}=\frac{\Gamma(n/2)}{(2\pi)^{n/2}}\left[\int_{0}^{\infty}s^{n/2-1}\overline{\mathcal{G}}_{n}(s)\mathrm{d}s\right]^{-1}.
\end{align*}
For those density generators, it is necessary to meet the following conditions
\begin{align}\label{(4)}
\int_{0}^{\infty}s^{n/2-1}\overline{G}_{n}(s)\mathrm{d}s<\infty
\end{align}
and
\begin{align}\label{(5)}
\int_{0}^{\infty}s^{n/2-1}\overline{\mathcal{G}}_{n}(s)\mathrm{d}s<\infty.
\end{align}

Now, we define elliptical random vectors $\mathbf{X}^{\ast}\sim E_{n}(\boldsymbol{\mu},~\boldsymbol{\Sigma},~\overline{G}_{n})$ and $\mathbf{X}^{\ast\ast}\sim E_{n}(\boldsymbol{\mu},~\boldsymbol{\Sigma},~\overline{\mathcal{G}}_{n})$. Their form of pdf (if them exist) are, respectively:
\begin{align*}
f_{\boldsymbol{X}^{\ast}}(\boldsymbol{x})=\frac{c_{n}^{\ast}}{\sqrt{|\boldsymbol{\Sigma}|}}\overline{G}_{n}\left\{\frac{1}{2}(\boldsymbol{x}-\boldsymbol{\mu})^{T}\mathbf{\Sigma}^{-1}(\boldsymbol{x}-\boldsymbol{\mu})\right\},~\boldsymbol{x}\in\mathbb{R}^{n}
\end{align*}
\begin{align*}
f_{\boldsymbol{X}^{\ast\ast}}(\boldsymbol{x})=\frac{c_{n}^{\ast\ast}}{\sqrt{|\boldsymbol{\Sigma}|}}\overline{\mathcal{G}}_{n}\left\{\frac{1}{2}(\boldsymbol{x}-\boldsymbol{\mu})^{T}\mathbf{\Sigma}^{-1}(\boldsymbol{x}-\boldsymbol{\mu})\right\},~\boldsymbol{x}\in\mathbb{R}^{n}.
\end{align*}
So
 $\mathbf{Y}^{\ast}\sim GSE_{n}(\boldsymbol{\mu},~\mathbf{\Sigma},~\overline{G}_{n},~H)$ and $\mathbf{Y}^{\ast\ast}\sim GSE_{n}(\boldsymbol{\mu},~\mathbf{\Sigma},~\overline{\mathcal{G}}_{n},~H)$ are corresponding generalized skew-elliptical random vectors.

Next, we introduce some notations and present some concepts that will be used in our proposed
theory.\\
$\mathbf{2.2.~Notation}$

Assume $\mathbf{W}=(W_{1},~W_{2},\cdots,W_{n})^{T}\in\mathbb{R}^{n}$ is an arbitrary random vector with probability density function $f_{\boldsymbol{W}}(\boldsymbol{w})$, then
 for any $\boldsymbol{a}= (a_{1},~a_{2},\cdots, a_{n})^{T}$, $\boldsymbol{b}= (b_{1},~b_{2},\cdots, b_{n})^{T}$
and $\boldsymbol{a} < \boldsymbol{b}$, i.e., $a_{k} < b_{k}$, for $k\in\{1,~2, \cdots , n\}$, we denote the doubly truncated random
vector $\mathbf{W }|(\boldsymbol{a} < \mathbf{W} \leq \boldsymbol{b})$ by $\mathbf{W}_{(\boldsymbol{a},\boldsymbol{b})}$, and $\mathrm{Pr} (\boldsymbol{a} < \mathbf{W} \leq \boldsymbol{b}) = \mathrm{Pr} (a_{1} < W_{1} \leq b_{1}, a_{2} < W_{2} \leq b_{2},\cdots , a_{n} < W_{n} \leq b_{n})$
by $F_{\mathbf{W}} (\boldsymbol{a},\boldsymbol{b}).$

For any $\boldsymbol{v}= (v_{1},~v_{2},\cdots, v_{n})^{T}\in\mathbb{R}^{n}$, writing
$$\boldsymbol{v}_{-k}=(v_{1},~v_{2},\cdots,v_{k-1},~v_{k+1},\cdots,v_{n})^{T},~k\in\{1,~2,\cdots,n\}$$
and
$$\boldsymbol{v}_{-kl}=(v_{1},\cdots,v_{k-1},~v_{k+1},\cdots,v_{l-1},~v_{l+1},\cdots,v_{n})^{T},~k,l\in\{1,~2,\cdots,n\}.$$
 And denoting
$$\mathbf{W}_{v_{k}}=(W_{1},~W_{2},\cdots,W_{k-1},~v_{k},~W_{k+1},\cdots,W_{n})^{T},~k\in\{1,~2,\cdots,n\}$$
and
$$\mathbf{W}_{v_{k},v_{l}}=(W_{1},\cdots,W_{k-1},~v_{k},~W_{k+1},\cdots,W_{l-1},~v_{l},~W_{l+1},\cdots,W_{n})^{T},~k,l\in\{1,~2,\cdots,n\}.$$
In addition,
$$\mathbf{W}_{-v_{k}}=(W_{1},~W_{2},\cdots,W_{k-1},~W_{k+1},\cdots,W_{n})^{T}\in \mathbb{R}^{n-1},~k\in\{1,~2,\cdots,n\}$$
and
$$\mathbf{W}_{-v_{k},v_{l}}=(W_{1},\cdots,W_{k-1},~W_{k+1},\cdots,W_{l-1},~W_{l+1},\cdots,W_{n})^{T}\in \mathbb{R}^{n-2},~k,l\in\{1,~2,\cdots,n\}.$$

To give expression of multivariate doubly truncated moment, we denote doubly truncated expectation $\overline{\mathrm{E}}^{(\boldsymbol{t,v})}_{\mathbf{W}}[h(\boldsymbol{W})]$ of $n$-dimensional random vector $\mathbf{W}$ with pdf $f_{\mathbf{W}}(\boldsymbol{w})$ as
$$\overline{\mathrm{E}}^{(\boldsymbol{t,v})}_{\mathbf{W}}[h(\boldsymbol{W})]=\int_{\boldsymbol{t}}^{\boldsymbol{v}}h(\boldsymbol{w})f_{\mathbf{W}}(\boldsymbol{w})\mathrm{d}\boldsymbol{w},~\boldsymbol{w},\boldsymbol{t},~\boldsymbol{v}\in\mathbb{R}^{n},$$
 where $h:~\mathbb{R}^{n}\rightarrow\mathbb{R}$.\\
 $\mathbf{Remark~1.}$ When $\boldsymbol{v\rightarrow+\infty}$, the doubly truncated expectation will be tail expectation:
 $$\overline{\mathrm{E}}^{\boldsymbol{t}}_{\mathbf{W}}[h(\boldsymbol{W})]=\int_{\boldsymbol{t}}^{\boldsymbol{+\infty}}h(\boldsymbol{w})f_{\mathbf{W}}(\boldsymbol{w})\mathrm{d}\boldsymbol{w},~\boldsymbol{w},\boldsymbol{t}\in\mathbb{R}^{n},$$
 which is defined in Zuo and Yin (2021b).

 The following
notation will be used throughout this paper: $\Phi(\cdot)$, $T(\cdot)$ and $Lo(\cdot)$ denote the cumulative distribution function (cdf) of the univariate
standard normal, student-$t$ (with degrees of freedom $m$) and logistic distributions;  $\phi(\cdot)$, $t(\cdot)$ and $lo(\cdot)$ denote the pdf of the univariate
standard normal, student-$t$ (with degrees of freedom $m$) and logistic distributions, respectively.

 \section{Multivariate doubly truncated moments}

   Let $\mathbf{Y}\sim GSE_{n}\left(\boldsymbol{\mu},~\boldsymbol{\Sigma},~g_{n},~H\right)$ be a random vector with finite fixed vector $\boldsymbol{\mu}=(\mu_1,\cdots,\mu_n)^{T}$, positive defined fixed matrix ${\bf {\Sigma}}= (\sigma_{ij})_{i,j=1}^{n}$ and pdf $f_{\boldsymbol{Y}}(\boldsymbol{y})$.
Let $\mathbf{Z}=\mathbf{\Sigma}^{-\frac{1}{2}}(\mathbf{Y}-\boldsymbol{\mu})\sim GSE_{n}\left(\boldsymbol{0},~\boldsymbol{I_{n}},~g_{n},~H\right).$
Writing
$\boldsymbol{\xi_{v}}=\left(\xi_{\boldsymbol{v},1},~\xi_{\boldsymbol{v},2},\cdots,\xi_{\boldsymbol{v},n}\right)^{T}=\mathbf{\Sigma}^{-\frac{1}{2}}(\boldsymbol{v-\mu}),$
 $\boldsymbol{v}\in\{\boldsymbol{a},~\boldsymbol{b}\}$.
Denoting
$$\mathbf{M^{\ast}}_{-\xi_{\boldsymbol{s}k}}=(M_{1}^{\ast},\cdots,M_{k-1}^{\ast},M_{k+1}^{\ast},\cdots,M_{n}^{\ast})^{T}\in \mathbb{R}^{n-1},$$
 $$\mathbf{M^{\ast\ast}}_{-\xi_{\boldsymbol{s}k}}=(M_{1}^{\ast\ast},\cdots,M_{k-1}^{\ast\ast},M_{k+1}^{\ast\ast},\cdots,M_{n}^{\ast\ast})^{T}\in \mathbb{R}^{n-1}$$
and
$$\mathbf{M^{\ast\ast}}_{-\xi_{\boldsymbol{s}k},\xi_{\boldsymbol{t}l}}=(M_{1}^{\ast\ast},\cdots,M_{k-1}^{\ast\ast},M_{k+1}^{\ast\ast},\cdots,M_{l-1}^{\ast\ast},M_{l+1}^{\ast\ast},\cdots,M_{n}^{\ast\ast})^{T}\in \mathbb{R}^{n-2}.$$
Now, we define $f_{\mathbf{M^{\ast}}_{-\xi_{\boldsymbol{s}k}}}(\boldsymbol{w})$, $f_{\mathbf{M^{\ast\ast}}_{-\xi_{\boldsymbol{s}k}}}(\boldsymbol{v})$ and $f_{\mathbf{M^{\ast\ast}}_{-\xi_{\boldsymbol{s}k},\xi_{\boldsymbol{t}l}}}(\boldsymbol{u})$, the pdf associated with elliptical random vectors $\mathbf{M^{\ast}}_{-\xi_{\boldsymbol{s}k}}$, $\mathbf{M^{\ast\ast}}_{-\xi_{\boldsymbol{s}k}}$ and $\mathbf{M^{\ast\ast}}_{-\xi_{\boldsymbol{s}k},\xi_{\boldsymbol{t}l}}$, respectively:
\begin{align}\label{(6)}
f_{\mathbf{M^{\ast}}_{-\xi_{\boldsymbol{s}k}}}(\boldsymbol{w})=c_{n-1,\xi_{\boldsymbol{s}k}}^{\ast}\overline{G}_{n}\left\{\frac{1}{2}\boldsymbol{w}^{T}\boldsymbol{w}+\frac{1}{2}\xi_{\boldsymbol{s}k}^{2}\right\},~\boldsymbol{w}\in \mathbb{R}^{n-1},
\end{align}
\begin{align}\label{(7)}
f_{\mathbf{M^{\ast\ast}}_{-\xi_{\boldsymbol{s}k}}}(\boldsymbol{v})=c_{n-1,\xi_{\boldsymbol{s}k}}^{\ast\ast}\overline{\mathcal{G}}_{n}\left\{\frac{1}{2}\boldsymbol{v}^{T}\boldsymbol{v}+\frac{1}{2}\xi_{\boldsymbol{s}k}^{2}\right\},~\boldsymbol{v}\in \mathbb{R}^{n-1},
\end{align}
\begin{align}\label{(8)}
f_{\mathbf{M^{\ast\ast}}_{-\xi_{\boldsymbol{s}k},\xi_{\boldsymbol{t}l}}}(\boldsymbol{u})=c_{n-2,\xi_{\boldsymbol{s}k},\xi_{\boldsymbol{t}l}}^{\ast\ast}\overline{\mathcal{G}}_{n}\left\{\frac{1}{2}\boldsymbol{u}^{T}\boldsymbol{u}+\frac{1}{2}\xi_{\boldsymbol{s}k}^{2}+\frac{1}{2}\xi_{\boldsymbol{t}l}^{2}\right\},~\boldsymbol{u}\in \mathbb{R}^{n-2},
\end{align}
where $\boldsymbol{s},\boldsymbol{t}\in\{\boldsymbol{a,~b}\},~k,l\in\{1,~2,\cdots,n\}.$ In addition, $c_{n-1,\xi_{\boldsymbol{s}k}}^{\ast}$, $c_{n-1,\xi_{\boldsymbol{s}k}}^{\ast\ast}$
and
$c_{n-2,\xi_{\boldsymbol{s}k},\xi_{\boldsymbol{t}l}}^{\ast\ast}$ are corresponding normalizing constants of $f_{\mathbf{M^{\ast}}_{-\xi_{\boldsymbol{s}k}}}(\boldsymbol{w})$, $f_{\mathbf{M^{\ast\ast}}_{-\xi_{\boldsymbol{s}k}}}(\boldsymbol{v})$ and $f_{\mathbf{M^{\ast\ast}}_{-\xi_{\boldsymbol{s}k},\xi_{\boldsymbol{t}l}}}(\boldsymbol{u})$, which are written as:
\begin{align}\label{(9)}
c_{n-1,\xi_{\boldsymbol{s}k}}^{\ast}=\frac{\Gamma\left((n-1)/2\right)}{(2\pi)^{(n-1)/2}}\left[\int_{0}^{\infty}x^{(n-3)/2}\overline{G}_{n}\left(\frac{1}{2}\xi_{\boldsymbol{s}k}^{2}+x\right)\mathrm{d}x\right]^{-1},
\end{align}
\begin{align}\label{(10)}
c_{n-1,\xi_{\boldsymbol{s}k}}^{\ast\ast}=\frac{\Gamma\left((n-1)/2\right)}{(2\pi)^{(n-1)/2}}\left[\int_{0}^{\infty}x^{(n-3)/2}\overline{\mathcal{G}}_{n}\left(\frac{1}{2}\xi_{\boldsymbol{s}k}^{2}+x\right)\mathrm{d}x\right]^{-1}
\end{align}
and
\begin{align}\label{(11)}
c_{n-2,\xi_{\boldsymbol{s}k},\xi_{\boldsymbol{t}l}}^{\ast\ast}=\frac{\Gamma\left((n-2)/2\right)}{(2\pi)^{(n-2)/2}}\left[\int_{0}^{\infty}x^{(n-4)/2}\overline{\mathcal{G}}_{n}\left(\frac{1}{2}\xi_{\boldsymbol{s}k}^{2}+\frac{1}{2}\xi_{\boldsymbol{t}l}^{2}+x\right)\mathrm{d}x\right]^{-1}.
\end{align}

Next, we present explicit expressions of multivariate doubly truncated (first two) moments for generalized skew-elliptical distributions.
\begin{theorem}\label{th.1} Let $\mathbf{Y}\sim GSE_{n}(\boldsymbol{\mu},~\mathbf{\Sigma},~g_{n},~H)$ $(n\geq2)$ be as in $(\ref{(1)})$. Suppose that it satisfies conditions $(\ref{(3)})$, $(\ref{(4)})$ and $(\ref{(5)})$. Further, assume $\partial_{i}H$ and $\partial_{ij}H$ exist for $i,j\in\{1,2,\cdots,n\}$. Then
\begin{align}\label{(12)}
\mathrm{(i)}~~\mathrm{E}[\mathbf{Y}|\boldsymbol{a}<\mathbf{Y}\leq \boldsymbol{b}]=\boldsymbol{\mu}+\frac{\mathbf{\Sigma}^{\frac{1}{2}}\boldsymbol{\delta}}{F_{\mathbf{Z}}(\boldsymbol{\xi_{a}},\boldsymbol{\xi_{b}})},
\end{align}
 \begin{align}\label{(13)}
\mathrm{(ii)}~~&\mathrm{E}[\mathbf{Y}\mathbf{Y}^{T}|\boldsymbol{a}<\mathbf{Y}\leq \boldsymbol{b}]=\boldsymbol{\mu}\boldsymbol{\mu}^{T}+\frac{\mathbf{\Sigma}^{\frac{1}{2}}\boldsymbol{\delta}\boldsymbol{\mu}^{T}}{F_{\mathbf{Z}}(\boldsymbol{\xi_{a}},\boldsymbol{\xi_{b}})}+\frac{\boldsymbol{\mu}\boldsymbol{\delta}^{T}\mathbf{\Sigma}^{\frac{1}{2}}}{F_{\mathbf{Z}}(\boldsymbol{\xi_{a}},\boldsymbol{\xi_{b}})}+ \frac{\mathbf{\Sigma}^{\frac{1}{2}}\boldsymbol{\Omega}\mathbf{\Sigma}^{\frac{1}{2}}}{F_{\mathbf{Z}}(\boldsymbol{\xi_{a}},\boldsymbol{\xi_{b}})}   ,
\end{align}
where $\mathbf{\Omega}=(\Omega_{ij})_{i,j=1}^{n}$ is an $n\times n$ symmetric matrix, and $\boldsymbol{\delta}=(\delta_{1},\delta_{2},\cdots,\delta_{n})^{T}$ is an $n\times1$ vector.
Here \begin{align}\label{(14)}
\nonumber\Omega_{ij}=& 2\bigg\{\frac{c_{n}}{c_{n-2,\xi_{\boldsymbol{a}i},\xi_{\boldsymbol{a}j}}^{\ast\ast}}\overline{\mathrm{E}}^{(\boldsymbol{\xi}_{\boldsymbol{a},-ij},\boldsymbol{\xi}_{\boldsymbol{b},-ij})}_{\boldsymbol{M}_{-\xi_{\boldsymbol{a}i},\xi_{\boldsymbol{a}j}}^{\ast\ast}}[H(\boldsymbol{M}_{\xi_{\boldsymbol{a}i},\xi_{\boldsymbol{a}j}}^{\ast\ast})]-\frac{c_{n}}{c_{n-2,\xi_{\boldsymbol{a}i},\xi_{\boldsymbol{b}j}}^{\ast\ast}}\overline{\mathrm{E}}^{(\boldsymbol{\xi}_{\boldsymbol{a},-ij},\boldsymbol{\xi}_{\boldsymbol{b},-ij})}_{\boldsymbol{M}_{-\xi_{\boldsymbol{a}i},\xi_{\boldsymbol{b}j}}^{\ast\ast}}[H(\boldsymbol{M}_{\xi_{\boldsymbol{a}i},\xi_{\boldsymbol{b}j}}^{\ast\ast})]\\
\nonumber&+\frac{c_{n}}{c_{n-1,\xi_{\boldsymbol{a}i}}^{\ast\ast}}\overline{\mathrm{E}}^{(\boldsymbol{\xi}_{\boldsymbol{a},-i},\boldsymbol{\xi}_{\boldsymbol{b},-i})}_{\boldsymbol{M}_{-\xi_{\boldsymbol{a}i}}^{\ast\ast}}[\partial_{j}H(\boldsymbol{M}_{\xi_{\boldsymbol{a}i}}^{\ast\ast})]\\
\nonumber&+\frac{c_{n}}{c_{n-2,\xi_{\boldsymbol{b}i},\xi_{\boldsymbol{b}j}}^{\ast\ast}}\overline{\mathrm{E}}^{(\boldsymbol{\xi}_{\boldsymbol{a},-ij},\boldsymbol{\xi}_{\boldsymbol{b},-ij})}_{\boldsymbol{M}_{-\xi_{\boldsymbol{b}i},\xi_{\boldsymbol{b}j}}^{\ast\ast}}[H(\boldsymbol{M}_{\xi_{\boldsymbol{b}i},\xi_{\boldsymbol{b}j}}^{\ast\ast})]-\frac{c_{n}}{c_{n-2,\xi_{\boldsymbol{b}i},\xi_{\boldsymbol{a}j}}^{\ast\ast}}\overline{\mathrm{E}}^{(\boldsymbol{\xi}_{\boldsymbol{a},-ij},\boldsymbol{\xi}_{\boldsymbol{b},-ij})}_{\boldsymbol{M}_{-\xi_{\boldsymbol{b}i},\xi_{\boldsymbol{a}j}}^{\ast\ast}}[H(\boldsymbol{M}_{\xi_{\boldsymbol{b}i},\xi_{\boldsymbol{a}j}}^{\ast\ast})]\\
\nonumber&-\frac{c_{n}}{c_{n-1,\xi_{\boldsymbol{b}i}}^{\ast\ast}}\overline{\mathrm{E}}^{(\boldsymbol{\xi}_{\boldsymbol{a},-i},\boldsymbol{\xi}_{\boldsymbol{b},-i})}_{\boldsymbol{M}_{-\xi_{\boldsymbol{b}i}}^{\ast\ast}}[\partial_{j}H(\boldsymbol{M}_{\xi_{\boldsymbol{b}i}}^{\ast\ast})]\\
\nonumber&+\frac{c_{n}}{c_{n-1,\xi_{\boldsymbol{a}j}}^{\ast\ast}}\overline{\mathrm{E}}^{(\boldsymbol{\xi}_{\boldsymbol{a},-j},\boldsymbol{\xi}_{\boldsymbol{b},-j})}_{\boldsymbol{M}_{-\xi_{\boldsymbol{a}j}}^{\ast\ast}}[\partial_{i}H(\boldsymbol{M}_{\xi_{\boldsymbol{a}j}}^{\ast\ast})]-\frac{c_{n}}{c_{n-1,\xi_{\boldsymbol{b}j}}^{\ast\ast}}\overline{\mathrm{E}}^{(\boldsymbol{\xi}_{\boldsymbol{a},-j},\boldsymbol{\xi}_{\boldsymbol{b},-j})}_{\boldsymbol{M}_{-\xi_{\boldsymbol{b}j}}^{\ast\ast}}[\partial_{i}H(\boldsymbol{M}_{\xi_{\boldsymbol{b}j}}^{\ast\ast})]\\
&+\frac{c_{n}}{c_{n}^{\ast\ast}}\overline{\mathrm{E}}^{(\boldsymbol{\xi}_{\boldsymbol{a}},\boldsymbol{\xi}_{\boldsymbol{b}})}_{\mathbf{M}^{\ast\ast}}[\partial_{ij}H(\boldsymbol{M}^{\ast\ast})]\bigg\},~i\neq j,
\end{align}
\begin{align}\label{(15)}
\nonumber\Omega_{ii}=& 2\bigg\{\frac{c_{n}}{c_{n-1,\xi_{\boldsymbol{a}i}}^{\ast}}\xi_{\boldsymbol{a},i}\overline{\mathrm{E}}^{(\boldsymbol{\xi}_{\boldsymbol{a},-i},\boldsymbol{\xi}_{\boldsymbol{b},-i})}_{\boldsymbol{M}_{-\xi_{\boldsymbol{a}i}}^{\ast}}[H(\boldsymbol{M}_{\xi_{\boldsymbol{a}i}}^{\ast})]-\frac{c_{n}}{c_{n-1,\xi_{\boldsymbol{b}i}}^{\ast}}\xi_{\boldsymbol{b},i}\overline{\mathrm{E}}^{(\boldsymbol{\xi}_{\boldsymbol{a},-i},\boldsymbol{\xi}_{\boldsymbol{b},-i})}_{\boldsymbol{M}_{-\xi_{\boldsymbol{b}i}}^{\ast}}[H(\boldsymbol{M}_{\xi_{\boldsymbol{b}i}}^{\ast})]\\ \nonumber&+\frac{c_{n}}{c_{n-1,\xi_{\boldsymbol{a}i}}^{\ast\ast}}\overline{\mathrm{E}}^{(\boldsymbol{\xi}_{\boldsymbol{a},-i},\boldsymbol{\xi}_{\boldsymbol{b},-i})}_{\boldsymbol{M}_{-\xi_{\boldsymbol{a}i}}^{\ast\ast}}[\partial_{i}H(\boldsymbol{M}_{\xi_{\boldsymbol{a}i}}^{\ast\ast})]-\frac{c_{n}}{c_{n-1,\xi_{\boldsymbol{b}i}}^{\ast\ast}}\overline{\mathrm{E}}^{(\boldsymbol{\xi}_{\boldsymbol{a},-i},\boldsymbol{\xi}_{\boldsymbol{b},-i})}_{\boldsymbol{M}_{-\xi_{\boldsymbol{b}i}}^{\ast\ast}}[\partial_{i}H(\boldsymbol{M}_{\xi_{\boldsymbol{b}i}}^{\ast\ast})]\\
&+\frac{c_{n}}{c_{n}^{\ast\ast}}\overline{\mathrm{E}}^{(\boldsymbol{\xi}_{\boldsymbol{a}},\boldsymbol{\xi}_{\boldsymbol{b}})}_{\mathbf{M}^{\ast\ast}}[\partial_{ii}H(\boldsymbol{M}^{\ast\ast})]\bigg\}
+\frac{c_{n}}{c_{n}^{\ast}}F_{\mathbf{Z^{\ast}}}(\boldsymbol{\xi}_{\boldsymbol{a}},\boldsymbol{\xi}_{\boldsymbol{b}})
,
\end{align}
\begin{align}\label{(16)}
\delta_{k}=&2\bigg\{\frac{c_{n}}{c_{n-1,\xi_{\boldsymbol{a}k}}^{\ast}}\overline{\mathrm{E}}^{(\boldsymbol{\xi}_{\boldsymbol{a},-k},\boldsymbol{\xi}_{\boldsymbol{b},-k})}_{\boldsymbol{M}_{-\xi_{\boldsymbol{a}k}}^{\ast}}[H(\boldsymbol{M}_{\xi_{\boldsymbol{a}k}}^{\ast})]-\frac{c_{n}}{c_{n-1,\xi_{\boldsymbol{b}k}}^{\ast}}\overline{\mathrm{E}}^{(\boldsymbol{\xi}_{\boldsymbol{a},-k},\boldsymbol{\xi}_{\boldsymbol{b},-k})}_{\boldsymbol{M}_{-\xi_{\boldsymbol{b}k}}^{\ast}}[H(\boldsymbol{M}_{\xi_{\boldsymbol{b}k}}^{\ast})]+\frac{c_{n}}{c_{n}^{\ast}}\overline{\mathrm{E}}^{(\boldsymbol{\xi}_{\boldsymbol{a}},\boldsymbol{\xi}_{\boldsymbol{b}})}_{\mathbf{M}^{\ast}}[\partial_{k}H(\boldsymbol{M}^{\ast})]\bigg\},
\end{align}
 $~i,j,k\in\{1,~2,\cdots,n\},$ $\mathbf{Z}\sim GSE_{n}\left(\boldsymbol{0},~\boldsymbol{I_{n}},~g_{n},~H\right)$, $\mathbf{Z}^{\ast}\sim GSE_{n}\left(\boldsymbol{0},~\boldsymbol{I_{n}},~\overline{G}_{n},~H\right)$, $\mathbf{M}^{\ast}\sim E_{n}\left(\boldsymbol{0},~\boldsymbol{I_{n}},~\overline{G}_{n}\right)$,\\ $\mathbf{M}^{\ast\ast}\sim E_{n}(\boldsymbol{0},~\mathbf{I_{n}},~\overline{\mathcal{G}}_{n})$, and pdfs of $\mathbf{M^{\ast}}_{-\xi_{\boldsymbol{s}k}}$, $\mathbf{M^{\ast\ast}}_{-\xi_{\boldsymbol{s}k}}$ and $\mathbf{M^{\ast\ast}}_{-\xi_{\boldsymbol{s}k},\xi_{\boldsymbol{t}l}}$
   are same as in (\ref{(6)}), (\ref{(7)}) and (\ref{(8)}), respectively. The corresponding normalizing constants $c_{n-1,\xi_{\boldsymbol{s}k}}^{\ast}$, $c_{n-1,\xi_{\boldsymbol{s}k}}^{\ast\ast}$
and
$c_{n-2,\xi_{\boldsymbol{s}k},\xi_{\boldsymbol{t}l}}^{\ast\ast}$ are as in (\ref{(9)}), (\ref{(10)}) and (\ref{(11)}), respectively. In addition, $\partial_{i}H(\boldsymbol{z})=\frac{\partial H(\boldsymbol{z})}{\partial z_{i}}$ and $\partial_{ij}H(\boldsymbol{z})=\frac{\partial^{2} H(\boldsymbol{z})}{\partial z_{i}\partial z_{j}}$.
\end{theorem}
\noindent $\mathbf{Proof}$. (i) Using the transformation $\mathbf{Z}=\mathbf{\Sigma}^{-\frac{1}{2}}(\mathbf{Y}-\boldsymbol{\mu})$ and basic algebraic calculations, we have
\begin{align*}
\mathrm{E}[\mathbf{Y}|\boldsymbol{a}<\mathbf{Y}\leq \boldsymbol{b}]&=\mathrm{E}\left[\left(\mathbf{\Sigma}^{\frac{1}{2}}\mathbf{Z}+\boldsymbol{\mu}\right)|\boldsymbol{\xi_{a}}<\mathbf{Z}\leq \boldsymbol{\xi_{b}}\right]\\
&=\boldsymbol{\mu}+\mathbf{\Sigma}^{\frac{1}{2}}\mathrm{E}\left[\mathbf{Z}|\boldsymbol{\xi_{a}}<\mathbf{Z}\leq\boldsymbol{\xi_{b}}\right],
\end{align*}
where $\boldsymbol{\xi_{v}}=\left(\xi_{\boldsymbol{v},1},~\xi_{\boldsymbol{v},2},\cdots,\xi_{\boldsymbol{v},n}\right)^{T}=\mathbf{\Sigma}^{-\frac{1}{2}}(\boldsymbol{v-\mu}),$
 $\boldsymbol{v}\in\{\boldsymbol{a},~\boldsymbol{b}\}$.\\
By definition of conditional expectation, we obtain
\begin{align*}
\mathrm{E}\left[\mathbf{Z}|\boldsymbol{\xi_{a}}<\mathbf{Z}\leq\boldsymbol{\xi_{b}}\right]&=\frac{1}{F_{\mathbf{Z}}(\boldsymbol{\xi_{a}},\boldsymbol{\xi_{b}})}\int_{\boldsymbol{\xi_{a}}}^{\boldsymbol{\xi_{b}}}\boldsymbol{z}2c_{n}g_{n}\left(\frac{1}{2}\boldsymbol{z}^{T}\boldsymbol{z}\right)H(\boldsymbol{z})\mathrm{d}\boldsymbol{z}\\
&=\frac{\boldsymbol{\delta}}{F_{\mathbf{Z}}(\boldsymbol{\xi_{a}},\boldsymbol{\xi_{b}})},
\end{align*}
where $\boldsymbol{\delta}=(\delta_{1},\delta_{2},\cdots,\delta_{n})^{T}$.\\
 Note that
\begin{align*}
\delta_{k}=&\int_{\boldsymbol{\xi_{a}}}^{\boldsymbol{\xi_{b}}}z_{k}2c_{n}g_{n}\left(\frac{1}{2}\boldsymbol{z}^{T}\boldsymbol{z}\right)H(\boldsymbol{z})\mathrm{d}\boldsymbol{z}\\
=&2c_{n}\int_{\boldsymbol{\xi_{a}},-k}^{\boldsymbol{\xi_{b}},-k}\int_{\xi_{\boldsymbol{a},k}}^{\xi_{\boldsymbol{b},k}}-H(\boldsymbol{z})\partial_{k}\overline{G}_{n}\left(\frac{1}{2}\boldsymbol{z}_{-k}^{T}\boldsymbol{z}_{-k}+\frac{1}{2}z_{k}^{2}\right)\mathrm{d}\boldsymbol{z}_{-k}\\
=&2c_{n}\int_{\boldsymbol{\xi_{a}},-k}^{\boldsymbol{\xi_{b}},-k}\bigg\{[H(\boldsymbol{z}_{\xi_{\boldsymbol{a},k}})]\overline{G}_{n}\left(\frac{1}{2}\boldsymbol{z}_{-k}^{T}\boldsymbol{z}_{-k}+\frac{1}{2}\xi_{\boldsymbol{a},k}^{2}\right)-[H(\boldsymbol{z}_{\xi_{\boldsymbol{b},k}})]\overline{G}_{n}\left(\frac{1}{2}\boldsymbol{z}_{-k}^{T}\boldsymbol{z}_{-k}+\frac{1}{2}\xi_{\boldsymbol{b},k}^{2}\right)\bigg\}\mathrm{d}\boldsymbol{z}_{-k}\\
&+2c_{n}\int_{\boldsymbol{\xi_{a}}}^{\boldsymbol{\xi_{b}}}\partial_{k}H(\boldsymbol{z})\overline{G}_{n}\left(\frac{1}{2}\boldsymbol{z}^{T}\boldsymbol{z}\right)\mathrm{d}\boldsymbol{z}\\
=&2\bigg\{\frac{c_{n}}{c_{n-1,\xi_{\boldsymbol{a}k}}^{\ast}}\overline{\mathrm{E}}^{(\boldsymbol{\xi}_{\boldsymbol{a},-k},\boldsymbol{\xi}_{\boldsymbol{b},-k})}_{\boldsymbol{M}_{-\xi_{\boldsymbol{a}k}}^{\ast}}[H(\boldsymbol{M}_{\xi_{\boldsymbol{a}k}}^{\ast})]-\frac{c_{n}}{c_{n-1,\xi_{\boldsymbol{b}k}}^{\ast}}\overline{\mathrm{E}}^{(\boldsymbol{\xi}_{\boldsymbol{a},-k},\boldsymbol{\xi}_{\boldsymbol{b},-k})}_{\boldsymbol{M}_{-\xi_{\boldsymbol{b}k}}^{\ast}}[H(\boldsymbol{M}_{\xi_{\boldsymbol{b}k}}^{\ast})]+\frac{c_{n}}{c_{n}^{\ast}}\overline{\mathrm{E}}^{(\boldsymbol{\xi}_{\boldsymbol{a}},\boldsymbol{\xi}_{\boldsymbol{b}})}_{\mathbf{M}^{\ast}}[\partial_{k}H(\boldsymbol{M}^{\ast})]\bigg\},
\end{align*}
where $\boldsymbol{z}_{\xi_{\boldsymbol{v},k}}=(z_{1},\cdots,z_{k-1},\xi_{\boldsymbol{v},k},z_{k+1},\cdots,z_{n})^{T}$, $\boldsymbol{v}\in\{\boldsymbol{a,~b}\}$,
and we have used integration by parts in the third equality.
Therefore, we get (\ref{(12)}), as required.\\
(ii) Similarly, using the transformation $\mathbf{Z}=\mathbf{\Sigma}^{-\frac{1}{2}}(\mathbf{Y}-\boldsymbol{\mu})$ and basic algebraic calculations, we have
\begin{align*}
\mathrm{E}[\mathbf{Y}\mathbf{Y}^{T}|\boldsymbol{a}<\mathbf{Y}\leq \boldsymbol{b}]=&\mathrm{E}\left[\left(\mathbf{\Sigma}^{\frac{1}{2}}\mathbf{Z}+\boldsymbol{\mu}\right)\left(\mathbf{\Sigma}^{\frac{1}{2}}\mathbf{Z}+\boldsymbol{\mu}\right)^{T}|\boldsymbol{\xi_{a}}<\mathbf{Z}\leq \boldsymbol{\xi_{b}}\right]\\
=&\mathbf{\Sigma}^{\frac{1}{2}} \mathrm{E}\left[\mathbf{Z}\mathbf{Z}^{T}|\boldsymbol{\xi_{a}}<\mathbf{Z}\leq\boldsymbol{\xi_{b}}\right]   \mathbf{\Sigma}^{\frac{1}{2}}+\mathbf{\Sigma}^{\frac{1}{2}}\mathrm{E}\left[\mathbf{Z}|\boldsymbol{\xi_{a}}<\mathbf{Z}\leq\boldsymbol{\xi_{b}}\right]\boldsymbol{\mu}^{T}\\
&+\boldsymbol{\mu}\mathrm{E}\left[\mathbf{Z}^{T}|\boldsymbol{\xi_{a}}<\mathbf{Z}\leq\boldsymbol{\xi_{b}}\right]\mathbf{\Sigma}^{\frac{1}{2}}+\boldsymbol{\mu}\boldsymbol{\mu}^{T}.
\end{align*}
For $\mathrm{E}\left[\mathbf{Z}\mathbf{Z}^{T}|\boldsymbol{\xi_{a}}<\mathbf{Z}\leq\boldsymbol{\xi_{b}}\right]$, by the definition of conditional expectation we obtain
\begin{align*}
\mathrm{E}\left[\mathbf{Z}\mathbf{Z}^{T}|\boldsymbol{\xi_{a}}<\mathbf{Z}\leq\boldsymbol{\xi_{b}}\right]&=\frac{1}{F_{\mathbf{Z}}(\boldsymbol{\xi_{a}},\boldsymbol{\xi_{b}})}\int_{\boldsymbol{\xi_{a}}}^{\boldsymbol{\xi_{b}}}\boldsymbol{z}\boldsymbol{z}^{T}2c_{n}g_{n}\left(\frac{1}{2}\boldsymbol{z}^{T}\boldsymbol{z}\right)H(\boldsymbol{z})\mathrm{d}\boldsymbol{z}\\
&=\frac{\boldsymbol{\Omega}}{F_{\mathbf{Z}}(\boldsymbol{\xi_{a}},\boldsymbol{\xi_{b}})},
\end{align*}
where $\mathbf{\Omega}=(\Omega_{ij})_{i,j=1}^{n}$.\\
Note that, for $i\neq j$,
\begin{align*}
\Omega_{ij}=&\int_{\boldsymbol{\xi_{a}}}^{\boldsymbol{\xi_{b}}}z_{i}z_{j}2c_{n}g_{n}\left(\frac{1}{2}\boldsymbol{z}^{T}\boldsymbol{z}\right)H(\boldsymbol{z})\mathrm{d}\boldsymbol{z}\\
=&2c_{n}\int_{\boldsymbol{\xi_{a}},-i}^{\boldsymbol{\xi_{b}},-i}\int_{\xi_{\boldsymbol{a}},i}^{\xi_{\boldsymbol{b}},i}-z_{j}H(\boldsymbol{z})\partial_{i}\overline{G}_{n}\left(\frac{1}{2}\boldsymbol{z}_{-i}^{T}\boldsymbol{z}_{-i}+\frac{1}{2}z_{i}^{2}\right)\mathrm{d}\boldsymbol{z}_{-i}\\
=&2c_{n}\int_{\boldsymbol{\xi_{a}},-i}^{\boldsymbol{\xi_{b}},-i}\bigg[z_{j}H(\boldsymbol{z}_{\xi_{\boldsymbol{a},i}})\overline{G}_{n}\left(\frac{1}{2}\boldsymbol{z}_{-i}^{T}\boldsymbol{z}_{-i}+\frac{1}{2}\xi_{\boldsymbol{a},i}^{2}\right)- z_{j}H(\boldsymbol{z}_{\xi_{\boldsymbol{b},i}})\overline{G}_{n}\left(\frac{1}{2}\boldsymbol{z}_{-i}^{T}\boldsymbol{z}_{-i}+\frac{1}{2}\xi_{\boldsymbol{b},i}^{2}\right)  \bigg]\mathrm{d}\boldsymbol{z}_{-i}\\
&+2c_{n}\int_{\boldsymbol{\xi_{a}}}^{\boldsymbol{\xi_{b}}}z_{j}\partial_{i}H(\boldsymbol{z})\overline{G}_{n}\left(\frac{1}{2}\boldsymbol{z}^{T}\boldsymbol{z}\right)\mathrm{d}\boldsymbol{z}
\end{align*}
\begin{align*}
=&2c_{n}\int_{\boldsymbol{\xi_{a}},-ij}^{\boldsymbol{\xi_{b}},-ij}\int_{\xi_{\boldsymbol{a}},j}^{\xi_{\boldsymbol{b}},j}\bigg[-H(\boldsymbol{z}_{\xi_{\boldsymbol{a},i}})\partial_{j}\overline{\mathcal{G}}_{n}\left(\frac{1}{2}\boldsymbol{z}_{-i}^{T}\boldsymbol{z}_{-i}+\frac{1}{2}\xi_{\boldsymbol{a},i}^{2}\right)+H(\boldsymbol{z}_{\xi_{\boldsymbol{b},i}})\partial_{j}\overline{\mathcal{G}}_{n}\left(\frac{1}{2}\boldsymbol{z}_{-i}^{T}\boldsymbol{z}_{-i}+\frac{1}{2}\xi_{\boldsymbol{b},i}^{2}\right)\bigg]\mathrm{d}\boldsymbol{z}_{-ij}\\
&+2c_{n}\int_{\boldsymbol{\xi_{a}},-j}^{\boldsymbol{\xi_{b}},-j}\int_{\xi_{\boldsymbol{a}},j}^{\xi_{\boldsymbol{b}},j}-\partial_{i}H(\boldsymbol{z})\partial_{j}\overline{\mathcal{G}}_{n}\left(\frac{1}{2}\boldsymbol{z}_{-j}^{T}\boldsymbol{z}_{-j}+\frac{1}{2}z_{j}^{2}\right)\mathrm{d}\boldsymbol{z}_{-j}\\
=&2c_{n}\int_{\boldsymbol{\xi_{a}},-ij}^{\boldsymbol{\xi_{b}},-ij}\bigg[H(\boldsymbol{z}_{\xi_{\boldsymbol{a}i},\xi_{\boldsymbol{a}j}})\overline{\mathcal{G}}_{n}\left(\frac{1}{2}\boldsymbol{z}_{-ij}^{T}\boldsymbol{z}_{-ij}+\frac{1}{2}\xi_{\boldsymbol{a},i}^{2}+\frac{1}{2}\xi_{\boldsymbol{a},j}^{2}\right)-H(\boldsymbol{z}_{\xi_{\boldsymbol{a}i},\xi_{\boldsymbol{b}j}})\overline{\mathcal{G}}_{n}\left(\frac{1}{2}\boldsymbol{z}_{-ij}^{T}\boldsymbol{z}_{-ij}+\frac{1}{2}\xi_{\boldsymbol{a},i}^{2}+\frac{1}{2}\xi_{\boldsymbol{b},j}^{2}\right)\\
&+\int_{\xi_{\boldsymbol{a}},j}^{\xi_{\boldsymbol{b}},j}\partial_{j}H(\boldsymbol{z}_{\xi_{\boldsymbol{a},i}})\overline{\mathcal{G}}_{n}\left(\frac{1}{2}\boldsymbol{z}_{-i}^{T}\boldsymbol{z}_{-i}+\frac{1}{2}\xi_{\boldsymbol{a},i}^{2}\right)\mathrm{d}z_{j}\\
&+H(\boldsymbol{z}_{\xi_{\boldsymbol{b}i},\xi_{\boldsymbol{b}j}})\overline{\mathcal{G}}_{n}\left(\frac{1}{2}\boldsymbol{z}_{-ij}^{T}\boldsymbol{z}_{-ij}+\frac{1}{2}\xi_{\boldsymbol{b},i}^{2}+\frac{1}{2}\xi_{\boldsymbol{b}j}^{2}\right)-H(\boldsymbol{z}_{\xi_{\boldsymbol{b}i},\xi_{\boldsymbol{a}j}})\overline{\mathcal{G}}_{n}\left(\frac{1}{2}\boldsymbol{z}_{-ij}^{T}\boldsymbol{z}_{-ij}+\frac{1}{2}\xi_{\boldsymbol{b},i}^{2}+\frac{1}{2}\xi_{\boldsymbol{a}j}^{2}\right)\\
&-\int_{\xi_{\boldsymbol{a}},j}^{\xi_{\boldsymbol{b}},j}\partial_{j}H(\boldsymbol{z}_{\xi_{\boldsymbol{b},i}})\overline{\mathcal{G}}_{n}\left(\frac{1}{2}\boldsymbol{z}_{-i}^{T}\boldsymbol{z}_{-i}+\frac{1}{2}\xi_{\boldsymbol{b},i}^{2}\right)\mathrm{d}z_{j}\bigg]\mathrm{d}\boldsymbol{z}_{-ij}\\
&+2c_{n}\int_{\boldsymbol{\xi_{a}},-j}^{\boldsymbol{\xi_{b}},-j}\bigg[\partial_{i}H(\boldsymbol{z}_{\xi_{\boldsymbol{a},j}})\overline{\mathcal{G}}_{n}\left(\frac{1}{2}\boldsymbol{z}_{-j}^{T}\boldsymbol{z}_{-j}+\frac{1}{2}\xi_{\boldsymbol{a},j}^{2}\right)-\partial_{i}H(\boldsymbol{z}_{\xi_{\boldsymbol{b},j}})\overline{\mathcal{G}}_{n}\left(\frac{1}{2}\boldsymbol{z}_{-j}^{T}\boldsymbol{z}_{-j}+\frac{1}{2}\xi_{\boldsymbol{b},j}^{2}\right)\\
&+\int_{\xi_{\boldsymbol{a}},j}^{\xi_{\boldsymbol{b}},j}\partial_{ij}H(\boldsymbol{z})\overline{\mathcal{G}}_{n}\left(\frac{1}{2}\boldsymbol{z}^{T}\boldsymbol{z}\right)\mathrm{d}z_{j}
\bigg]\mathrm{d}\boldsymbol{z}_{-j}\\
=&2\bigg\{\frac{c_{n}}{c_{n-2,\xi_{\boldsymbol{a}i},\xi_{\boldsymbol{a}j}}^{\ast\ast}}\overline{\mathrm{E}}^{(\boldsymbol{\xi}_{\boldsymbol{a},-ij},\boldsymbol{\xi}_{\boldsymbol{b},-ij})}_{\boldsymbol{M}_{-\xi_{\boldsymbol{a}i},\xi_{\boldsymbol{a}j}}^{\ast\ast}}[H(\boldsymbol{M}_{\xi_{\boldsymbol{a}i},\xi_{\boldsymbol{a}j}}^{\ast\ast})]-\frac{c_{n}}{c_{n-2,\xi_{\boldsymbol{a}i},\xi_{\boldsymbol{b}j}}^{\ast\ast}}\overline{\mathrm{E}}^{(\boldsymbol{\xi}_{\boldsymbol{a},-ij},\boldsymbol{\xi}_{\boldsymbol{b},-ij})}_{\boldsymbol{M}_{-\xi_{\boldsymbol{a}i},\xi_{\boldsymbol{b}j}}^{\ast\ast}}[H(\boldsymbol{M}_{\xi_{\boldsymbol{a}i},\xi_{\boldsymbol{b}j}}^{\ast\ast})]\\
\nonumber&+\frac{c_{n}}{c_{n-1,\xi_{\boldsymbol{a}i}}^{\ast\ast}}\overline{\mathrm{E}}^{(\boldsymbol{\xi}_{\boldsymbol{a},-i},\boldsymbol{\xi}_{\boldsymbol{b},-i})}_{\boldsymbol{M}_{-\xi_{\boldsymbol{a}i}}^{\ast\ast}}[\partial_{j}H(\boldsymbol{M}_{\xi_{\boldsymbol{a}i}}^{\ast\ast})]\\
\nonumber&+\frac{c_{n}}{c_{n-2,\xi_{\boldsymbol{b}i},\xi_{\boldsymbol{b}j}}^{\ast\ast}}\overline{\mathrm{E}}^{(\boldsymbol{\xi}_{\boldsymbol{a},-ij},\boldsymbol{\xi}_{\boldsymbol{b},-ij})}_{\boldsymbol{M}_{-\xi_{\boldsymbol{b}i},\xi_{\boldsymbol{b}j}}^{\ast\ast}}[H(\boldsymbol{M}_{\xi_{\boldsymbol{b}i},\xi_{\boldsymbol{b}j}}^{\ast\ast})]-\frac{c_{n}}{c_{n-2,\xi_{\boldsymbol{b}i},\xi_{\boldsymbol{a}j}}^{\ast\ast}}\overline{\mathrm{E}}^{(\boldsymbol{\xi}_{\boldsymbol{a},-ij},\boldsymbol{\xi}_{\boldsymbol{b},-ij})}_{\boldsymbol{M}_{-\xi_{\boldsymbol{b}i},\xi_{\boldsymbol{a}j}}^{\ast\ast}}[H(\boldsymbol{M}_{\xi_{\boldsymbol{b}i},\xi_{\boldsymbol{a}j}}^{\ast\ast})]\\
\nonumber&-\frac{c_{n}}{c_{n-1,\xi_{\boldsymbol{b}i}}^{\ast\ast}}\overline{\mathrm{E}}^{(\boldsymbol{\xi}_{\boldsymbol{a},-i},\boldsymbol{\xi}_{\boldsymbol{b},-i})}_{\boldsymbol{M}_{-\xi_{\boldsymbol{b}i}}^{\ast\ast}}[\partial_{j}H(\boldsymbol{M}_{\xi_{\boldsymbol{b}i}}^{\ast\ast})]\\
\nonumber&+\frac{c_{n}}{c_{n-1,\xi_{\boldsymbol{a}j}}^{\ast\ast}}\overline{\mathrm{E}}^{(\boldsymbol{\xi}_{\boldsymbol{a},-j},\boldsymbol{\xi}_{\boldsymbol{b},-j})}_{\boldsymbol{M}_{-\xi_{\boldsymbol{a}j}}^{\ast\ast}}[\partial_{i}H(\boldsymbol{M}_{\xi_{\boldsymbol{a}j}}^{\ast\ast})]-\frac{c_{n}}{c_{n-1,\xi_{\boldsymbol{b}j}}^{\ast\ast}}\overline{\mathrm{E}}^{(\boldsymbol{\xi}_{\boldsymbol{a},-j},\boldsymbol{\xi}_{\boldsymbol{b},-j})}_{\boldsymbol{M}_{-\xi_{\boldsymbol{b}j}}^{\ast\ast}}[\partial_{i}H(\boldsymbol{M}_{\xi_{\boldsymbol{b}j}}^{\ast\ast})]\\
&+\frac{c_{n}}{c_{n}^{\ast\ast}}\overline{\mathrm{E}}^{(\boldsymbol{\xi}_{\boldsymbol{a}},\boldsymbol{\xi}_{\boldsymbol{b}})}_{\mathbf{M}^{\ast\ast}}[\partial_{ij}H(\boldsymbol{M}^{\ast\ast})]\bigg\},
\end{align*}
where $$\boldsymbol{z}_{\xi_{\boldsymbol{v}i},\xi_{\boldsymbol{v}j}}=(z_{1},\cdots,z_{i-1},\xi_{\boldsymbol{v},i},z_{i+1},\cdots,z_{j-1},\xi_{\boldsymbol{v},j},z_{j+1},\cdots,z_{n})^{T},~\boldsymbol{v}\in\{\boldsymbol{a,~b}\},$$
and we have used integration by parts in the third and fifth equalities.\\
While
\begin{align*}
\Omega_{ii}=&\int_{\boldsymbol{\xi_{a}}}^{\boldsymbol{\xi_{b}}}z_{i}^{2}2c_{n}g_{n}\left(\frac{1}{2}\boldsymbol{z}^{T}\boldsymbol{z}\right)H(\boldsymbol{z})\mathrm{d}\boldsymbol{z}\\
=&2c_{n}\int_{\boldsymbol{\xi}_{\boldsymbol{a},-i}}^{\boldsymbol{\xi}_{\boldsymbol{b},-i}}\int_{\xi_{\boldsymbol{a},i}}^{\xi_{\boldsymbol{b},i}}-z_{i}H(\boldsymbol{z})\partial_{i}\overline{G}_{n}\left(\frac{1}{2}\boldsymbol{z}_{-i}^{T}\boldsymbol{z}_{-i}+\frac{1}{2}z_{i}^{2}\right)\mathrm{d}\boldsymbol{z}_{-i}\\
=&2c_{n}\int_{\boldsymbol{\xi}_{\boldsymbol{a},-i}}^{\boldsymbol{\xi}_{\boldsymbol{b},-i}}\bigg\{\xi_{\boldsymbol{a},i}H(\boldsymbol{z}_{\xi_{\boldsymbol{a},i}})\overline{G}_{n}\left(\frac{1}{2}\boldsymbol{z}_{-i}^{T}\boldsymbol{z}_{-i}+\frac{1}{2}\xi_{\boldsymbol{a},i}^{2}\right)-\xi_{\boldsymbol{b},i}H(\boldsymbol{z}_{\xi_{\boldsymbol{b},i}})\overline{G}_{n}\left(\frac{1}{2}\boldsymbol{z}_{-i}^{T}\boldsymbol{z}_{-i}+\frac{1}{2}\xi_{\boldsymbol{b},i}^{2}\right)\\
&+\int_{\xi_{\boldsymbol{a},i}}^{\xi_{\boldsymbol{b},i}}-\partial_{i}H(\boldsymbol{z})\partial_{i}\overline{\mathcal{G}}_{n}\left(\frac{1}{2}\boldsymbol{z}_{-i}^{T}\boldsymbol{z}_{-i}+\frac{1}{2}z_{i}^{2}\right)
\bigg\}\mathrm{d}\boldsymbol{z}_{-i}
+2c_{n}\int_{\boldsymbol{\xi}_{\boldsymbol{a}}}^{\boldsymbol{\xi}_{\boldsymbol{b}}}H(\boldsymbol{z})\overline{G}_{n}\left(\frac{1}{2}\boldsymbol{z}^{T}\boldsymbol{z}\right)\mathrm{d}\boldsymbol{z}
\end{align*}
\begin{align*}
=&2c_{n}\int_{\boldsymbol{\xi}_{\boldsymbol{a},-i}}^{\boldsymbol{\xi}_{\boldsymbol{b},-i}}\bigg\{\xi_{\boldsymbol{a},i}H(\boldsymbol{z}_{\xi_{\boldsymbol{a},i}})\overline{G}_{n}\left(\frac{1}{2}\boldsymbol{z}_{-i}^{T}\boldsymbol{z}_{-i}+\frac{1}{2}\xi_{\boldsymbol{a},i}^{2}\right)-\xi_{\boldsymbol{b},i}H(\boldsymbol{z}_{\xi_{\boldsymbol{b},i}})\overline{G}_{n}\left(\frac{1}{2}\boldsymbol{z}_{-i}^{T}\boldsymbol{z}_{-i}+\frac{1}{2}\xi_{\boldsymbol{b},i}^{2}\right)\\
&+\partial_{i}H(\boldsymbol{z}_{\xi_{\boldsymbol{a},i}})\overline{\mathcal{G}}_{n}\left(\frac{1}{2}\boldsymbol{z}_{-i}^{T}\boldsymbol{z}_{-i}+\frac{1}{2}\xi_{\boldsymbol{a},i}^{2}\right)
-\partial_{i}H(\boldsymbol{z}_{\xi_{\boldsymbol{b},i}})\overline{\mathcal{G}}_{n}\left(\frac{1}{2}\boldsymbol{z}_{-i}^{T}\boldsymbol{z}_{-i}+\frac{1}{2}\xi_{\boldsymbol{b},i}^{2}\right)\bigg\}\mathrm{d}\boldsymbol{z}_{-i}\\
&+2c_{n}\int_{\boldsymbol{\xi}_{\boldsymbol{a}}}^{\boldsymbol{\xi}_{\boldsymbol{b}}}\partial_{ii}H(\boldsymbol{z})\overline{\mathcal{G}}_{n}\left(\frac{1}{2}\boldsymbol{z}^{T}\boldsymbol{z}\right)\mathrm{d}\boldsymbol{z}
+2c_{n}\int_{\boldsymbol{\xi}_{\boldsymbol{a}}}^{\boldsymbol{\xi}_{\boldsymbol{b}}}H(\boldsymbol{z})\overline{G}_{n}\left(\frac{1}{2}\boldsymbol{z}^{T}\boldsymbol{z}\right)\mathrm{d}\boldsymbol{z}
\\
=&2\bigg\{\frac{c_{n}}{c_{n-1,\xi_{\boldsymbol{a}i}}^{\ast}}\xi_{\boldsymbol{a},i}\overline{\mathrm{E}}^{(\boldsymbol{\xi}_{\boldsymbol{a},-i},\boldsymbol{\xi}_{\boldsymbol{b},-i})}_{\boldsymbol{M}_{-\xi_{\boldsymbol{a}i}}^{\ast}}[H(\boldsymbol{M}_{\xi_{\boldsymbol{a}i}}^{\ast})]-\frac{c_{n}}{c_{n-1,\xi_{\boldsymbol{b}i}}^{\ast}}\xi_{\boldsymbol{b},i}\overline{\mathrm{E}}^{(\boldsymbol{\xi}_{\boldsymbol{a},-i},\boldsymbol{\xi}_{\boldsymbol{b},-i})}_{\boldsymbol{M}_{-\xi_{\boldsymbol{b}i}}^{\ast}}[H(\boldsymbol{M}_{\xi_{\boldsymbol{b}i}}^{\ast})]\\ \nonumber&+\frac{c_{n}}{c_{n-1,\xi_{\boldsymbol{a}i}}^{\ast\ast}}\overline{\mathrm{E}}^{(\boldsymbol{\xi}_{\boldsymbol{a},-i},\boldsymbol{\xi}_{\boldsymbol{b},-i})}_{\boldsymbol{M}_{-\xi_{\boldsymbol{a}i}}^{\ast\ast}}[\partial_{i}H(\boldsymbol{M}_{\xi_{\boldsymbol{a}i}}^{\ast\ast})]-\frac{c_{n}}{c_{n-1,\xi_{\boldsymbol{b}i}}^{\ast\ast}}\overline{\mathrm{E}}^{(\boldsymbol{\xi}_{\boldsymbol{a},-i},\boldsymbol{\xi}_{\boldsymbol{b},-i})}_{\boldsymbol{M}_{-\xi_{\boldsymbol{b}i}}^{\ast\ast}}[\partial_{i}H(\boldsymbol{M}_{\xi_{\boldsymbol{b}i}}^{\ast\ast})]\\
&+\frac{c_{n}}{c_{n}^{\ast\ast}}\overline{\mathrm{E}}^{(\boldsymbol{\xi}_{\boldsymbol{a}},\boldsymbol{\xi}_{\boldsymbol{b}})}_{\mathbf{M}^{\ast\ast}}[\partial_{ii}H(\boldsymbol{M}^{\ast\ast})]\bigg\}
+\frac{c_{n}}{c_{n}^{\ast}}F_{\mathbf{Z^{\ast}}}(\boldsymbol{\xi}_{\boldsymbol{a}},\boldsymbol{\xi}_{\boldsymbol{b}})
,
\end{align*}
where
 we have used integration by parts in the third and fourth equalities.\\
As for $\mathrm{E}\left[\mathbf{Z}|\boldsymbol{\xi_{a}}<\mathbf{Z}\leq\boldsymbol{\xi_{b}}\right]$, using (i) we directly obtain
\begin{align*}
\mathrm{E}\left[\mathbf{Z}|\boldsymbol{\xi_{a}}<\mathbf{Z}\leq\boldsymbol{\xi_{b}}\right]=\frac{\boldsymbol{\delta}}{F_{\mathbf{Z}}(\boldsymbol{\xi_{a}},\boldsymbol{\xi_{b}})}.
\end{align*}
Consequently, we obtain $(\ref{(13)})$, ending the proof of (ii).

Zuo and Yin (2021c) introduced multivariate doubly truncated expectation (MDTE) and covariance (MDTCov) for any an $n\times1$ vector $\mathbf{X}$ as follows, respectively:
\begin{align*}
\mathrm{MDTE}_{(\boldsymbol{a},\boldsymbol{b})}(\mathbf{X})&=\mathrm{E}\left[\mathbf{X}|\boldsymbol{a}<\mathbf{X}\leq\boldsymbol{b}\right]
\end{align*}
and
\begin{align*}
&\mathrm{MDTCov}_{(\boldsymbol{a},\boldsymbol{b})}(\mathbf{X})=\mathrm{E}\left[(\mathbf{X}-\mathrm{MDTE}_{(\boldsymbol{a},\boldsymbol{b})}(\mathbf{X}))(\mathbf{X}-\mathrm{MDTE}_{(\boldsymbol{a},\boldsymbol{b})}(\mathbf{X}))^{T}|\boldsymbol{a}<\mathbf{X}\leq\boldsymbol{b}\right].
\end{align*}

Now, we can give following proposition of MDTE and MDTCov for GSE distributions.
\begin{proposition}
 Under the conditions of Theorem 1, we have
  \begin{align}\label{(17)}
\mathrm{(I)}~\mathrm{MDTE}_{(\boldsymbol{a},\boldsymbol{b})}(\mathbf{Y})&=\boldsymbol{\mu}+\frac{\mathbf{\Sigma}^{\frac{1}{2}}\boldsymbol{\delta}}{F_{\mathbf{Z}}(\boldsymbol{\xi_{a}},\boldsymbol{\xi_{b}})},
\end{align}
\begin{align}\label{(18)}
\mathrm{(II)}~&\mathrm{MDTCov}_{(\boldsymbol{a},\boldsymbol{b})}(\mathbf{Y})=\mathbf{\Sigma}^{\frac{1}{2}}\left[ \frac{\boldsymbol{\Omega}}{F_{\mathbf{Z}}(\boldsymbol{\xi_{a}},\boldsymbol{\xi_{b}})}-\frac{\boldsymbol{\delta}\boldsymbol{\delta}^{T}}{F_{\mathbf{Z}}^{2}(\boldsymbol{\xi_{a}},\boldsymbol{\xi_{b}})} \right]  \mathbf{\Sigma}^{\frac{1}{2}},
\end{align}
where $\boldsymbol{\Omega}$ and $\boldsymbol{\delta}$ are the same as those in Theorem 1.
\end{proposition}
\noindent $\mathbf{Proof}$. (I) Using (i) of Theorem 1, we immediate obtain (\ref{(17)}).\\ (II) By the definition of MDTCov, we have
\begin{align*}
\mathrm{MDTCov}_{(\boldsymbol{a},\boldsymbol{b})}(\mathbf{Y})&=\mathrm{E}\left[(\mathbf{Y}-\mathrm{MDTE}_{(\boldsymbol{a},\boldsymbol{b})}(\mathbf{Y}))(\mathbf{Y}-\mathrm{MDTE}_{(\boldsymbol{a},\boldsymbol{b})}(\mathbf{Y}))^{T}|\boldsymbol{a}<\mathbf{Y}\leq\boldsymbol{b}\right]\\
&=\mathrm{E}\left[\mathbf{Y}\mathbf{Y}^{T}|\boldsymbol{a}<\mathbf{Y}\leq\boldsymbol{b}\right]-\mathrm{MDTE}_{(\boldsymbol{a},\boldsymbol{b})}(\mathbf{Y})\mathrm{MDTE}_{(\boldsymbol{a},\boldsymbol{b})}^{T}(\mathbf{Y}).
\end{align*}
Applying (i) and (ii) of Theorem 1, and using basic algebraic calculations we instantly obtain (\ref{(18)}).\\
$\mathbf{Remark~2.}$ When $H(\cdot)=\frac{1}{2}$, we note that above results coincide with the results of Theorem 1 and Theorem 2 in  Zuo and Yin (2021c).\\
$\mathbf{Remark~3.}$ When $\boldsymbol{a\rightarrow-\infty}$ and $\boldsymbol{b\rightarrow+\infty}$, MDTE and MDTCov are reduced to expectation and covariance, respectively.
\section{Special cases}
This section focuses $\mathrm{E}[\mathbf{Y}|\boldsymbol{a}<\mathbf{Y}\leq \boldsymbol{b}]$, $\mathrm{E}[\mathbf{Y}\mathbf{Y}^{T}|\boldsymbol{a}<\mathbf{Y}\leq \boldsymbol{b}]$ and $\mathrm{MDTCov}_{(\boldsymbol{a},\boldsymbol{b})}(\mathbf{Y})$ for important cases of the GSE distributions, including GSN, GSSt, GSLo and GSLa distributions. Their forms of $\mathrm{E}[\mathbf{Y}|\boldsymbol{a}<\mathbf{Y}\leq \boldsymbol{b}]$, $\mathrm{E}[\mathbf{Y}\mathbf{Y}^{T}|\boldsymbol{a}<\mathbf{Y}\leq \boldsymbol{b}]$ and $\mathrm{MDTCov}_{(\boldsymbol{a},\boldsymbol{b})}(\mathbf{Y})$ are given in (\ref{(12)}), (\ref{(13)}) and (\ref{(18)}), respectively, so that we merely present $\delta_{k}$, $k\in\{1,2,\cdots,n\}$, $\Omega_{ij},~i=j$ and $i\neq j$.\\
\noindent$\mathbf{Corollary~1}$ (GSN distribution). Let $\mathbf{Y}\sim GSN_{n}(\boldsymbol{\mu},~\mathbf{\Sigma},~\boldsymbol{\gamma},~J)$. In this case,\\ $\overline{\mathcal{G}}_{n}(u)=\overline{G}_{n}(u)=g_{n}(u)=\exp(-u)$, $c_{n}^{\ast\ast}=c_{n}^{\ast}=c_{n}=(2\pi)^{-\frac{n}{2}}$ and $H\left(\mathbf{\Sigma}^{-\frac{1}{2}}(\boldsymbol{y}-\boldsymbol{\mu})\right)=J\left(\boldsymbol{\gamma}^{T}\mathbf{\Sigma}^{-\frac{1}{2}}(\boldsymbol{y}-\boldsymbol{\mu})\right),$ where $\boldsymbol{\gamma}=(\gamma_{1},~\gamma_{2},~\cdots,~\gamma_{n})^{T}\in\mathbb{R}^{n}$ and $J: \mathbb{R}\rightarrow\mathbb{R}$.
Since
 \begin{align*}
f_{\mathbf{M^{\ast}}_{-\xi_{\boldsymbol{s}k}}}(\boldsymbol{w})=c_{n-1,\xi_{\boldsymbol{s}k}}^{\ast}\exp\left\{-\frac{1}{2}\boldsymbol{w}^{T}\boldsymbol{w}-\frac{1}{2}\xi_{\boldsymbol{s}k}^{2}\right\}=\phi_{n-1}(\boldsymbol{w}),~\boldsymbol{w}\in \mathbb{R}^{n-1},
\end{align*}
\begin{align*}
f_{\mathbf{M^{\ast\ast}}_{-\xi_{\boldsymbol{s}k}}}(\boldsymbol{v})=c_{n-1,\xi_{\boldsymbol{s}k}}^{\ast\ast}\exp\left\{-\frac{1}{2}\boldsymbol{v}^{T}\boldsymbol{v}-\frac{1}{2}\xi_{\boldsymbol{s}k}^{2}\right\}=\phi_{n-1}(\boldsymbol{v}),~\boldsymbol{v}\in \mathbb{R}^{n-1},
\end{align*}
\begin{align*}
f_{\mathbf{M^{\ast\ast}}_{-\xi_{\boldsymbol{s}k},\xi_{\boldsymbol{t}l}}}(\boldsymbol{u})=c_{n-2,\xi_{\boldsymbol{s}k},\xi_{\boldsymbol{t}l}}^{\ast\ast}\exp\left\{-\frac{1}{2}\boldsymbol{u}^{T}\boldsymbol{u}-\frac{1}{2}\xi_{\boldsymbol{s}k}^{2}-\frac{1}{2}\xi_{\boldsymbol{t}l}^{2}\right\}=\phi_{n-2}(\boldsymbol{u}),~\boldsymbol{u}\in \mathbb{R}^{n-2},
\end{align*}
$\phi_{k}(\cdot)$ denotes the pdf of $N_{k}(\boldsymbol{0}, \Sigma)$ (the
$k$-variate normal distribution with mean vector $\boldsymbol{0}$ and covariance matrix $\Sigma$).
So
$c_{n-1,\xi_{\boldsymbol{s}k}}^{\ast}=c_{n-1,\xi_{\boldsymbol{s}k}}^{\ast\ast}=\frac{(2\pi)^{-\frac{n}{2}}}{\phi(\xi_{\boldsymbol{s}k})}$ and
$c_{n-2,\xi_{\boldsymbol{s}k},\xi_{\boldsymbol{t}l}}^{\ast\ast}=\frac{(2\pi)^{-\frac{n}{2}}}{\phi(\xi_{\boldsymbol{s}k})\phi(\xi_{\boldsymbol{t}l})}$.
Thus,
\begin{align*}
\nonumber\Omega_{ij}=& 2\bigg\{\phi(\xi_{\boldsymbol{a},i})\phi(\xi_{\boldsymbol{a},j})\overline{\mathrm{E}}^{(\boldsymbol{\xi}_{\boldsymbol{a},-ij},\boldsymbol{\xi}_{\boldsymbol{b},-ij})}_{\boldsymbol{M}_{-\xi_{\boldsymbol{a}i},\xi_{\boldsymbol{a}j}}}[J(\boldsymbol{\gamma}^{T}\boldsymbol{M}_{\xi_{\boldsymbol{a}i},\xi_{\boldsymbol{a}j}})]-\phi(\xi_{\boldsymbol{a},i})\phi(\xi_{\boldsymbol{b},j})\overline{\mathrm{E}}^{(\boldsymbol{\xi}_{\boldsymbol{a},-ij},\boldsymbol{\xi}_{\boldsymbol{b},-ij})}_{\boldsymbol{M}_{-\xi_{\boldsymbol{a}i},\xi_{\boldsymbol{b}j}}}[J(\boldsymbol{\gamma}^{T}\boldsymbol{M}_{\xi_{\boldsymbol{a}i},\xi_{\boldsymbol{b}j}})]\\
\nonumber&+\gamma_{j}\phi(\xi_{\boldsymbol{a},i})\overline{\mathrm{E}}^{(\boldsymbol{\xi}_{\boldsymbol{a},-i},\boldsymbol{\xi}_{\boldsymbol{b},-i})}_{\boldsymbol{M}_{-\xi_{\boldsymbol{a}i}}}[J'(\boldsymbol{\gamma}^{T}\boldsymbol{M}_{\xi_{\boldsymbol{a}i}})]\\
\nonumber&+\phi(\xi_{\boldsymbol{b},i})\phi(\xi_{\boldsymbol{b},j})\overline{\mathrm{E}}^{(\boldsymbol{\xi}_{\boldsymbol{a},-ij},\boldsymbol{\xi}_{\boldsymbol{b},-ij})}_{\boldsymbol{M}_{-\xi_{\boldsymbol{b}i},\xi_{\boldsymbol{b}j}}}[J(\boldsymbol{\gamma}^{T}\boldsymbol{M}_{\xi_{\boldsymbol{b}i},\xi_{\boldsymbol{b}j}})]-\phi(\xi_{\boldsymbol{b},i})\phi(\xi_{\boldsymbol{a},j})\overline{\mathrm{E}}^{(\boldsymbol{\xi}_{\boldsymbol{a},-ij},\boldsymbol{\xi}_{\boldsymbol{b},-ij})}_{\boldsymbol{M}_{-\xi_{\boldsymbol{b}i},\xi_{\boldsymbol{a}j}}}[J(\boldsymbol{\gamma}^{T}\boldsymbol{M}_{\xi_{\boldsymbol{b}i},\xi_{\boldsymbol{a}j}})]\\
\nonumber&-\gamma_{j}\phi(\xi_{\boldsymbol{b},i})\overline{\mathrm{E}}^{(\boldsymbol{\xi}_{\boldsymbol{a},-i},\boldsymbol{\xi}_{\boldsymbol{b},-i})}_{\boldsymbol{M}_{-\xi_{\boldsymbol{b}i}}}[J'(\boldsymbol{\gamma}^{T}\boldsymbol{M}_{\xi_{\boldsymbol{b}i}})]\\
\nonumber&+\gamma_{i}\phi(\xi_{\boldsymbol{a},j})\overline{\mathrm{E}}^{(\boldsymbol{\xi}_{\boldsymbol{a},-j},\boldsymbol{\xi}_{\boldsymbol{b},-j})}_{\boldsymbol{M}_{-\xi_{\boldsymbol{a}j}}}[J'(\boldsymbol{\gamma}^{T}\boldsymbol{M}_{\xi_{\boldsymbol{a}j}})]-\gamma_{i}\phi(\xi_{\boldsymbol{b},j})\overline{\mathrm{E}}^{(\boldsymbol{\xi}_{\boldsymbol{a},-j},\boldsymbol{\xi}_{\boldsymbol{b},-j})}_{\boldsymbol{M}_{-\xi_{\boldsymbol{b}j}}}[J'(\boldsymbol{\gamma}^{T}\boldsymbol{M}_{\xi_{\boldsymbol{b}j}})]\\
&+\gamma_{i}\gamma_{j}\overline{\mathrm{E}}^{(\boldsymbol{\xi}_{\boldsymbol{a}},\boldsymbol{\xi}_{\boldsymbol{b}})}_{\mathbf{M}}[J''(\boldsymbol{\gamma}^{T}\boldsymbol{M})]\bigg\},~i\neq j,
\end{align*}
 \begin{align*}
\nonumber\Omega_{ii}=& 2\bigg\{\phi(\xi_{\boldsymbol{a},i})\xi_{\boldsymbol{a},i}\overline{\mathrm{E}}^{(\boldsymbol{\xi}_{\boldsymbol{a},-i},\boldsymbol{\xi}_{\boldsymbol{b},-i})}_{\boldsymbol{M}_{-\xi_{\boldsymbol{a}i}}}[J(\boldsymbol{\gamma}^{T}\boldsymbol{M}_{\xi_{\boldsymbol{a}i}})]-\phi(\xi_{\boldsymbol{b},i})\xi_{\boldsymbol{b},i}\overline{\mathrm{E}}^{(\boldsymbol{\xi}_{\boldsymbol{a},-i},\boldsymbol{\xi}_{\boldsymbol{b},-i})}_{\boldsymbol{M}_{-\xi_{\boldsymbol{b}i}}}[J(\boldsymbol{\gamma}^{T}\boldsymbol{M}_{\xi_{\boldsymbol{b}i}})]\\ \nonumber&+\gamma_{i}\phi(\xi_{\boldsymbol{a},i})\overline{\mathrm{E}}^{(\boldsymbol{\xi}_{\boldsymbol{a},-i},\boldsymbol{\xi}_{\boldsymbol{b},-i})}_{\boldsymbol{M}_{-\xi_{\boldsymbol{a}i}}}[J'(\boldsymbol{\gamma}^{T}\boldsymbol{M}_{\xi_{\boldsymbol{a}i}})]-\gamma_{i}\phi(\xi_{\boldsymbol{b},i})\overline{\mathrm{E}}^{(\boldsymbol{\xi}_{\boldsymbol{a},-i},\boldsymbol{\xi}_{\boldsymbol{b},-i})}_{\boldsymbol{M}_{-\xi_{\boldsymbol{b}i}}}[J'(\boldsymbol{\gamma}^{T}\boldsymbol{M}_{\xi_{\boldsymbol{b}i}})]\\
&+\gamma_{i}^{2}\overline{\mathrm{E}}^{(\boldsymbol{\xi}_{\boldsymbol{a}},\boldsymbol{\xi}_{\boldsymbol{b}})}_{\mathbf{M}}[J''(\boldsymbol{\gamma}^{T}\boldsymbol{M})]\bigg\}
+F_{\mathbf{Z}}(\boldsymbol{\xi}_{\boldsymbol{a}},\boldsymbol{\xi}_{\boldsymbol{b}})
,
\end{align*}
\begin{align*}
\nonumber\delta_{k}=&2\bigg\{\phi(\xi_{\boldsymbol{a},k})\overline{\mathrm{E}}^{(\boldsymbol{\xi}_{\boldsymbol{a},-k},\boldsymbol{\xi}_{\boldsymbol{b},-k})}_{\boldsymbol{M}_{-\xi_{\boldsymbol{a}k}}}[J(\boldsymbol{\gamma}^{T}\boldsymbol{M}_{\xi_{\boldsymbol{a}k}})]-\phi(\xi_{\boldsymbol{b},k})\overline{\mathrm{E}}^{(\boldsymbol{\xi}_{\boldsymbol{a},-k},\boldsymbol{\xi}_{\boldsymbol{b},-k})}_{\boldsymbol{M}_{-\xi_{\boldsymbol{b}k}}}[J(\boldsymbol{\gamma}^{T}\boldsymbol{M}_{\xi_{\boldsymbol{b}k}})]+\gamma_{k}\overline{\mathrm{E}}^{(\boldsymbol{\xi}_{\boldsymbol{a}},\boldsymbol{\xi}_{\boldsymbol{b}})}_{\mathbf{M}}[J'(\boldsymbol{\gamma}^{T}\boldsymbol{M})]\bigg\},
\end{align*}
$~i,j,k\in\{1,~2,\cdots,n\},$ $\mathbf{Z}\sim GSN_{n}\left(\boldsymbol{0},~\boldsymbol{I_{n}},~\boldsymbol{\gamma},~J\right)$, $\mathbf{M}\sim N_{n}(\boldsymbol{0},~\mathbf{I_{n}})$, $\mathbf{M}_{-\xi_{\boldsymbol{s}k}}\sim N_{n-1}\left(\boldsymbol{0},~\boldsymbol{I_{n-1}}\right)$, \\ $\mathbf{M}_{-\xi_{\boldsymbol{s}k},\xi_{\boldsymbol{t}l}}\sim N_{n-2}\left(\boldsymbol{0},~\boldsymbol{I_{n-2}}\right)$, and  $J'(\cdot)$ simply denotes the derivative of $J(\cdot)$.\\
$\mathbf{Example~1}$ (Skew-normal distribution).
  Letting $J(\cdot)=\Phi(\cdot)$ in Corollary~1. Thus,
\begin{align*}
\nonumber\Omega_{ij}=& 2\bigg\{\phi(\xi_{\boldsymbol{a},i})\phi(\xi_{\boldsymbol{a},j})\overline{\mathrm{E}}^{(\boldsymbol{\xi}_{\boldsymbol{a},-ij},\boldsymbol{\xi}_{\boldsymbol{b},-ij})}_{\boldsymbol{M}_{-\xi_{\boldsymbol{a}i},\xi_{\boldsymbol{a}j}}}[\Phi(\boldsymbol{\gamma}^{T}\boldsymbol{M}_{\xi_{\boldsymbol{a}i},\xi_{\boldsymbol{a}j}})]-\phi(\xi_{\boldsymbol{a},i})\phi(\xi_{\boldsymbol{b},j})\overline{\mathrm{E}}^{(\boldsymbol{\xi}_{\boldsymbol{a},-ij},\boldsymbol{\xi}_{\boldsymbol{b},-ij})}_{\boldsymbol{M}_{-\xi_{\boldsymbol{a}i},\xi_{\boldsymbol{b}j}}}[\Phi(\boldsymbol{\gamma}^{T}\boldsymbol{M}_{\xi_{\boldsymbol{a}i},\xi_{\boldsymbol{b}j}})]\\
\nonumber&+\gamma_{j}\phi(\xi_{\boldsymbol{a},i})\overline{\mathrm{E}}^{(\boldsymbol{\xi}_{\boldsymbol{a},-i},\boldsymbol{\xi}_{\boldsymbol{b},-i})}_{\boldsymbol{M}_{-\xi_{\boldsymbol{a}i}}}[\phi(\boldsymbol{\gamma}^{T}\boldsymbol{M}_{\xi_{\boldsymbol{a}i}})]\\
\nonumber&+\phi(\xi_{\boldsymbol{b},i})\phi(\xi_{\boldsymbol{b},j})\overline{\mathrm{E}}^{(\boldsymbol{\xi}_{\boldsymbol{a},-ij},\boldsymbol{\xi}_{\boldsymbol{b},-ij})}_{\boldsymbol{M}_{-\xi_{\boldsymbol{b}i},\xi_{\boldsymbol{b}j}}}[\Phi(\boldsymbol{\gamma}^{T}\boldsymbol{M}_{\xi_{\boldsymbol{b}i},\xi_{\boldsymbol{b}j}})]-\phi(\xi_{\boldsymbol{b},i})\phi(\xi_{\boldsymbol{a},j})\overline{\mathrm{E}}^{(\boldsymbol{\xi}_{\boldsymbol{a},-ij},\boldsymbol{\xi}_{\boldsymbol{b},-ij})}_{\boldsymbol{M}_{-\xi_{\boldsymbol{b}i},\xi_{\boldsymbol{a}j}}}[\Phi(\boldsymbol{\gamma}^{T}\boldsymbol{M}_{\xi_{\boldsymbol{b}i},\xi_{\boldsymbol{a}j}})]\\
\nonumber&-\gamma_{j}\phi(\xi_{\boldsymbol{b},i})\overline{\mathrm{E}}^{(\boldsymbol{\xi}_{\boldsymbol{a},-i},\boldsymbol{\xi}_{\boldsymbol{b},-i})}_{\boldsymbol{M}_{-\xi_{\boldsymbol{b}i}}}[\phi(\boldsymbol{\gamma}^{T}\boldsymbol{M}_{\xi_{\boldsymbol{b}i}})]\\
\nonumber&+\gamma_{i}\phi(\xi_{\boldsymbol{a},j})\overline{\mathrm{E}}^{(\boldsymbol{\xi}_{\boldsymbol{a},-j},\boldsymbol{\xi}_{\boldsymbol{b},-j})}_{\boldsymbol{M}_{-\xi_{\boldsymbol{a}j}}}[\phi(\boldsymbol{\gamma}^{T}\boldsymbol{M}_{\xi_{\boldsymbol{a}j}})]-\gamma_{i}\phi(\xi_{\boldsymbol{b},j})\overline{\mathrm{E}}^{(\boldsymbol{\xi}_{\boldsymbol{a},-j},\boldsymbol{\xi}_{\boldsymbol{b},-j})}_{\boldsymbol{M}_{-\xi_{\boldsymbol{b}j}}}[\phi(\boldsymbol{\gamma}^{T}\boldsymbol{M}_{\xi_{\boldsymbol{b}j}})]\\
&-\gamma_{i}\gamma_{j}\overline{\mathrm{E}}^{(\boldsymbol{\xi}_{\boldsymbol{a}},\boldsymbol{\xi}_{\boldsymbol{b}})}_{\mathbf{M}}[\boldsymbol{\gamma}^{T}\boldsymbol{M}\phi(\boldsymbol{\gamma}^{T}\boldsymbol{M})]\bigg\},~i\neq j,
\end{align*}
\begin{align*}
\nonumber\Omega_{ii}=& 2\bigg\{\phi(\xi_{\boldsymbol{a},i})\xi_{\boldsymbol{a},i}\overline{\mathrm{E}}^{(\boldsymbol{\xi}_{\boldsymbol{a},-i},\boldsymbol{\xi}_{\boldsymbol{b},-i})}_{\boldsymbol{M}_{-\xi_{\boldsymbol{a}i}}}[\Phi(\boldsymbol{\gamma}^{T}\boldsymbol{M}_{\xi_{\boldsymbol{a}i}})]-\phi(\xi_{\boldsymbol{b},i})\xi_{\boldsymbol{b},i}\overline{\mathrm{E}}^{(\boldsymbol{\xi}_{\boldsymbol{a},-i},\boldsymbol{\xi}_{\boldsymbol{b},-i})}_{\boldsymbol{M}_{-\xi_{\boldsymbol{b}i}}}[\Phi(\boldsymbol{\gamma}^{T}\boldsymbol{M}_{\xi_{\boldsymbol{b}i}})]\\ \nonumber&+\gamma_{i}\phi(\xi_{\boldsymbol{a},i})\overline{\mathrm{E}}^{(\boldsymbol{\xi}_{\boldsymbol{a},-i},\boldsymbol{\xi}_{\boldsymbol{b},-i})}_{\boldsymbol{M}_{-\xi_{\boldsymbol{a}i}}}[\phi(\boldsymbol{\gamma}^{T}\boldsymbol{M}_{\xi_{\boldsymbol{a}i}})]-\gamma_{i}\phi(\xi_{\boldsymbol{b},i})\overline{\mathrm{E}}^{(\boldsymbol{\xi}_{\boldsymbol{a},-i},\boldsymbol{\xi}_{\boldsymbol{b},-i})}_{\boldsymbol{M}_{-\xi_{\boldsymbol{b}i}}}[\phi(\boldsymbol{\gamma}^{T}\boldsymbol{M}_{\xi_{\boldsymbol{b}i}})]\\
&-\gamma_{i}^{2}\overline{\mathrm{E}}^{(\boldsymbol{\xi}_{\boldsymbol{a}},\boldsymbol{\xi}_{\boldsymbol{b}})}_{\mathbf{M}}[\boldsymbol{\gamma}^{T}\boldsymbol{M}\phi(\boldsymbol{\gamma}^{T}\boldsymbol{M})]\bigg\}
+F_{\mathbf{Z}}(\boldsymbol{\xi}_{\boldsymbol{a}},\boldsymbol{\xi}_{\boldsymbol{b}})
,
\end{align*}
\begin{align*}
\nonumber\delta_{k}=&2\bigg\{\phi(\xi_{\boldsymbol{a},k})\overline{\mathrm{E}}^{(\boldsymbol{\xi}_{\boldsymbol{a},-k},\boldsymbol{\xi}_{\boldsymbol{b},-k})}_{\boldsymbol{M}_{-\xi_{\boldsymbol{a}k}}}[\Phi(\boldsymbol{\gamma}^{T}\boldsymbol{M}_{\xi_{\boldsymbol{a}k}})]-\phi(\xi_{\boldsymbol{b},k})\overline{\mathrm{E}}^{(\boldsymbol{\xi}_{\boldsymbol{a},-k},\boldsymbol{\xi}_{\boldsymbol{b},-k})}_{\boldsymbol{M}_{-\xi_{\boldsymbol{b}k}}}[\Phi(\boldsymbol{\gamma}^{T}\boldsymbol{M}_{\xi_{\boldsymbol{b}k}})]+\gamma_{k}\overline{\mathrm{E}}^{(\boldsymbol{\xi}_{\boldsymbol{a}},\boldsymbol{\xi}_{\boldsymbol{b}})}_{\mathbf{M}}[\phi(\boldsymbol{\gamma}^{T}\boldsymbol{M})]\bigg\},
\end{align*}
$~i,j,k\in\{1,~2,\cdots,n\},$ $\mathbf{Z}\sim SN_{n}\left(\boldsymbol{0},~\boldsymbol{I_{n}},~\boldsymbol{\gamma}\right)$, $\mathbf{M}$, $\mathbf{M}_{-\xi_{\boldsymbol{s}k}}$ and $\mathbf{M}_{-\xi_{\boldsymbol{s}k},\xi_{\boldsymbol{t}l}}$ are given in Corollary~1.\\
\noindent$\mathbf{Corollary~2}$ (GSSt distribution). Let $\mathbf{Y}\sim GSSt_{n}(\boldsymbol{\mu},~\mathbf{\Sigma},~m,~\boldsymbol{\gamma},~J)$. Here $m>0$ is degrees of freedom, $J: \mathbb{R}\rightarrow\mathbb{R}$  and $\boldsymbol{\gamma}=(\gamma_{1},~\gamma_{2},~\cdots,~\gamma_{n})^{T}\in\mathbb{R}^{n}$. In this case,  $$g_{n}(u)=\left(1+\frac{2u}{m}\right)^{-(m+n)/2},~c_{n}=\frac{\Gamma\left(\frac{m+n}{2}\right)}{\Gamma(m/2)(m\pi)^{\frac{n}{2}}},$$
and $$H\left(\mathbf{\Sigma}^{-\frac{1}{2}}(\boldsymbol{y}-\boldsymbol{\mu})\right)=J\left(\boldsymbol{\gamma}^{T}\mathbf{\Sigma}^{-\frac{1}{2}}(\boldsymbol{y}-\boldsymbol{\mu})\right).$$
 So $\overline{G}_{n}(t)$ and $\overline{\mathcal{G}}_{n}(t)$ can be expressed, respectively, as
 $$\overline{G}_{n}(t)=\frac{m}{m+n-2}\left(1+\frac{2t}{m}\right)^{-(m+n-2)/2}$$
  and
 $$\overline{\mathcal{G}}_{n}(t)=\frac{m^{2}}{(m+n-2)(m+n-4)}\left(1+\frac{2t}{m}\right)^{-(m+n-4)/2}.$$ In addition, \begin{align*}
 c_{n}^{\ast}&=\frac{(m+n-2)\Gamma(n/2)}{(2\pi)^{n/2}m}\left[\int_{0}^{\infty}t^{n/2-1}\left(1+\frac{2t}{m}\right)^{-(m+n-2)/2}\mathrm{d}t\right]^{-1}\\
 &=\frac{(m+n-2)\Gamma(n/2)}{(m\pi)^{n/2}mB(\frac{n}{2},~\frac{m-2}{2})},~if~m>2
 \end{align*} and
 \begin{align*}
 c_{n}^{\ast\ast}&=\frac{(m+n-2)(m+n-4)\Gamma(n/2)}{(2\pi)^{n/2}m^{2}}\left[\int_{0}^{\infty}t^{n/2-1}\left(1+\frac{2t}{m}\right)^{-(m+n-4)/2}\mathrm{d}t\right]^{-1}\\
 &=\frac{(m+n-2)(m+n-4)\Gamma(n/2)}{(m\pi)^{n/2}m^{2}B(\frac{n}{2},~\frac{m-4}{2})},~if~m>4,
 \end{align*}
 where $\Gamma(\cdot)$ and $B(\cdot,\cdot)$ denote Gamma function and Beta function, respectively. Since
  \begin{align*}
f_{\mathbf{M^{\ast}}_{-\xi_{\boldsymbol{s}k}}}(\boldsymbol{w})&=\frac{c_{n-1,\xi_{\boldsymbol{s}k}}^{\ast}m}{m+n-2}\left(1+\frac{\xi_{\boldsymbol{s},k}^{2}}{m}\right)^{-\frac{m+n-2}{2}}\left(1+\frac{\boldsymbol{w}^{T}\Delta_{\xi_{\boldsymbol{s}k}}^{-1}\boldsymbol{w}}{m-1}\right)^{-\frac{m+n-2}{2}}\\
 &=St_{n-1}(\boldsymbol{0},~\Delta_{\xi_{\boldsymbol{s}k}},~m-1),~\boldsymbol{w}\in \mathbb{R}^{n-1},
\end{align*}
\begin{align*}
f_{\mathbf{M^{\ast\ast}}_{-\xi_{\boldsymbol{s}k}}}(\boldsymbol{v})&=\frac{c_{n-1,\xi_{\boldsymbol{s}k}}^{\ast\ast}m^{2}}{(m+n-2)(m+n-4)}\left(1+\frac{\xi_{\boldsymbol{s},k}^{2}}{m}\right)^{-\frac{m+n-4}{2}}\left(1+\frac{\boldsymbol{v}^{T}\Lambda_{\xi_{\boldsymbol{s}k}}^{-1}\boldsymbol{v}}{m-3}\right)^{-\frac{m+n-4}{2}}\\
 &=St_{n-1}(\boldsymbol{0},~\Lambda_{\xi_{\boldsymbol{s}k}},~m-3),~\boldsymbol{v}\in \mathbb{R}^{n-1}
\end{align*}
and
\begin{align*}
f_{\mathbf{M^{\ast\ast}}_{-\xi_{\boldsymbol{s}k},\xi_{\boldsymbol{t}l}}}(\boldsymbol{u})&=\frac{c_{n-2,\xi_{\boldsymbol{s}k},\xi_{\boldsymbol{t}l}}^{\ast\ast}m^{2}}{(m+n-2)(m+n-4)}\left(1+\frac{\xi_{\boldsymbol{s},k}^{2}+\xi_{\boldsymbol{t},l}^{2}}{m}\right)^{-\frac{m+n-4}{2}}\left(1+\frac{\boldsymbol{u}^{T}\Theta_{\xi_{\boldsymbol{s}k},\xi_{\boldsymbol{t}l}}^{-1}\boldsymbol{u}}{m-2}\right)^{-\frac{m+n-4}{2}}\\
 &=St_{n-2}(\boldsymbol{0},~\Theta_{\xi_{\boldsymbol{s}k},\xi_{\boldsymbol{t}l}},~m-2),~\boldsymbol{u}\in \mathbb{R}^{n-2},
\end{align*}
so that
$$c_{n-1,\xi_{\boldsymbol{s}k}}^{\ast}=\frac{\Gamma\left(\frac{m+n-2}{2}\right)(m+n-2)}{\Gamma\left(\frac{m-1}{2}\right)\pi^{\frac{n-1}{2}}(m-1)^{\frac{n-1}{2}}m^{\frac{n+m}{2}}}\left(m+\xi_{\boldsymbol{s},k}^{2}\right)^{\frac{m+n-2}{2}},$$
$$c_{n-1,\xi_{\boldsymbol{s}k}}^{\ast\ast}=\frac{\Gamma\left(\frac{m+n-4}{2}\right)(m+n-2)(m+n-4)}{\Gamma\left(\frac{m-3}{2}\right)\pi^{\frac{n-1}{2}}m^{\frac{m+n}{2}}}\left(m+\xi_{\boldsymbol{s},k}^{2}\right)^{\frac{m+n-4}{2}}$$
 and
 $$c_{n-2,\xi_{\boldsymbol{s}k},\xi_{\boldsymbol{t}l}}^{\ast\ast}=\frac{\Gamma\left(\frac{m+n-4}{2}\right)(m+n-2)(m+n-4)}{\Gamma\left(\frac{m-2}{2}\right)\pi^{\frac{n-2}{2}}m^{\frac{m+n}{2}}}\left(m+\xi_{\boldsymbol{s},k}^{2}+\xi_{\boldsymbol{t},l}^{2}\right)^{\frac{m+n-4}{2}}.$$
Then
\begin{align*}
\nonumber\Omega_{ij}=& 2\bigg\{\frac{c_{n}}{c_{n-2,\xi_{\boldsymbol{a}i},\xi_{\boldsymbol{a}j}}^{\ast\ast}}\overline{\mathrm{E}}^{(\boldsymbol{\xi}_{\boldsymbol{a},-ij},\boldsymbol{\xi}_{\boldsymbol{b},-ij})}_{\boldsymbol{M}_{-\xi_{\boldsymbol{a}i},\xi_{\boldsymbol{a}j}}^{\ast\ast}}[J(\boldsymbol{\gamma}^{T}\boldsymbol{M}_{\xi_{\boldsymbol{a}i},\xi_{\boldsymbol{a}j}}^{\ast\ast})]-\frac{c_{n}}{c_{n-2,\xi_{\boldsymbol{a}i},\xi_{\boldsymbol{b}j}}^{\ast\ast}}\overline{\mathrm{E}}^{(\boldsymbol{\xi}_{\boldsymbol{a},-ij},\boldsymbol{\xi}_{\boldsymbol{b},-ij})}_{\boldsymbol{M}_{-\xi_{\boldsymbol{a}i},\xi_{\boldsymbol{b}j}}^{\ast\ast}}[J(\boldsymbol{\gamma}^{T}\boldsymbol{M}_{\xi_{\boldsymbol{a}i},\xi_{\boldsymbol{b}j}}^{\ast\ast})]\\
\nonumber&+\frac{c_{n}\gamma_{j}}{c_{n-1,\xi_{\boldsymbol{a}i}}^{\ast\ast}}\overline{\mathrm{E}}^{(\boldsymbol{\xi}_{\boldsymbol{a},-i},\boldsymbol{\xi}_{\boldsymbol{b},-i})}_{\boldsymbol{M}_{-\xi_{\boldsymbol{a}i}}^{\ast\ast}}[J'(\boldsymbol{\gamma}^{T}\boldsymbol{M}_{\xi_{\boldsymbol{a}i}}^{\ast\ast})]\\
\nonumber&+\frac{c_{n}}{c_{n-2,\xi_{\boldsymbol{b}i},\xi_{\boldsymbol{b}j}}^{\ast\ast}}\overline{\mathrm{E}}^{(\boldsymbol{\xi}_{\boldsymbol{a},-ij},\boldsymbol{\xi}_{\boldsymbol{b},-ij})}_{\boldsymbol{M}_{-\xi_{\boldsymbol{b}i},\xi_{\boldsymbol{b}j}}^{\ast\ast}}[J(\boldsymbol{\gamma}^{T}\boldsymbol{M}_{\xi_{\boldsymbol{b}i},\xi_{\boldsymbol{b}j}}^{\ast\ast})]-\frac{c_{n}}{c_{n-2,\xi_{\boldsymbol{b}i},\xi_{\boldsymbol{a}j}}^{\ast\ast}}\overline{\mathrm{E}}^{(\boldsymbol{\xi}_{\boldsymbol{a},-ij},\boldsymbol{\xi}_{\boldsymbol{b},-ij})}_{\boldsymbol{M}_{-\xi_{\boldsymbol{b}i},\xi_{\boldsymbol{a}j}}^{\ast\ast}}[J(\boldsymbol{\gamma}^{T}\boldsymbol{M}_{\xi_{\boldsymbol{b}i},\xi_{\boldsymbol{a}j}}^{\ast\ast})]\\
\nonumber&-\frac{c_{n}\gamma_{j}}{c_{n-1,\xi_{\boldsymbol{b}i}}^{\ast\ast}}\overline{\mathrm{E}}^{(\boldsymbol{\xi}_{\boldsymbol{a},-i},\boldsymbol{\xi}_{\boldsymbol{b},-i})}_{\boldsymbol{M}_{-\xi_{\boldsymbol{b}i}}^{\ast\ast}}[J'(\boldsymbol{\gamma}^{T}\boldsymbol{M}_{\xi_{\boldsymbol{b}i}}^{\ast\ast})]\\
\nonumber&+\frac{c_{n}\gamma_{i}}{c_{n-1,\xi_{\boldsymbol{a}j}}^{\ast\ast}}\overline{\mathrm{E}}^{(\boldsymbol{\xi}_{\boldsymbol{a},-j},\boldsymbol{\xi}_{\boldsymbol{b},-j})}_{\boldsymbol{M}_{-\xi_{\boldsymbol{a}j}}^{\ast\ast}}[J'(\boldsymbol{\gamma}^{T}\boldsymbol{M}_{\xi_{\boldsymbol{a}j}}^{\ast\ast})]-\frac{c_{n}\gamma_{i}}{c_{n-1,\xi_{\boldsymbol{b}j}}^{\ast\ast}}\overline{\mathrm{E}}^{(\boldsymbol{\xi}_{\boldsymbol{a},-j},\boldsymbol{\xi}_{\boldsymbol{b},-j})}_{\boldsymbol{M}_{-\xi_{\boldsymbol{b}j}}^{\ast\ast}}[J'(\boldsymbol{\gamma}^{T}\boldsymbol{M}_{\xi_{\boldsymbol{b}j}}^{\ast\ast})]\\
&+\frac{m^{2}\gamma_{i}\gamma_{j}}{(m-2)(m-4)}\overline{\mathrm{E}}^{(\boldsymbol{\xi}_{\boldsymbol{a}},\boldsymbol{\xi}_{\boldsymbol{b}})}_{\mathbf{M}^{\ast\ast}}[J''(\boldsymbol{\gamma}^{T}\boldsymbol{M}^{\ast\ast})]\bigg\},~i\neq j,~m>4,
\end{align*}
\begin{align*}
\nonumber\Omega_{ii}=& 2\bigg\{\frac{c_{n}\xi_{\boldsymbol{a},i}}{c_{n-1,\xi_{\boldsymbol{a}i}}^{\ast}}\overline{\mathrm{E}}^{(\boldsymbol{\xi}_{\boldsymbol{a},-i},\boldsymbol{\xi}_{\boldsymbol{b},-i})}_{\boldsymbol{M}_{-\xi_{\boldsymbol{a}i}}^{\ast}}[J(\boldsymbol{\gamma}^{T}\boldsymbol{M}_{\xi_{\boldsymbol{a}i}}^{\ast})]-\frac{c_{n}\xi_{\boldsymbol{b},i}}{c_{n-1,\xi_{\boldsymbol{b}i}}^{\ast}}\overline{\mathrm{E}}^{(\boldsymbol{\xi}_{\boldsymbol{a},-i},\boldsymbol{\xi}_{\boldsymbol{b},-i})}_{\boldsymbol{M}_{-\xi_{\boldsymbol{b}i}}^{\ast}}[J(\boldsymbol{\gamma}^{T}\boldsymbol{M}_{\xi_{\boldsymbol{b}i}}^{\ast})]\\ \nonumber&+\frac{c_{n}\gamma_{i}}{c_{n-1,\xi_{\boldsymbol{a}i}}^{\ast\ast}}\overline{\mathrm{E}}^{(\boldsymbol{\xi}_{\boldsymbol{a},-i},\boldsymbol{\xi}_{\boldsymbol{b},-i})}_{\boldsymbol{M}_{-\xi_{\boldsymbol{a}i}}^{\ast\ast}}[J'(\boldsymbol{\gamma}^{T}\boldsymbol{M}_{\xi_{\boldsymbol{a}i}}^{\ast\ast})]-\frac{c_{n}\gamma_{i}}{c_{n-1,\xi_{\boldsymbol{b}i}}^{\ast\ast}}\overline{\mathrm{E}}^{(\boldsymbol{\xi}_{\boldsymbol{a},-i},\boldsymbol{\xi}_{\boldsymbol{b},-i})}_{\boldsymbol{M}_{-\xi_{\boldsymbol{b}i}}^{\ast\ast}}[J'(\boldsymbol{\gamma}^{T}\boldsymbol{M}_{\xi_{\boldsymbol{b}i}}^{\ast\ast})]\\
&+\frac{m^{2}\gamma_{i}^{2}}{(m-2)(m-4)}\overline{\mathrm{E}}^{(\boldsymbol{\xi}_{\boldsymbol{a}},\boldsymbol{\xi}_{\boldsymbol{b}})}_{\mathbf{M}^{\ast\ast}}[J''(\boldsymbol{\gamma}^{T}\boldsymbol{M}^{\ast\ast})]\bigg\}
+\frac{m}{m-2}F_{\mathbf{Z^{\ast}}}(\boldsymbol{\xi}_{\boldsymbol{a}},\boldsymbol{\xi}_{\boldsymbol{b}})
,~m>4,
\end{align*}
\begin{align*}
\nonumber\delta_{k}=&2\bigg\{\frac{c_{n}}{c_{n-1,\xi_{\boldsymbol{a}k}}^{\ast}}\overline{\mathrm{E}}^{(\boldsymbol{\xi}_{\boldsymbol{a},-k},\boldsymbol{\xi}_{\boldsymbol{b},-k})}_{\boldsymbol{M}_{-\xi_{\boldsymbol{a}k}}^{\ast}}[J(\boldsymbol{\gamma}^{T}\boldsymbol{M}_{\xi_{\boldsymbol{a}k}}^{\ast})]-\frac{c_{n}}{c_{n-1,\xi_{\boldsymbol{b}k}}^{\ast}}\overline{\mathrm{E}}^{(\boldsymbol{\xi}_{\boldsymbol{a},-k},\boldsymbol{\xi}_{\boldsymbol{b},-k})}_{\boldsymbol{M}_{-\xi_{\boldsymbol{b}k}}^{\ast}}[J(\boldsymbol{\gamma}^{T}\boldsymbol{M}_{\xi_{\boldsymbol{b}k}}^{\ast})]\\
&+\frac{m\gamma_{k}}{m-2}\overline{\mathrm{E}}^{(\boldsymbol{\xi}_{\boldsymbol{a}},\boldsymbol{\xi}_{\boldsymbol{b}})}_{\mathbf{M}^{\ast}}[J'(\boldsymbol{\gamma}^{T}\boldsymbol{M}^{\ast})]\bigg\},~m>2,
\end{align*}
 $i,~j,~k,l\in\{1,~2,\cdots,n\}$, $\mathbf{Z}\sim GSSt_{n}(\boldsymbol{0},~\mathbf{I_{n}},~m,~\boldsymbol{\gamma},~J),$ $\mathbf{Z}^{\ast}\sim GSE_{n}\left(\boldsymbol{0},~\boldsymbol{I_{n}},~m,~\overline{G}_{n},~\boldsymbol{\gamma},~J\right)$,\\
 $\mathbf{M}^{\ast}\sim E_{n}(\boldsymbol{0},~\mathbf{I_{n}},~\overline{G}_{n})$, $\mathbf{M}^{\ast\ast}\sim E_{n}(\boldsymbol{0},~\mathbf{I_{n}},~\overline{\mathcal{G}}_{n})$,  $\mathbf{M^{\ast\ast}}_{-\xi_{\boldsymbol{s}k},\xi_{\boldsymbol{t}l}}\sim St_{n-2}\left(\boldsymbol{0},~\Theta_{\xi_{\boldsymbol{s}k},\xi_{\boldsymbol{t}l}},~m-2\right)$,\\ $\mathbf{M^{\ast\ast}}_{-\xi_{\boldsymbol{s}k}}\sim St_{n-1}\left(\boldsymbol{0},~\Lambda_{\xi_{\boldsymbol{s}k}},~m-3\right)$, $\mathbf{M^{\ast}}_{-\xi_{\boldsymbol{s}k}}\sim St_{n-1}(\boldsymbol{0},~\mathbf{\bigtriangleup}_{\xi_{\boldsymbol{s}k}},~m-1),$
$\bigtriangleup_{\xi_{\boldsymbol{s}k}}=\left(\frac{m+\xi_{\boldsymbol{s},k}^{2}}{m-1}\right)\mathbf{I}_{n-1} $,\\ $\Lambda_{\xi_{\boldsymbol{s}k}}=\left(\frac{m+\xi_{\boldsymbol{s},k}^{2}}{m-3}\right)\mathbf{I}_{n-1}$, $\Theta_{\xi_{\boldsymbol{s}k},\xi_{\boldsymbol{t}l}}=\left(\frac{m+\xi_{\boldsymbol{s},k}^{2}+\xi_{\boldsymbol{t},l}^{2}}{m-2}\right)\mathbf{I}_{n-2}$ and $\mathbf{M}\sim St_{n}(\boldsymbol{0},~\mathbf{I_{n}},m)$.\\
Therefore, $\frac{c_{n}}{c_{n-1,\xi_{\boldsymbol{s}k}}^{\ast}}$, $\frac{c_{n}}{c_{n-1,\xi_{\boldsymbol{s}k}}^{\ast\ast}}$ and $\frac{c_{n}}{c_{n-2,\xi_{\boldsymbol{s}k},\xi_{\boldsymbol{t}l}}^{\ast\ast}}$ can be further simplified as
$$\frac{c_{n}}{c_{n-1,\xi_{\boldsymbol{s}k}}^{\ast}}=\frac{\Gamma(\frac{m-1}{2})\left(\frac{m-1}{m}\right)^{(n-1)/2}}{\Gamma(\frac{m}{2})\sqrt{\pi/m}}\left(1+\frac{\xi_{\boldsymbol{s},k}^{2}}{m}\right)^{-(m+n-2)/2},$$
$$\frac{c_{n}}{c_{n-1,\xi_{\boldsymbol{s}k}}^{\ast\ast}}=\frac{\Gamma(\frac{m-3}{2})(m+n)m^{\frac{m}{2}}}{4\Gamma(\frac{m}{2})(m+n-4)\sqrt{\pi}}\left(m+\xi_{\boldsymbol{s},k}^{2}\right)^{-(m+n-4)/2}$$
and
$$\frac{c_{n}}{c_{n-2,\xi_{\boldsymbol{s}k},\xi_{\boldsymbol{t}l}}^{\ast\ast}}=\frac{(m+n)m^{\frac{m-2}{2}}}{2(m+n-4)\pi}\left(m+\xi_{\boldsymbol{s},k}^{2}+\xi_{\boldsymbol{t},l}^{2}\right)^{-(m+n-4)/2}.$$
$\mathbf{Example~2}$ (Skew student-$t$ distribution) Letting $J(\cdot)=T(\cdot)\sim St_{1}(0,~1,~m)$ in Corollary~2. Thus,
\begin{align*}
\nonumber\Omega_{ij}=& 2\bigg\{\frac{c_{n}}{c_{n-2,\xi_{\boldsymbol{a}i},\xi_{\boldsymbol{a}j}}^{\ast\ast}}\overline{\mathrm{E}}^{(\boldsymbol{\xi}_{\boldsymbol{a},-ij},\boldsymbol{\xi}_{\boldsymbol{b},-ij})}_{\boldsymbol{M}_{-\xi_{\boldsymbol{a}i},\xi_{\boldsymbol{a}j}}^{\ast\ast}}[T(\boldsymbol{\gamma}^{T}\boldsymbol{M}_{\xi_{\boldsymbol{a}i},\xi_{\boldsymbol{a}j}}^{\ast\ast})]-\frac{c_{n}}{c_{n-2,\xi_{\boldsymbol{a}i},\xi_{\boldsymbol{b}j}}^{\ast\ast}}\overline{\mathrm{E}}^{(\boldsymbol{\xi}_{\boldsymbol{a},-ij},\boldsymbol{\xi}_{\boldsymbol{b},-ij})}_{\boldsymbol{M}_{-\xi_{\boldsymbol{a}i},\xi_{\boldsymbol{b}j}}^{\ast\ast}}[T(\boldsymbol{\gamma}^{T}\boldsymbol{M}_{\xi_{\boldsymbol{a}i},\xi_{\boldsymbol{b}j}}^{\ast\ast})]\\
\nonumber&+\frac{c_{n}\gamma_{j}}{c_{n-1,\xi_{\boldsymbol{a}i}}^{\ast\ast}}\overline{\mathrm{E}}^{(\boldsymbol{\xi}_{\boldsymbol{a},-i},\boldsymbol{\xi}_{\boldsymbol{b},-i})}_{\boldsymbol{M}_{-\xi_{\boldsymbol{a}i}}^{\ast\ast}}[\boldsymbol{t}(\boldsymbol{\gamma}^{T}\boldsymbol{M}_{\xi_{\boldsymbol{a}i}}^{\ast\ast})]\\
\nonumber&+\frac{c_{n}}{c_{n-2,\xi_{\boldsymbol{b}i},\xi_{\boldsymbol{b}j}}^{\ast\ast}}\overline{\mathrm{E}}^{(\boldsymbol{\xi}_{\boldsymbol{a},-ij},\boldsymbol{\xi}_{\boldsymbol{b},-ij})}_{\boldsymbol{M}_{-\xi_{\boldsymbol{b}i},\xi_{\boldsymbol{b}j}}^{\ast\ast}}[T(\boldsymbol{\gamma}^{T}\boldsymbol{M}_{\xi_{\boldsymbol{b}i},\xi_{\boldsymbol{b}j}}^{\ast\ast})]-\frac{c_{n}}{c_{n-2,\xi_{\boldsymbol{b}i},\xi_{\boldsymbol{a}j}}^{\ast\ast}}\overline{\mathrm{E}}^{(\boldsymbol{\xi}_{\boldsymbol{a},-ij},\boldsymbol{\xi}_{\boldsymbol{b},-ij})}_{\boldsymbol{M}_{-\xi_{\boldsymbol{b}i},\xi_{\boldsymbol{a}j}}^{\ast\ast}}[T(\boldsymbol{\gamma}^{T}\boldsymbol{M}_{\xi_{\boldsymbol{b}i},\xi_{\boldsymbol{a}j}}^{\ast\ast})]\\
\nonumber&-\frac{c_{n}\gamma_{j}}{c_{n-1,\xi_{\boldsymbol{b}i}}^{\ast\ast}}\overline{\mathrm{E}}^{(\boldsymbol{\xi}_{\boldsymbol{a},-i},\boldsymbol{\xi}_{\boldsymbol{b},-i})}_{\boldsymbol{M}_{-\xi_{\boldsymbol{b}i}}^{\ast\ast}}[\boldsymbol{t}(\boldsymbol{\gamma}^{T}\boldsymbol{M}_{\xi_{\boldsymbol{b}i}}^{\ast\ast})]\\
\nonumber&+\frac{c_{n}\gamma_{i}}{c_{n-1,\xi_{\boldsymbol{a}j}}^{\ast\ast}}\overline{\mathrm{E}}^{(\boldsymbol{\xi}_{\boldsymbol{a},-j},\boldsymbol{\xi}_{\boldsymbol{b},-j})}_{\boldsymbol{M}_{-\xi_{\boldsymbol{a}j}}^{\ast\ast}}[\boldsymbol{t}(\boldsymbol{\gamma}^{T}\boldsymbol{M}_{\xi_{\boldsymbol{a}j}}^{\ast\ast})]-\frac{c_{n}\gamma_{i}}{c_{n-1,\xi_{\boldsymbol{b}j}}^{\ast\ast}}\overline{\mathrm{E}}^{(\boldsymbol{\xi}_{\boldsymbol{a},-j},\boldsymbol{\xi}_{\boldsymbol{b},-j})}_{\boldsymbol{M}_{-\xi_{\boldsymbol{b}j}}^{\ast\ast}}[\boldsymbol{t}(\boldsymbol{\gamma}^{T}\boldsymbol{M}_{\xi_{\boldsymbol{b}j}}^{\ast\ast})]\\
&-\frac{m(m+1)\gamma_{i}\gamma_{j}}{(m-2)(m-4)}\overline{\mathrm{E}}^{(\boldsymbol{\xi}_{\boldsymbol{a}},\boldsymbol{\xi}_{\boldsymbol{b}})}_{\mathbf{M}^{\ast\ast}}\left[\left(1+\frac{(\boldsymbol{\gamma}^{T}\boldsymbol{M}^{\ast\ast})^{2}}{m}\right)^{-1}\boldsymbol{\gamma}^{T}\boldsymbol{M}^{\ast\ast}\boldsymbol{t}(\boldsymbol{\gamma}^{T}\boldsymbol{M}^{\ast\ast})\right]\bigg\},~i\neq j,~m>4,
\end{align*}
\begin{align*}
\nonumber\Omega_{ii}=& 2\bigg\{\frac{c_{n}\xi_{\boldsymbol{a},i}}{c_{n-1,\xi_{\boldsymbol{a}i}}^{\ast}}\overline{\mathrm{E}}^{(\boldsymbol{\xi}_{\boldsymbol{a},-i},\boldsymbol{\xi}_{\boldsymbol{b},-i})}_{\boldsymbol{M}_{-\xi_{\boldsymbol{a}i}}^{\ast}}[T(\boldsymbol{\gamma}^{T}\boldsymbol{M}_{\xi_{\boldsymbol{a}i}}^{\ast})]-\frac{c_{n}\xi_{\boldsymbol{b},i}}{c_{n-1,\xi_{\boldsymbol{b}i}}^{\ast}}\overline{\mathrm{E}}^{(\boldsymbol{\xi}_{\boldsymbol{a},-i},\boldsymbol{\xi}_{\boldsymbol{b},-i})}_{\boldsymbol{M}_{-\xi_{\boldsymbol{b}i}}^{\ast}}[T(\boldsymbol{\gamma}^{T}\boldsymbol{M}_{\xi_{\boldsymbol{b}i}}^{\ast})]\\ \nonumber&+\frac{c_{n}\gamma_{i}}{c_{n-1,\xi_{\boldsymbol{a}i}}^{\ast\ast}}\overline{\mathrm{E}}^{(\boldsymbol{\xi}_{\boldsymbol{a},-i},\boldsymbol{\xi}_{\boldsymbol{b},-i})}_{\boldsymbol{M}_{-\xi_{\boldsymbol{a}i}}^{\ast\ast}}[\boldsymbol{t}(\boldsymbol{\gamma}^{T}\boldsymbol{M}_{\xi_{\boldsymbol{a}i}}^{\ast\ast})]-\frac{c_{n}\gamma_{i}}{c_{n-1,\xi_{\boldsymbol{b}i}}^{\ast\ast}}\overline{\mathrm{E}}^{(\boldsymbol{\xi}_{\boldsymbol{a},-i},\boldsymbol{\xi}_{\boldsymbol{b},-i})}_{\boldsymbol{M}_{-\xi_{\boldsymbol{b}i}}^{\ast\ast}}[\boldsymbol{t}(\boldsymbol{\gamma}^{T}\boldsymbol{M}_{\xi_{\boldsymbol{b}i}}^{\ast\ast})]\\
&-\frac{m(m+1)\gamma_{i}^{2}}{(m-2)(m-4)}\overline{\mathrm{E}}^{(\boldsymbol{\xi}_{\boldsymbol{a}},\boldsymbol{\xi}_{\boldsymbol{b}})}_{\mathbf{M}^{\ast\ast}}\left[\left(1+\frac{(\boldsymbol{\gamma}^{T}\boldsymbol{M}^{\ast\ast})^{2}}{m}\right)^{-1}\boldsymbol{\gamma}^{T}\boldsymbol{M}^{\ast\ast}\boldsymbol{t}(\boldsymbol{\gamma}^{T}\boldsymbol{M}^{\ast\ast})\right]\bigg\}\\
&+\frac{m}{m-2}F_{\mathbf{Z^{\ast}}}(\boldsymbol{\xi}_{\boldsymbol{a}},\boldsymbol{\xi}_{\boldsymbol{b}})
,~m>4,
\end{align*}
\begin{align*}
\nonumber\delta_{k}=&2\bigg\{\frac{c_{n}}{c_{n-1,\xi_{\boldsymbol{a}k}}^{\ast}}\overline{\mathrm{E}}^{(\boldsymbol{\xi}_{\boldsymbol{a},-k},\boldsymbol{\xi}_{\boldsymbol{b},-k})}_{\boldsymbol{M}_{-\xi_{\boldsymbol{a}k}}^{\ast}}[T(\boldsymbol{\gamma}^{T}\boldsymbol{M}_{\xi_{\boldsymbol{a}k}}^{\ast})]-\frac{c_{n}}{c_{n-1,\xi_{\boldsymbol{b}k}}^{\ast}}\overline{\mathrm{E}}^{(\boldsymbol{\xi}_{\boldsymbol{a},-k},\boldsymbol{\xi}_{\boldsymbol{b},-k})}_{\boldsymbol{M}_{-\xi_{\boldsymbol{b}k}}^{\ast}}[T(\boldsymbol{\gamma}^{T}\boldsymbol{M}_{\xi_{\boldsymbol{b}k}}^{\ast})]\\
&+\frac{m\gamma_{k}}{m-2}\overline{\mathrm{E}}^{(\boldsymbol{\xi}_{\boldsymbol{a}},\boldsymbol{\xi}_{\boldsymbol{b}})}_{\mathbf{M}^{\ast}}[\boldsymbol{t}(\boldsymbol{\gamma}^{T}\boldsymbol{M}^{\ast})]\bigg\},~m>2,
\end{align*}
 $i,~j,~k\in\{1,~2,\cdots,n\}$,
$\mathbf{Z}\sim SSt_{n}(\boldsymbol{0},~\mathbf{I_{n}},~m,~\boldsymbol{\gamma}),$ $\mathbf{Z}^{\ast}\sim GSE_{n}\left(\boldsymbol{0},~\boldsymbol{I_{n}},~m,~\overline{G}_{n},~\boldsymbol{\gamma},~T\right)$. In addition, $\mathbf{M}^{\ast}$, $\mathbf{M}^{\ast\ast}$, $\mathbf{M^{\ast\ast}}_{-\xi_{\boldsymbol{s}k},\xi_{\boldsymbol{t}l}}$, $\mathbf{M^{\ast\ast}}_{-\xi_{\boldsymbol{s}k}}$ and $\mathbf{M^{\ast}}_{-\xi_{\boldsymbol{s}k}}$ are given in Corollary~2.\\
\noindent$\mathbf{Corollary~3}$ (GSLo distribution). Let $\mathbf{Y}\sim GSLo_{n}(\boldsymbol{\mu},~\mathbf{\Sigma},~\boldsymbol{\gamma},~J)$, where $J: \mathbb{R}\rightarrow\mathbb{R}$ and $\boldsymbol{\gamma}=(\gamma_{1},~\gamma_{2},~\cdots,~\gamma_{n})^{T}\in \mathbb{R}^{n}$.
In this case,
 $$g_{n}(u)=\frac{\exp\{-u\}}{[1+\exp\{-u\}]^{2}},~
c_{n}=\frac{1}{(2\pi)^{n/2}\Psi_{2}^{\ast}(-1,\frac{n}{2},1)},$$
  and $H\left(\mathbf{\Sigma}^{-\frac{1}{2}}(\boldsymbol{y}-\boldsymbol{\mu})\right)=J\left(\boldsymbol{\gamma}^{T}\mathbf{\Sigma}^{-\frac{1}{2}}(\boldsymbol{y}-\boldsymbol{\mu})\right)$.
 $\overline{G}_{n}(t)$ and $\overline{\mathcal{G}}_{n}(t)$ are expressed as
 $$\overline{G}_{n}(t)=\frac{\exp(-t)}{1+\exp(-t)},\;\;
 \overline{\mathcal{G}}_{n}(t)=\ln\left[1+\exp(-t)\right].$$ In addition,
 \begin{align*}
 c_{n}^{\ast}=\frac{1}{(2\pi)^{n/2}\Psi_{1}^{\ast}(-1,\frac{n}{2},1)}
\end{align*}
 and
 \begin{align*}
 c_{n}^{\ast\ast}=\frac{1}{(2\pi)^{n/2}\Psi_{1}^{\ast}(-1,\frac{n}{2}+1,1)}.
\end{align*}
Since
 $$f_{\mathbf{M^{\ast}}_{-\xi_{\boldsymbol{s}k}}}(\boldsymbol{w})=c_{n-1,\xi_{\boldsymbol{s}k}}^{\ast}\frac{\exp\left(-\frac{1}{2}\boldsymbol{w}^{T}\boldsymbol{w}-\frac{1}{2}\xi_{\boldsymbol{s},k}^{2}\right)}{\left[1+\exp\left(-\frac{1}{2}\boldsymbol{w}^{T}\boldsymbol{w}-\frac{1}{2}\xi_{\boldsymbol{s},k}^{2}\right)\right]^{2}},~\boldsymbol{w}\in \mathbb{R}^{n-1},$$
 $$f_{\mathbf{M^{\ast\ast}}_{-\xi_{\boldsymbol{s}k}}}(\boldsymbol{v})=c_{n-1,\xi_{\boldsymbol{s}k}}^{\ast\ast}\ln\left[1+\exp\left(-\frac{1}{2}\boldsymbol{v}^{T}\boldsymbol{v}-\frac{1}{2}\xi_{\boldsymbol{s},k}^{2}\right)\right],~\boldsymbol{v}\in \mathbb{R}^{n-1}$$
 and
$$f_{\mathbf{M^{\ast\ast}}_{-\xi_{\boldsymbol{s}k},\xi_{\boldsymbol{t}l}}}(\boldsymbol{u})=c_{n-2,\xi_{\boldsymbol{s}k},\xi_{\boldsymbol{t}l}}^{\ast\ast}\ln\left[1+\exp\left(-\frac{1}{2}\boldsymbol{u}^{T}\boldsymbol{u}-\frac{1}{2}\xi_{\boldsymbol{s},k}^{2}-\frac{1}{2}\xi_{\boldsymbol{t},l}^{2}\right)\right],~\boldsymbol{u}\in \mathbb{R}^{n-2},$$
so that
\begin{align*}
c_{n-1,\xi_{\boldsymbol{s}k}}^{\ast}&=\frac{\Gamma((n-1)/2)\exp\{\frac{\xi_{\boldsymbol{s},k}^{2}}{2}\}}{(2\pi)^{(n-1)/2}}\left[\int_{0}^{\infty}\frac{t^{(n-3)/2}\exp\{-t\}}{1+\exp\{-\frac{\xi_{\boldsymbol{s},k}^{2}}{2}\}\exp\{-t\}}\mathrm{d}t\right]^{-1}\\
&=\frac{\exp\{\frac{\xi_{\boldsymbol{s},k}^{2}}{2}\}}{(2\pi)^{(n-1)/2}\Psi_{1}^{\ast}\left(-\exp\{-\frac{\xi_{\boldsymbol{s},k}^{2}}{2}\},~\frac{n-1}{2},~1\right)},
\end{align*}
\begin{align*}
c_{n-1,\xi_{\boldsymbol{s}k}}^{\ast\ast}&=\frac{\Gamma((n-1)/2)}{(2\pi)^{(n-1)/2}}\left\{\int_{0}^{\infty}t^{(n-3)/2}\ln\left[1+\exp\left(-\frac{1}{2}\xi_{\boldsymbol{s},k}^{2}\right)\exp(-t)\right]\mathrm{d}t\right\}^{-1},
\end{align*}
and
\begin{align*}
&c_{n-2,\xi_{\boldsymbol{s}k},\xi_{\boldsymbol{t}l}}^{\ast\ast}=\frac{\Gamma((n-2)/2)}{(2\pi)^{(n-2)/2}}\left\{\int_{0}^{\infty}t^{(n-4)/2}\ln\left[1+\exp\left(-\frac{1}{2}\xi_{\boldsymbol{s},k}^{2}-\frac{1}{2}\xi_{\boldsymbol{t},l}^{2}\right)\exp(-t)\right]\mathrm{d}t\right\}^{-1}.
\end{align*}
Then
\begin{align*}
\nonumber\Omega_{ij}=& 2\bigg\{\frac{c_{n}}{c_{n-2,\xi_{\boldsymbol{a}i},\xi_{\boldsymbol{a}j}}^{\ast\ast}}\overline{\mathrm{E}}^{(\boldsymbol{\xi}_{\boldsymbol{a},-ij},\boldsymbol{\xi}_{\boldsymbol{b},-ij})}_{\boldsymbol{M}_{-\xi_{\boldsymbol{a}i},\xi_{\boldsymbol{a}j}}^{\ast\ast}}[J(\boldsymbol{\gamma}^{T}\boldsymbol{M}_{\xi_{\boldsymbol{a}i},\xi_{\boldsymbol{a}j}}^{\ast\ast})]-\frac{c_{n}}{c_{n-2,\xi_{\boldsymbol{a}i},\xi_{\boldsymbol{b}j}}^{\ast\ast}}\overline{\mathrm{E}}^{(\boldsymbol{\xi}_{\boldsymbol{a},-ij},\boldsymbol{\xi}_{\boldsymbol{b},-ij})}_{\boldsymbol{M}_{-\xi_{\boldsymbol{a}i},\xi_{\boldsymbol{b}j}}^{\ast\ast}}[J(\boldsymbol{\gamma}^{T}\boldsymbol{M}_{\xi_{\boldsymbol{a}i},\xi_{\boldsymbol{b}j}}^{\ast\ast})]\\
\nonumber&+\frac{c_{n}\gamma_{j}}{c_{n-1,\xi_{\boldsymbol{a}i}}^{\ast\ast}}\overline{\mathrm{E}}^{(\boldsymbol{\xi}_{\boldsymbol{a},-i},\boldsymbol{\xi}_{\boldsymbol{b},-i})}_{\boldsymbol{M}_{-\xi_{\boldsymbol{a}i}}^{\ast\ast}}[J'(\boldsymbol{\gamma}^{T}\boldsymbol{M}_{\xi_{\boldsymbol{a}i}}^{\ast\ast})]\\
\nonumber&+\frac{c_{n}}{c_{n-2,\xi_{\boldsymbol{b}i},\xi_{\boldsymbol{b}j}}^{\ast\ast}}\overline{\mathrm{E}}^{(\boldsymbol{\xi}_{\boldsymbol{a},-ij},\boldsymbol{\xi}_{\boldsymbol{b},-ij})}_{\boldsymbol{M}_{-\xi_{\boldsymbol{b}i},\xi_{\boldsymbol{b}j}}^{\ast\ast}}[J(\boldsymbol{\gamma}^{T}\boldsymbol{M}_{\xi_{\boldsymbol{b}i},\xi_{\boldsymbol{b}j}}^{\ast\ast})]-\frac{c_{n}}{c_{n-2,\xi_{\boldsymbol{b}i},\xi_{\boldsymbol{a}j}}^{\ast\ast}}\overline{\mathrm{E}}^{(\boldsymbol{\xi}_{\boldsymbol{a},-ij},\boldsymbol{\xi}_{\boldsymbol{b},-ij})}_{\boldsymbol{M}_{-\xi_{\boldsymbol{b}i},\xi_{\boldsymbol{a}j}}^{\ast\ast}}[J(\boldsymbol{\gamma}^{T}\boldsymbol{M}_{\xi_{\boldsymbol{b}i},\xi_{\boldsymbol{a}j}}^{\ast\ast})]\\
\nonumber&-\frac{c_{n}\gamma_{j}}{c_{n-1,\xi_{\boldsymbol{b}i}}^{\ast\ast}}\overline{\mathrm{E}}^{(\boldsymbol{\xi}_{\boldsymbol{a},-i},\boldsymbol{\xi}_{\boldsymbol{b},-i})}_{\boldsymbol{M}_{-\xi_{\boldsymbol{b}i}}^{\ast\ast}}[J'(\boldsymbol{\gamma}^{T}\boldsymbol{M}_{\xi_{\boldsymbol{b}i}}^{\ast\ast})]\\
\nonumber&+\frac{c_{n}\gamma_{i}}{c_{n-1,\xi_{\boldsymbol{a}j}}^{\ast\ast}}\overline{\mathrm{E}}^{(\boldsymbol{\xi}_{\boldsymbol{a},-j},\boldsymbol{\xi}_{\boldsymbol{b},-j})}_{\boldsymbol{M}_{-\xi_{\boldsymbol{a}j}}^{\ast\ast}}[J'(\boldsymbol{\gamma}^{T}\boldsymbol{M}_{\xi_{\boldsymbol{a}j}}^{\ast\ast})]-\frac{c_{n}\gamma_{i}}{c_{n-1,\xi_{\boldsymbol{b}j}}^{\ast\ast}}\overline{\mathrm{E}}^{(\boldsymbol{\xi}_{\boldsymbol{a},-j},\boldsymbol{\xi}_{\boldsymbol{b},-j})}_{\boldsymbol{M}_{-\xi_{\boldsymbol{b}j}}^{\ast\ast}}[J'(\boldsymbol{\gamma}^{T}\boldsymbol{M}_{\xi_{\boldsymbol{b}j}}^{\ast\ast})]\\
&+\frac{\Psi_{1}^{\ast}(-1,\frac{n}{2}+1,1)\gamma_{i}\gamma_{j}}{\Psi_{2}^{\ast}(-1,\frac{n}{2},1)}\overline{\mathrm{E}}^{(\boldsymbol{\xi}_{\boldsymbol{a}},\boldsymbol{\xi}_{\boldsymbol{b}})}_{\mathbf{M}^{\ast\ast}}[J''(\boldsymbol{\gamma}^{T}\boldsymbol{M}^{\ast\ast})]\bigg\},~i\neq j,
\end{align*}
\begin{align*}
\nonumber\Omega_{ii}=& 2\bigg\{\frac{c_{n}\xi_{\boldsymbol{a},i}}{c_{n-1,\xi_{\boldsymbol{a}i}}^{\ast}}\overline{\mathrm{E}}^{(\boldsymbol{\xi}_{\boldsymbol{a},-i},\boldsymbol{\xi}_{\boldsymbol{b},-i})}_{\boldsymbol{M}_{-\xi_{\boldsymbol{a}i}}^{\ast}}[J(\boldsymbol{\gamma}^{T}\boldsymbol{M}_{\xi_{\boldsymbol{a}i}}^{\ast})]-\frac{c_{n}\xi_{\boldsymbol{b},i}}{c_{n-1,\xi_{\boldsymbol{b}i}}^{\ast}}\overline{\mathrm{E}}^{(\boldsymbol{\xi}_{\boldsymbol{a},-i},\boldsymbol{\xi}_{\boldsymbol{b},-i})}_{\boldsymbol{M}_{-\xi_{\boldsymbol{b}i}}^{\ast}}[J(\boldsymbol{\gamma}^{T}\boldsymbol{M}_{\xi_{\boldsymbol{b}i}}^{\ast})]\\ \nonumber&+\frac{c_{n}\gamma_{i}}{c_{n-1,\xi_{\boldsymbol{a}i}}^{\ast\ast}}\overline{\mathrm{E}}^{(\boldsymbol{\xi}_{\boldsymbol{a},-i},\boldsymbol{\xi}_{\boldsymbol{b},-i})}_{\boldsymbol{M}_{-\xi_{\boldsymbol{a}i}}^{\ast\ast}}[J'(\boldsymbol{\gamma}^{T}\boldsymbol{M}_{\xi_{\boldsymbol{a}i}}^{\ast\ast})]-\frac{c_{n}\gamma_{i}}{c_{n-1,\xi_{\boldsymbol{b}i}}^{\ast\ast}}\overline{\mathrm{E}}^{(\boldsymbol{\xi}_{\boldsymbol{a},-i},\boldsymbol{\xi}_{\boldsymbol{b},-i})}_{\boldsymbol{M}_{-\xi_{\boldsymbol{b}i}}^{\ast\ast}}[J'(\boldsymbol{\gamma}^{T}\boldsymbol{M}_{\xi_{\boldsymbol{b}i}}^{\ast\ast})]\\
&+\frac{\Psi_{1}^{\ast}(-1,\frac{n}{2}+1,1)\gamma_{i}^{2}}{\Psi_{2}^{\ast}(-1,\frac{n}{2},1)}\overline{\mathrm{E}}^{(\boldsymbol{\xi}_{\boldsymbol{a}},\boldsymbol{\xi}_{\boldsymbol{b}})}_{\mathbf{M}^{\ast\ast}}[J''(\boldsymbol{\gamma}^{T}\boldsymbol{M}^{\ast\ast})]\bigg\}
+\frac{\Psi_{1}^{\ast}(-1,\frac{n}{2},1)}{\Psi_{2}^{\ast}(-1,\frac{n}{2},1)}F_{\mathbf{Z^{\ast}}}(\boldsymbol{\xi}_{\boldsymbol{a}},\boldsymbol{\xi}_{\boldsymbol{b}})
,
\end{align*}
\begin{align*}
\nonumber\delta_{k}=&2\bigg\{\frac{c_{n}}{c_{n-1,\xi_{\boldsymbol{a}k}}^{\ast}}\overline{\mathrm{E}}^{(\boldsymbol{\xi}_{\boldsymbol{a},-k},\boldsymbol{\xi}_{\boldsymbol{b},-k})}_{\boldsymbol{M}_{-\xi_{\boldsymbol{a}k}}^{\ast}}[J(\boldsymbol{\gamma}^{T}\boldsymbol{M}_{\xi_{\boldsymbol{a}k}}^{\ast})]-\frac{c_{n}}{c_{n-1,\xi_{\boldsymbol{b}k}}^{\ast}}\overline{\mathrm{E}}^{(\boldsymbol{\xi}_{\boldsymbol{a},-k},\boldsymbol{\xi}_{\boldsymbol{b},-k})}_{\boldsymbol{M}_{-\xi_{\boldsymbol{b}k}}^{\ast}}[J(\boldsymbol{\gamma}^{T}\boldsymbol{M}_{\xi_{\boldsymbol{b}k}}^{\ast})]\\
&+\frac{\Psi_{1}^{\ast}(-1,\frac{n}{2},1)\gamma_{k}}{\Psi_{2}^{\ast}(-1,\frac{n}{2},1)}\overline{\mathrm{E}}^{(\boldsymbol{\xi}_{\boldsymbol{a}},\boldsymbol{\xi}_{\boldsymbol{b}})}_{\mathbf{M}^{\ast}}[J'(\boldsymbol{\gamma}^{T}\boldsymbol{M}^{\ast})]\bigg\},
\end{align*}
 $i,~j,~k\in\{1,~2,\cdots,n\}$, $\mathbf{Z}\sim GSLo_{n}(\boldsymbol{0},~\mathbf{I_{n}},~\boldsymbol{\gamma},~J),$ $\mathbf{Z}^{\ast}\sim GSE_{n}\left(\boldsymbol{0},~\boldsymbol{I_{n}},~\overline{G}_{n},~\boldsymbol{\gamma},~J\right)$, $\mathbf{M}^{\ast}\sim E_{n}(\boldsymbol{0},~\mathbf{I_{n}},~\overline{G}_{n})$, \\ $\mathbf{M}^{\ast\ast}\sim E_{n}(\boldsymbol{0},~\mathbf{I_{n}},~\overline{\mathcal{G}}_{n})$ and $\mathbf{M}\sim Lo_{n}(\boldsymbol{0},~\mathbf{I_{n}})$.
 $\Psi_{\mu}^{\ast}(z,s,a)$ denotes generalized Hurwitz-Lerch zeta function (see, for instance,
Lin et al., 2006).\\
Therefore, $$\frac{c_{n}}{c_{n-1,\xi_{\boldsymbol{s}k}}^{\ast}}=\frac{\Psi_{1}^{\ast}\left(-\exp\{-\frac{\xi_{\boldsymbol{s},k}^{2}}{2}\},~\frac{n-1}{2},~1\right)\phi(\xi_{\boldsymbol{s},k})}{\Psi_{2}^{\ast}(-1,\frac{n}{2},1)},$$ $$\frac{c_{n}}{c_{n-1,\xi_{\boldsymbol{s}k}}^{\ast\ast}}=\frac{\int_{0}^{\infty}t^{(n-3)/2}\ln\left[1+\exp\left(-\frac{1}{2}\xi_{\boldsymbol{s},k}^{2}\right)\exp(-t)\right]\mathrm{d}t}{\Gamma((n-1)/2)\sqrt{2\pi}\Psi_{2}^{\ast}(-1,\frac{n}{2},1)}$$
and
$$\frac{c_{n}}{c_{n-2,\xi_{\boldsymbol{s}k},\xi_{\boldsymbol{t}l}}^{\ast\ast}}=\frac{\int_{0}^{\infty}t^{(n-4)/2}\ln\left[1+\exp\left(-\frac{1}{2}\xi_{\boldsymbol{s},k}^{2}-\frac{1}{2}\xi_{\boldsymbol{t},l}^{2}\right)\exp(-t)\right]\mathrm{d}t}{2\Gamma((n-1)/2)\pi\Psi_{2}^{\ast}(-1,\frac{n}{2},1)}.$$
$\mathbf{Example ~3}$ (Skew-logistic distribution) Letting $J(\cdot)=Lo(\cdot)$ in Corollary~3. Thus,
\begin{align*}
\nonumber\Omega_{ij}=& 2\bigg\{\frac{c_{n}}{c_{n-2,\xi_{\boldsymbol{a}i},\xi_{\boldsymbol{a}j}}^{\ast\ast}}\overline{\mathrm{E}}^{(\boldsymbol{\xi}_{\boldsymbol{a},-ij},\boldsymbol{\xi}_{\boldsymbol{b},-ij})}_{\boldsymbol{M}_{-\xi_{\boldsymbol{a}i},\xi_{\boldsymbol{a}j}}^{\ast\ast}}[Lo(\boldsymbol{\gamma}^{T}\boldsymbol{M}_{\xi_{\boldsymbol{a}i},\xi_{\boldsymbol{a}j}}^{\ast\ast})]-\frac{c_{n}}{c_{n-2,\xi_{\boldsymbol{a}i},\xi_{\boldsymbol{b}j}}^{\ast\ast}}\overline{\mathrm{E}}^{(\boldsymbol{\xi}_{\boldsymbol{a},-ij},\boldsymbol{\xi}_{\boldsymbol{b},-ij})}_{\boldsymbol{M}_{-\xi_{\boldsymbol{a}i},\xi_{\boldsymbol{b}j}}^{\ast\ast}}[Lo(\boldsymbol{\gamma}^{T}\boldsymbol{M}_{\xi_{\boldsymbol{a}i},\xi_{\boldsymbol{b}j}}^{\ast\ast})]\\
\nonumber&+\frac{c_{n}\gamma_{j}}{c_{n-1,\xi_{\boldsymbol{a}i}}^{\ast\ast}}\overline{\mathrm{E}}^{(\boldsymbol{\xi}_{\boldsymbol{a},-i},\boldsymbol{\xi}_{\boldsymbol{b},-i})}_{\boldsymbol{M}_{-\xi_{\boldsymbol{a}i}}^{\ast\ast}}[lo(\boldsymbol{\gamma}^{T}\boldsymbol{M}_{\xi_{\boldsymbol{a}i}}^{\ast\ast})]\\
\nonumber&+\frac{c_{n}}{c_{n-2,\xi_{\boldsymbol{b}i},\xi_{\boldsymbol{b}j}}^{\ast\ast}}\overline{\mathrm{E}}^{(\boldsymbol{\xi}_{\boldsymbol{a},-ij},\boldsymbol{\xi}_{\boldsymbol{b},-ij})}_{\boldsymbol{M}_{-\xi_{\boldsymbol{b}i},\xi_{\boldsymbol{b}j}}^{\ast\ast}}[Lo(\boldsymbol{\gamma}^{T}\boldsymbol{M}_{\xi_{\boldsymbol{b}i},\xi_{\boldsymbol{b}j}}^{\ast\ast})]-\frac{c_{n}}{c_{n-2,\xi_{\boldsymbol{b}i},\xi_{\boldsymbol{a}j}}^{\ast\ast}}\overline{\mathrm{E}}^{(\boldsymbol{\xi}_{\boldsymbol{a},-ij},\boldsymbol{\xi}_{\boldsymbol{b},-ij})}_{\boldsymbol{M}_{-\xi_{\boldsymbol{b}i},\xi_{\boldsymbol{a}j}}^{\ast\ast}}[Lo(\boldsymbol{\gamma}^{T}\boldsymbol{M}_{\xi_{\boldsymbol{b}i},\xi_{\boldsymbol{a}j}}^{\ast\ast})]\\
\nonumber&-\frac{c_{n}\gamma_{j}}{c_{n-1,\xi_{\boldsymbol{b}i}}^{\ast\ast}}\overline{\mathrm{E}}^{(\boldsymbol{\xi}_{\boldsymbol{a},-i},\boldsymbol{\xi}_{\boldsymbol{b},-i})}_{\boldsymbol{M}_{-\xi_{\boldsymbol{b}i}}^{\ast\ast}}[lo(\boldsymbol{\gamma}^{T}\boldsymbol{M}_{\xi_{\boldsymbol{b}i}}^{\ast\ast})]\\
\nonumber&+\frac{c_{n}\gamma_{i}}{c_{n-1,\xi_{\boldsymbol{a}j}}^{\ast\ast}}\overline{\mathrm{E}}^{(\boldsymbol{\xi}_{\boldsymbol{a},-j},\boldsymbol{\xi}_{\boldsymbol{b},-j})}_{\boldsymbol{M}_{-\xi_{\boldsymbol{a}j}}^{\ast\ast}}[lo(\boldsymbol{\gamma}^{T}\boldsymbol{M}_{\xi_{\boldsymbol{a}j}}^{\ast\ast})]-\frac{c_{n}\gamma_{i}}{c_{n-1,\xi_{\boldsymbol{b}j}}^{\ast\ast}}\overline{\mathrm{E}}^{(\boldsymbol{\xi}_{\boldsymbol{a},-j},\boldsymbol{\xi}_{\boldsymbol{b},-j})}_{\boldsymbol{M}_{-\xi_{\boldsymbol{b}j}}^{\ast\ast}}[lo(\boldsymbol{\gamma}^{T}\boldsymbol{M}_{\xi_{\boldsymbol{b}j}}^{\ast\ast})]\\
&+\frac{\Psi_{1}^{\ast}(-1,\frac{n}{2}+1,1)\gamma_{i}\gamma_{j}}{\Psi_{2}^{\ast}(-1,\frac{n}{2},1)}\overline{\mathrm{E}}^{(\boldsymbol{\xi}_{\boldsymbol{a}},\boldsymbol{\xi}_{\boldsymbol{b}})}_{\mathbf{M}^{\ast\ast}}\left[\boldsymbol{\gamma}^{T}\boldsymbol{M}^{\ast\ast}lo(\boldsymbol{\gamma}^{T}\boldsymbol{M}^{\ast\ast})\left(\frac{2\phi(\boldsymbol{\gamma}^{T}\boldsymbol{M}^{\ast\ast})}{(2\pi)^{-1/2}+\phi(\boldsymbol{\gamma}^{T}\boldsymbol{M}^{\ast\ast})}-1\right)\right]\bigg\},~i\neq j,
\end{align*}
\begin{align*}
\nonumber\Omega_{ii}=& 2\bigg\{\frac{c_{n}\xi_{\boldsymbol{a},i}}{c_{n-1,\xi_{\boldsymbol{a}i}}^{\ast}}\overline{\mathrm{E}}^{(\boldsymbol{\xi}_{\boldsymbol{a},-i},\boldsymbol{\xi}_{\boldsymbol{b},-i})}_{\boldsymbol{M}_{-\xi_{\boldsymbol{a}i}}^{\ast}}[Lo(\boldsymbol{\gamma}^{T}\boldsymbol{M}_{\xi_{\boldsymbol{a}i}}^{\ast})]-\frac{c_{n}\xi_{\boldsymbol{b},i}}{c_{n-1,\xi_{\boldsymbol{b}i}}^{\ast}}\overline{\mathrm{E}}^{(\boldsymbol{\xi}_{\boldsymbol{a},-i},\boldsymbol{\xi}_{\boldsymbol{b},-i})}_{\boldsymbol{M}_{-\xi_{\boldsymbol{b}i}}^{\ast}}[Lo(\boldsymbol{\gamma}^{T}\boldsymbol{M}_{\xi_{\boldsymbol{b}i}}^{\ast})]\\ \nonumber&+\frac{c_{n}\gamma_{i}}{c_{n-1,\xi_{\boldsymbol{a}i}}^{\ast\ast}}\overline{\mathrm{E}}^{(\boldsymbol{\xi}_{\boldsymbol{a},-i},\boldsymbol{\xi}_{\boldsymbol{b},-i})}_{\boldsymbol{M}_{-\xi_{\boldsymbol{a}i}}^{\ast\ast}}[lo(\boldsymbol{\gamma}^{T}\boldsymbol{M}_{\xi_{\boldsymbol{a}i}}^{\ast\ast})]-\frac{c_{n}\gamma_{i}}{c_{n-1,\xi_{\boldsymbol{b}i}}^{\ast\ast}}\overline{\mathrm{E}}^{(\boldsymbol{\xi}_{\boldsymbol{a},-i},\boldsymbol{\xi}_{\boldsymbol{b},-i})}_{\boldsymbol{M}_{-\xi_{\boldsymbol{b}i}}^{\ast\ast}}[lo(\boldsymbol{\gamma}^{T}\boldsymbol{M}_{\xi_{\boldsymbol{b}i}}^{\ast\ast})]\\
&+\frac{\Psi_{1}^{\ast}(-1,\frac{n}{2}+1,1)\gamma_{i}^{2}}{\Psi_{2}^{\ast}(-1,\frac{n}{2},1)}\overline{\mathrm{E}}^{(\boldsymbol{\xi}_{\boldsymbol{a}},\boldsymbol{\xi}_{\boldsymbol{b}})}_{\mathbf{M}^{\ast\ast}}\left[\boldsymbol{\gamma}^{T}\boldsymbol{M}^{\ast\ast}lo(\boldsymbol{\gamma}^{T}\boldsymbol{M}^{\ast\ast})\left(\frac{2\phi(\boldsymbol{\gamma}^{T}\boldsymbol{M}^{\ast\ast})}{(2\pi)^{-1/2}+\phi(\boldsymbol{\gamma}^{T}\boldsymbol{M}^{\ast\ast})}-1\right)\right]\bigg\}\\
&+\frac{\Psi_{1}^{\ast}(-1,\frac{n}{2},1)}{\Psi_{2}^{\ast}(-1,\frac{n}{2},1)}F_{\mathbf{Z^{\ast}}}(\boldsymbol{\xi}_{\boldsymbol{a}},\boldsymbol{\xi}_{\boldsymbol{b}})
,
\end{align*}
\begin{align*}
\nonumber\delta_{k}=&2\bigg\{\frac{c_{n}}{c_{n-1,\xi_{\boldsymbol{a}k}}^{\ast}}\overline{\mathrm{E}}^{(\boldsymbol{\xi}_{\boldsymbol{a},-k},\boldsymbol{\xi}_{\boldsymbol{b},-k})}_{\boldsymbol{M}_{-\xi_{\boldsymbol{a}k}}^{\ast}}[Lo(\boldsymbol{\gamma}^{T}\boldsymbol{M}_{\xi_{\boldsymbol{a}k}}^{\ast})]-\frac{c_{n}}{c_{n-1,\xi_{\boldsymbol{b}k}}^{\ast}}\overline{\mathrm{E}}^{(\boldsymbol{\xi}_{\boldsymbol{a},-k},\boldsymbol{\xi}_{\boldsymbol{b},-k})}_{\boldsymbol{M}_{-\xi_{\boldsymbol{b}k}}^{\ast}}[Lo(\boldsymbol{\gamma}^{T}\boldsymbol{M}_{\xi_{\boldsymbol{b}k}}^{\ast})]\\
&+\frac{\Psi_{1}^{\ast}(-1,\frac{n}{2},1)\gamma_{k}}{\Psi_{2}^{\ast}(-1,\frac{n}{2},1)}\overline{\mathrm{E}}^{(\boldsymbol{\xi}_{\boldsymbol{a}},\boldsymbol{\xi}_{\boldsymbol{b}})}_{\mathbf{M}^{\ast}}[lo(\boldsymbol{\gamma}^{T}\boldsymbol{M}^{\ast})]\bigg\},
\end{align*}
 $~i,j,k\in\{1,~2,\cdots,n\},$
$\mathbf{Z}\sim SLo_{n}(\boldsymbol{0},~\mathbf{I_{n}},~\boldsymbol{\gamma}),$ $\mathbf{Z}^{\ast}\sim GSE_{n}\left(\boldsymbol{0},~\boldsymbol{I_{n}},~\overline{G}_{n},~\boldsymbol{\gamma},~Lo\right)$.
Moreover, $\mathbf{M}^{\ast}$, $\mathbf{M}^{\ast\ast}$, $\mathbf{M^{\ast\ast}}_{-\xi_{\boldsymbol{s}k},\xi_{\boldsymbol{t}l}}$, $\mathbf{M^{\ast\ast}}_{-\xi_{\boldsymbol{s}k}}$ and $\mathbf{M^{\ast}}_{-\xi_{\boldsymbol{s}k}}$ are given in Corollary~3.\\
$\mathbf{Corollary~4}$ (GSLa distribution). Let $\mathbf{Y}\sim GSLa_{n}(\boldsymbol{\mu},~\mathbf{\Sigma},~\boldsymbol{\gamma},~J)$, where $J: \mathbb{R}\rightarrow\mathbb{R}$ and $\boldsymbol{\gamma}=(\gamma_{1},~\gamma_{2},~\cdots,~\gamma_{n})^{T}$.
     In this case,  $g_{n}(u)=\exp\{-\sqrt{2u}\},~c_{n}=\frac{\Gamma(n/2)}{2\pi^{n/2}\Gamma(n)}$ and  $$H\left(\mathbf{\Sigma}^{-\frac{1}{2}}(\boldsymbol{y}-\boldsymbol{\mu})\right)=J\left(\boldsymbol{\gamma}^{T}\mathbf{\Sigma}^{-\frac{1}{2}}(\boldsymbol{y}-\boldsymbol{\mu})\right).$$
  So
 $$\overline{G}_{n}(t)=(1+\sqrt{2t})\exp(-\sqrt{2t}),$$
 $$\overline{\mathcal{G}}_{n}(t)=(3+2t+3\sqrt{2t})\exp(-\sqrt{2t}).$$
 In addition,
$$c_{n}^{\ast}=\frac{n\Gamma(n/2)}{2\pi^{n/2}\Gamma(n+2)},\;\;
 c_{n}^{\ast\ast}=\frac{n(n+2)\Gamma(n/2)}{2\pi^{n/2}\Gamma(n+4)}.$$
  Since
 \begin{align*}
f_{\mathbf{M^{\ast}}_{-\xi_{\boldsymbol{s}k}}}(\boldsymbol{w})=c_{n-1,\xi_{\boldsymbol{s}k}}^{\ast}\left(1+\sqrt{\boldsymbol{w}^{T}\boldsymbol{w}+\xi_{\boldsymbol{s},k}^{2}}\right)\exp\left\{-\sqrt{\boldsymbol{w}^{T}\boldsymbol{w}+\xi_{\boldsymbol{s},k}^{2}}\right\},~\boldsymbol{w}\in \mathbb{R}^{n-1},
\end{align*}
 $$f_{\mathbf{M^{\ast\ast}}_{-\xi_{\boldsymbol{s}k}}}(\boldsymbol{v})=c_{n-1,\xi_{\boldsymbol{s}k}}^{\ast\ast}\left[\left(1+\frac{3}{\sqrt{2}}\right)\left(\boldsymbol{v}^{T}\boldsymbol{v}+\xi_{\boldsymbol{s},k}^{2}\right)+3\right]\exp\left(-\sqrt{\boldsymbol{v}^{T}\boldsymbol{v}+\xi_{\boldsymbol{s},k}^{2}}\right),~\boldsymbol{v}\in \mathbb{R}^{n-1}$$
 and
\begin{align*}
f_{\mathbf{M^{\ast\ast}}_{-\xi_{\boldsymbol{s}k},\xi_{\boldsymbol{t}l}}}(\boldsymbol{u})=&c_{n-2,ij}^{\ast\ast}\left[\left(1+\frac{3}{\sqrt{2}}\right)\left(\boldsymbol{u}^{T}\boldsymbol{u}+\xi_{\boldsymbol{s},k}^{2}+\xi_{\boldsymbol{t},l}^{2}\right)+3\right]\exp\left(-\sqrt{\boldsymbol{u}^{T}\boldsymbol{u}+\xi_{\boldsymbol{s},k}^{2}+\xi_{\boldsymbol{t},l}^{2}}\right),~\boldsymbol{u}\in \mathbb{R}^{n-2},
\end{align*}
thus
$$c_{n-1,\xi_{\boldsymbol{s}k}}^{\ast}=\frac{\Gamma\left(\frac{n-1}{2}\right)}{(2\pi)^{(n-1)/2}}\left[\int_{0}^{\infty}t^{\frac{n-3}{2}}\left(1+\sqrt{2t+\xi_{\boldsymbol{s},k}^{2}}\right)\exp\left\{-\sqrt{2t+\xi_{\boldsymbol{s},k}^{2}}\right\}\mathrm{d}t\right]^{-1},$$
\begin{align*}
c_{n-1,\xi_{\boldsymbol{s}k}}^{\ast\ast}=&\frac{\Gamma\left(\frac{n-1}{2}\right)}{(2\pi)^{(n-1)/2}}\bigg\{\int_{0}^{\infty}t^{\frac{n-3}{2}}\left[(2+3\sqrt{2})\left(t+\frac{\xi_{\boldsymbol{s},k}^{2}}{2}\right)+3\right]\exp\left(-\sqrt{2t+\xi_{\boldsymbol{s},k}^{2}}\right)\mathrm{d}t\bigg\}^{-1}
\end{align*}
and
\begin{align*}
c_{n-2,\xi_{\boldsymbol{s}k},\xi_{\boldsymbol{t}l}}^{\ast\ast}=&\frac{\Gamma\left(\frac{n-2}{2}\right)}{(2\pi)^{(n-2)/2}}\bigg\{\int_{0}^{\infty}t^{\frac{n-4}{2}}\left[(2+3\sqrt{2})\left(t+\frac{\xi_{\boldsymbol{s},k}^{2}}{2}+\frac{\xi_{\boldsymbol{t},l}^{2}}{2}\right)+3\right]\exp\left(-\sqrt{2t+\xi_{\boldsymbol{s},k}^{2}+\xi_{\boldsymbol{t},l}^{2}}\right)\mathrm{d}t\bigg\}^{-1}.
\end{align*}
Then
\begin{align*}
\nonumber\Omega_{ij}=& 2\bigg\{\frac{c_{n}}{c_{n-2,\xi_{\boldsymbol{a}i},\xi_{\boldsymbol{a}j}}^{\ast\ast}}\overline{\mathrm{E}}^{(\boldsymbol{\xi}_{\boldsymbol{a},-ij},\boldsymbol{\xi}_{\boldsymbol{b},-ij})}_{\boldsymbol{M}_{-\xi_{\boldsymbol{a}i},\xi_{\boldsymbol{a}j}}^{\ast\ast}}[J(\boldsymbol{\gamma}^{T}\boldsymbol{M}_{\xi_{\boldsymbol{a}i},\xi_{\boldsymbol{a}j}}^{\ast\ast})]-\frac{c_{n}}{c_{n-2,\xi_{\boldsymbol{a}i},\xi_{\boldsymbol{b}j}}^{\ast\ast}}\overline{\mathrm{E}}^{(\boldsymbol{\xi}_{\boldsymbol{a},-ij},\boldsymbol{\xi}_{\boldsymbol{b},-ij})}_{\boldsymbol{M}_{-\xi_{\boldsymbol{a}i},\xi_{\boldsymbol{b}j}}^{\ast\ast}}[J(\boldsymbol{\gamma}^{T}\boldsymbol{M}_{\xi_{\boldsymbol{a}i},\xi_{\boldsymbol{b}j}}^{\ast\ast})]\\
\nonumber&+\frac{c_{n}\gamma_{j}}{c_{n-1,\xi_{\boldsymbol{a}i}}^{\ast\ast}}\overline{\mathrm{E}}^{(\boldsymbol{\xi}_{\boldsymbol{a},-i},\boldsymbol{\xi}_{\boldsymbol{b},-i})}_{\boldsymbol{M}_{-\xi_{\boldsymbol{a}i}}^{\ast\ast}}[J'(\boldsymbol{\gamma}^{T}\boldsymbol{M}_{\xi_{\boldsymbol{a}i}}^{\ast\ast})]\\
\nonumber&+\frac{c_{n}}{c_{n-2,\xi_{\boldsymbol{b}i},\xi_{\boldsymbol{b}j}}^{\ast\ast}}\overline{\mathrm{E}}^{(\boldsymbol{\xi}_{\boldsymbol{a},-ij},\boldsymbol{\xi}_{\boldsymbol{b},-ij})}_{\boldsymbol{M}_{-\xi_{\boldsymbol{b}i},\xi_{\boldsymbol{b}j}}^{\ast\ast}}[J(\boldsymbol{\gamma}^{T}\boldsymbol{M}_{\xi_{\boldsymbol{b}i},\xi_{\boldsymbol{b}j}}^{\ast\ast})]-\frac{c_{n}}{c_{n-2,\xi_{\boldsymbol{b}i},\xi_{\boldsymbol{a}j}}^{\ast\ast}}\overline{\mathrm{E}}^{(\boldsymbol{\xi}_{\boldsymbol{a},-ij},\boldsymbol{\xi}_{\boldsymbol{b},-ij})}_{\boldsymbol{M}_{-\xi_{\boldsymbol{b}i},\xi_{\boldsymbol{a}j}}^{\ast\ast}}[J(\boldsymbol{\gamma}^{T}\boldsymbol{M}_{\xi_{\boldsymbol{b}i},\xi_{\boldsymbol{a}j}}^{\ast\ast})]\\
\nonumber&-\frac{c_{n}\gamma_{j}}{c_{n-1,\xi_{\boldsymbol{b}i}}^{\ast\ast}}\overline{\mathrm{E}}^{(\boldsymbol{\xi}_{\boldsymbol{a},-i},\boldsymbol{\xi}_{\boldsymbol{b},-i})}_{\boldsymbol{M}_{-\xi_{\boldsymbol{b}i}}^{\ast\ast}}[J'(\boldsymbol{\gamma}^{T}\boldsymbol{M}_{\xi_{\boldsymbol{b}i}}^{\ast\ast})]\\
\nonumber&+\frac{c_{n}\gamma_{i}}{c_{n-1,\xi_{\boldsymbol{a}j}}^{\ast\ast}}\overline{\mathrm{E}}^{(\boldsymbol{\xi}_{\boldsymbol{a},-j},\boldsymbol{\xi}_{\boldsymbol{b},-j})}_{\boldsymbol{M}_{-\xi_{\boldsymbol{a}j}}^{\ast\ast}}[J'(\boldsymbol{\gamma}^{T}\boldsymbol{M}_{\xi_{\boldsymbol{a}j}}^{\ast\ast})]-\frac{c_{n}\gamma_{i}}{c_{n-1,\xi_{\boldsymbol{b}j}}^{\ast\ast}}\overline{\mathrm{E}}^{(\boldsymbol{\xi}_{\boldsymbol{a},-j},\boldsymbol{\xi}_{\boldsymbol{b},-j})}_{\boldsymbol{M}_{-\xi_{\boldsymbol{b}j}}^{\ast\ast}}[J'(\boldsymbol{\gamma}^{T}\boldsymbol{M}_{\xi_{\boldsymbol{b}j}}^{\ast\ast})]\\
&+(n+3)(n+1)\gamma_{i}\gamma_{j}\overline{\mathrm{E}}^{(\boldsymbol{\xi}_{\boldsymbol{a}},\boldsymbol{\xi}_{\boldsymbol{b}})}_{\mathbf{M}^{\ast\ast}}[J''(\boldsymbol{\gamma}^{T}\boldsymbol{M}^{\ast\ast})]\bigg\},~i\neq j,
\end{align*}
\begin{align*}
\nonumber\Omega_{ii}=& 2\bigg\{\frac{c_{n}\xi_{\boldsymbol{a},i}}{c_{n-1,\xi_{\boldsymbol{a}i}}^{\ast}}\overline{\mathrm{E}}^{(\boldsymbol{\xi}_{\boldsymbol{a},-i},\boldsymbol{\xi}_{\boldsymbol{b},-i})}_{\boldsymbol{M}_{-\xi_{\boldsymbol{a}i}}^{\ast}}[J(\boldsymbol{\gamma}^{T}\boldsymbol{M}_{\xi_{\boldsymbol{a}i}}^{\ast})]-\frac{c_{n}\xi_{\boldsymbol{b},i}}{c_{n-1,\xi_{\boldsymbol{b}i}}^{\ast}}\overline{\mathrm{E}}^{(\boldsymbol{\xi}_{\boldsymbol{a},-i},\boldsymbol{\xi}_{\boldsymbol{b},-i})}_{\boldsymbol{M}_{-\xi_{\boldsymbol{b}i}}^{\ast}}[J(\boldsymbol{\gamma}^{T}\boldsymbol{M}_{\xi_{\boldsymbol{b}i}}^{\ast})]\\ \nonumber&+\frac{c_{n}\gamma_{i}}{c_{n-1,\xi_{\boldsymbol{a}i}}^{\ast\ast}}\overline{\mathrm{E}}^{(\boldsymbol{\xi}_{\boldsymbol{a},-i},\boldsymbol{\xi}_{\boldsymbol{b},-i})}_{\boldsymbol{M}_{-\xi_{\boldsymbol{a}i}}^{\ast\ast}}[J'(\boldsymbol{\gamma}^{T}\boldsymbol{M}_{\xi_{\boldsymbol{a}i}}^{\ast\ast})]-\frac{c_{n}\gamma_{i}}{c_{n-1,\xi_{\boldsymbol{b}i}}^{\ast\ast}}\overline{\mathrm{E}}^{(\boldsymbol{\xi}_{\boldsymbol{a},-i},\boldsymbol{\xi}_{\boldsymbol{b},-i})}_{\boldsymbol{M}_{-\xi_{\boldsymbol{b}i}}^{\ast\ast}}[J'(\boldsymbol{\gamma}^{T}\boldsymbol{M}_{\xi_{\boldsymbol{b}i}}^{\ast\ast})]\\
&+(n+3)(n+1)\gamma_{i}^{2}\overline{\mathrm{E}}^{(\boldsymbol{\xi}_{\boldsymbol{a}},\boldsymbol{\xi}_{\boldsymbol{b}})}_{\mathbf{M}^{\ast\ast}}[J''(\boldsymbol{\gamma}^{T}\boldsymbol{M}^{\ast\ast})]\bigg\}
+(n+1)F_{\mathbf{Z^{\ast}}}(\boldsymbol{\xi}_{\boldsymbol{a}},\boldsymbol{\xi}_{\boldsymbol{b}})
,
\end{align*}
\begin{align*}
\nonumber\delta_{k}=&2\bigg\{\frac{c_{n}}{c_{n-1,\xi_{\boldsymbol{a}k}}^{\ast}}\overline{\mathrm{E}}^{(\boldsymbol{\xi}_{\boldsymbol{a},-k},\boldsymbol{\xi}_{\boldsymbol{b},-k})}_{\boldsymbol{M}_{-\xi_{\boldsymbol{a}k}}^{\ast}}[J(\boldsymbol{\gamma}^{T}\boldsymbol{M}_{\xi_{\boldsymbol{a}k}}^{\ast})]-\frac{c_{n}}{c_{n-1,\xi_{\boldsymbol{b}k}}^{\ast}}\overline{\mathrm{E}}^{(\boldsymbol{\xi}_{\boldsymbol{a},-k},\boldsymbol{\xi}_{\boldsymbol{b},-k})}_{\boldsymbol{M}_{-\xi_{\boldsymbol{b}k}}^{\ast}}[J(\boldsymbol{\gamma}^{T}\boldsymbol{M}_{\xi_{\boldsymbol{b}k}}^{\ast})]\\
&+(n+1)\gamma_{k}\overline{\mathrm{E}}^{(\boldsymbol{\xi}_{\boldsymbol{a}},\boldsymbol{\xi}_{\boldsymbol{b}})}_{\mathbf{M}^{\ast}}[J'(\boldsymbol{\gamma}^{T}\boldsymbol{M}^{\ast})]\bigg\},
\end{align*}
$i,~j,~k\in\{1,~2,\cdots,n\}$, $\mathbf{Z}\sim GSLa_{n}(\boldsymbol{0},~\mathbf{I_{n}},~\boldsymbol{\gamma},~J),$ $\mathbf{Z}^{\ast}\sim GSE_{n}\left(\boldsymbol{0},~\boldsymbol{I_{n}},~\overline{G}_{n},~\boldsymbol{\gamma},~J\right)$, $\mathbf{M}^{\ast}\sim E_{n}(\boldsymbol{0},~\mathbf{I_{n}},~\overline{G}_{n})$,\\
$\mathbf{M}^{\ast\ast}\sim E_{n}(\boldsymbol{0},~\mathbf{I_{n}},~\overline{\mathcal{G}}_{n})$ and $\mathbf{M}\sim La_{n}(\boldsymbol{0},~\mathbf{I_{n}})$.
 $\frac{c_{n}}{c_{n-1,\xi_{\boldsymbol{s}k}}^{\ast}}$, $\frac{c_{n}}{c_{n-1,\xi_{\boldsymbol{s}k}}^{\ast\ast}}$ and $\frac{c_{n}}{c_{n-2,\xi_{\boldsymbol{s}k},\xi_{\boldsymbol{t}l}}^{\ast\ast}}$ can be further simplified as
$$\frac{c_{n}}{c_{n-1,\xi_{\boldsymbol{s}k}}^{\ast}}=\frac{\Gamma\left(\frac{n}{2}\right)2^{(n-3)/2}}{\Gamma(n)\Gamma\left(\frac{n-1}{2}\right)\sqrt{\pi}}\left[\int_{0}^{\infty}t^{\frac{n-3}{2}}\left(1+\sqrt{2t+\xi_{\boldsymbol{s},k}^{2}}\right)\exp\left\{-\sqrt{2t+\xi_{\boldsymbol{s},k}^{2}}\right\}\mathrm{d}t\right],$$
\begin{align*}
\frac{c_{n}}{c_{n-1,\xi_{\boldsymbol{s}k}}^{\ast\ast}}=&\frac{\Gamma\left(\frac{n}{2}\right)2^{(n-3)/2}}{\Gamma(n)\Gamma\left(\frac{n-1}{2}\right)\sqrt{\pi}}\bigg\{\int_{0}^{\infty}t^{\frac{n-3}{2}}\left[(2+3\sqrt{2})\left(t+\frac{\xi_{\boldsymbol{s},k}^{2}}{2}\right)+3\right]\exp\left(-\sqrt{2t+\xi_{\boldsymbol{s},k}^{2}}\right)\mathrm{d}t\bigg\}
\end{align*}
and
\begin{align*}
\frac{c_{n}}{c_{n-2,\xi_{\boldsymbol{s}k},\xi_{\boldsymbol{t}l}}^{\ast\ast}}=&\frac{2^{(n-6)/2}}{\Gamma(n-1)\pi}\bigg\{\int_{0}^{\infty}t^{\frac{n-4}{2}}\left[(2+3\sqrt{2})\left(t+\frac{\xi_{\boldsymbol{s},k}^{2}}{2}+\frac{\xi_{\boldsymbol{t},l}^{2}}{2}\right)+3\right]\exp\left(-\sqrt{2t+\xi_{\boldsymbol{s},k}^{2}+\xi_{\boldsymbol{t},l}^{2}}\right)\mathrm{d}t\bigg\}.
\end{align*}
$\mathbf{Example~4}$ (Skew-Laplace-normal distribution) Letting $J(\cdot)=\Phi(\cdot)$ in Corollary~4. Thus,
\begin{align*}
\nonumber\Omega_{ij}=& 2\bigg\{\frac{c_{n}}{c_{n-2,\xi_{\boldsymbol{a}i},\xi_{\boldsymbol{a}j}}^{\ast\ast}}\overline{\mathrm{E}}^{(\boldsymbol{\xi}_{\boldsymbol{a},-ij},\boldsymbol{\xi}_{\boldsymbol{b},-ij})}_{\boldsymbol{M}_{-\xi_{\boldsymbol{a}i},\xi_{\boldsymbol{a}j}}^{\ast\ast}}[\Phi(\boldsymbol{\gamma}^{T}\boldsymbol{M}_{\xi_{\boldsymbol{a}i},\xi_{\boldsymbol{a}j}}^{\ast\ast})]-\frac{c_{n}}{c_{n-2,\xi_{\boldsymbol{a}i},\xi_{\boldsymbol{b}j}}^{\ast\ast}}\overline{\mathrm{E}}^{(\boldsymbol{\xi}_{\boldsymbol{a},-ij},\boldsymbol{\xi}_{\boldsymbol{b},-ij})}_{\boldsymbol{M}_{-\xi_{\boldsymbol{a}i},\xi_{\boldsymbol{b}j}}^{\ast\ast}}[\Phi(\boldsymbol{\gamma}^{T}\boldsymbol{M}_{\xi_{\boldsymbol{a}i},\xi_{\boldsymbol{b}j}}^{\ast\ast})]\\
\nonumber&+\frac{c_{n}\gamma_{j}}{c_{n-1,\xi_{\boldsymbol{a}i}}^{\ast\ast}}\overline{\mathrm{E}}^{(\boldsymbol{\xi}_{\boldsymbol{a},-i},\boldsymbol{\xi}_{\boldsymbol{b},-i})}_{\boldsymbol{M}_{-\xi_{\boldsymbol{a}i}}^{\ast\ast}}[\phi(\boldsymbol{\gamma}^{T}\boldsymbol{M}_{\xi_{\boldsymbol{a}i}}^{\ast\ast})]\\
\nonumber&+\frac{c_{n}}{c_{n-2,\xi_{\boldsymbol{b}i},\xi_{\boldsymbol{b}j}}^{\ast\ast}}\overline{\mathrm{E}}^{(\boldsymbol{\xi}_{\boldsymbol{a},-ij},\boldsymbol{\xi}_{\boldsymbol{b},-ij})}_{\boldsymbol{M}_{-\xi_{\boldsymbol{b}i},\xi_{\boldsymbol{b}j}}^{\ast\ast}}[\Phi(\boldsymbol{\gamma}^{T}\boldsymbol{M}_{\xi_{\boldsymbol{b}i},\xi_{\boldsymbol{b}j}}^{\ast\ast})]-\frac{c_{n}}{c_{n-2,\xi_{\boldsymbol{b}i},\xi_{\boldsymbol{a}j}}^{\ast\ast}}\overline{\mathrm{E}}^{(\boldsymbol{\xi}_{\boldsymbol{a},-ij},\boldsymbol{\xi}_{\boldsymbol{b},-ij})}_{\boldsymbol{M}_{-\xi_{\boldsymbol{b}i},\xi_{\boldsymbol{a}j}}^{\ast\ast}}[\Phi(\boldsymbol{\gamma}^{T}\boldsymbol{M}_{\xi_{\boldsymbol{b}i},\xi_{\boldsymbol{a}j}}^{\ast\ast})]\\
\nonumber&-\frac{c_{n}\gamma_{j}}{c_{n-1,\xi_{\boldsymbol{b}i}}^{\ast\ast}}\overline{\mathrm{E}}^{(\boldsymbol{\xi}_{\boldsymbol{a},-i},\boldsymbol{\xi}_{\boldsymbol{b},-i})}_{\boldsymbol{M}_{-\xi_{\boldsymbol{b}i}}^{\ast\ast}}[\phi(\boldsymbol{\gamma}^{T}\boldsymbol{M}_{\xi_{\boldsymbol{b}i}}^{\ast\ast})]\\
\nonumber&+\frac{c_{n}\gamma_{i}}{c_{n-1,\xi_{\boldsymbol{a}j}}^{\ast\ast}}\overline{\mathrm{E}}^{(\boldsymbol{\xi}_{\boldsymbol{a},-j},\boldsymbol{\xi}_{\boldsymbol{b},-j})}_{\boldsymbol{M}_{-\xi_{\boldsymbol{a}j}}^{\ast\ast}}[\phi(\boldsymbol{\gamma}^{T}\boldsymbol{M}_{\xi_{\boldsymbol{a}j}}^{\ast\ast})]-\frac{c_{n}\gamma_{i}}{c_{n-1,\xi_{\boldsymbol{b}j}}^{\ast\ast}}\overline{\mathrm{E}}^{(\boldsymbol{\xi}_{\boldsymbol{a},-j},\boldsymbol{\xi}_{\boldsymbol{b},-j})}_{\boldsymbol{M}_{-\xi_{\boldsymbol{b}j}}^{\ast\ast}}[\phi(\boldsymbol{\gamma}^{T}\boldsymbol{M}_{\xi_{\boldsymbol{b}j}}^{\ast\ast})]\\
&-(n+3)(n+1)\gamma_{i}\gamma_{j}\overline{\mathrm{E}}^{(\boldsymbol{\xi}_{\boldsymbol{a}},\boldsymbol{\xi}_{\boldsymbol{b}})}_{\mathbf{M}^{\ast\ast}}[\boldsymbol{\gamma}^{T}\boldsymbol{M}^{\ast\ast}\phi(\boldsymbol{\gamma}^{T}\boldsymbol{M}^{\ast\ast})]\bigg\},~i\neq j,
\end{align*}
\begin{align*}
\nonumber\Omega_{ii}=& 2\bigg\{\frac{c_{n}\xi_{\boldsymbol{a},i}}{c_{n-1,\xi_{\boldsymbol{a}i}}^{\ast}}\overline{\mathrm{E}}^{(\boldsymbol{\xi}_{\boldsymbol{a},-i},\boldsymbol{\xi}_{\boldsymbol{b},-i})}_{\boldsymbol{M}_{-\xi_{\boldsymbol{a}i}}^{\ast}}[\Phi(\boldsymbol{\gamma}^{T}\boldsymbol{M}_{\xi_{\boldsymbol{a}i}}^{\ast})]-\frac{c_{n}\xi_{\boldsymbol{b},i}}{c_{n-1,\xi_{\boldsymbol{b}i}}^{\ast}}\overline{\mathrm{E}}^{(\boldsymbol{\xi}_{\boldsymbol{a},-i},\boldsymbol{\xi}_{\boldsymbol{b},-i})}_{\boldsymbol{M}_{-\xi_{\boldsymbol{b}i}}^{\ast}}[\Phi(\boldsymbol{\gamma}^{T}\boldsymbol{M}_{\xi_{\boldsymbol{b}i}}^{\ast})]\\ \nonumber&+\frac{c_{n}\gamma_{i}}{c_{n-1,\xi_{\boldsymbol{a}i}}^{\ast\ast}}\overline{\mathrm{E}}^{(\boldsymbol{\xi}_{\boldsymbol{a},-i},\boldsymbol{\xi}_{\boldsymbol{b},-i})}_{\boldsymbol{M}_{-\xi_{\boldsymbol{a}i}}^{\ast\ast}}[\phi(\boldsymbol{\gamma}^{T}\boldsymbol{M}_{\xi_{\boldsymbol{a}i}}^{\ast\ast})]-\frac{c_{n}\gamma_{i}}{c_{n-1,\xi_{\boldsymbol{b}i}}^{\ast\ast}}\overline{\mathrm{E}}^{(\boldsymbol{\xi}_{\boldsymbol{a},-i},\boldsymbol{\xi}_{\boldsymbol{b},-i})}_{\boldsymbol{M}_{-\xi_{\boldsymbol{b}i}}^{\ast\ast}}[\phi(\boldsymbol{\gamma}^{T}\boldsymbol{M}_{\xi_{\boldsymbol{b}i}}^{\ast\ast})]\\
&-(n+3)(n+1)\gamma_{i}^{2}\overline{\mathrm{E}}^{(\boldsymbol{\xi}_{\boldsymbol{a}},\boldsymbol{\xi}_{\boldsymbol{b}})}_{\mathbf{M}^{\ast\ast}}[\boldsymbol{\gamma}^{T}\boldsymbol{M}^{\ast\ast}\phi(\boldsymbol{\gamma}^{T}\boldsymbol{M}^{\ast\ast})]\bigg\}
+(n+1)F_{\mathbf{Z^{\ast}}}(\boldsymbol{\xi}_{\boldsymbol{a}},\boldsymbol{\xi}_{\boldsymbol{b}})
,
\end{align*}
\begin{align*}
\nonumber\delta_{k}=&2\bigg\{\frac{c_{n}}{c_{n-1,\xi_{\boldsymbol{a}k}}^{\ast}}\overline{\mathrm{E}}^{(\boldsymbol{\xi}_{\boldsymbol{a},-k},\boldsymbol{\xi}_{\boldsymbol{b},-k})}_{\boldsymbol{M}_{-\xi_{\boldsymbol{a}k}}^{\ast}}[\Phi(\boldsymbol{\gamma}^{T}\boldsymbol{M}_{\xi_{\boldsymbol{a}k}}^{\ast})]-\frac{c_{n}}{c_{n-1,\xi_{\boldsymbol{b}k}}^{\ast}}\overline{\mathrm{E}}^{(\boldsymbol{\xi}_{\boldsymbol{a},-k},\boldsymbol{\xi}_{\boldsymbol{b},-k})}_{\boldsymbol{M}_{-\xi_{\boldsymbol{b}k}}^{\ast}}[\Phi(\boldsymbol{\gamma}^{T}\boldsymbol{M}_{\xi_{\boldsymbol{b}k}}^{\ast})]\\
&+(n+1)\gamma_{k}\overline{\mathrm{E}}^{(\boldsymbol{\xi}_{\boldsymbol{a}},\boldsymbol{\xi}_{\boldsymbol{b}})}_{\mathbf{M}^{\ast}}[\phi(\boldsymbol{\gamma}^{T}\boldsymbol{M}^{\ast})]\bigg\},
\end{align*}
$i,~j,~k\in\{1,~2,\cdots,n\}$,
$\mathbf{Z}\sim SLaN_{n}(\boldsymbol{0},~\mathbf{I_{n}},~\boldsymbol{\gamma}),$ $\mathbf{Z}^{\ast}\sim GSE_{n}\left(\boldsymbol{0},~\boldsymbol{I_{n}},~\overline{G}_{n},~\boldsymbol{\gamma},~\Phi\right)$. In addition, $\mathbf{M}^{\ast}$, $\mathbf{M}^{\ast\ast}$, $\mathbf{M^{\ast\ast}}_{-\xi_{\boldsymbol{s}k},\xi_{\boldsymbol{t}l}}$, $\mathbf{M^{\ast\ast}}_{-\xi_{\boldsymbol{s}k}}$ and $\mathbf{M^{\ast}}_{-\xi_{\boldsymbol{s}k}}$ are given in Corollary~4.
\section{ Multivariate tail conditional risk measures for GSE distributions}
In this section, we mainly consider multivariate tail expectation (MTCE) and multivariate tail covariance (MTCov) for GSE distributions.

Landsman et al. (2018) proposed a new multivariate tail conditional expectation (MTCE) for an $n\times1$ vector of risks $\mathbf{X}=(X_{1},~X_{2},\cdots,X_{n})^{T}$ with cdf $F_{\mathbf{X}}(\boldsymbol{x})$:
\begin{align}\label{(19)}
\nonumber\mathrm{MTCE}_{\boldsymbol{q}}(\mathbf{X})&=\mathrm{E}\left[\mathbf{X}|\mathbf{X}>VaR_{\boldsymbol{q}}(\mathbf{X})\right]\\
&=\mathrm{E}[\mathbf{X}|X_{1}>VaR_{q_{1}}(X_{1}),\cdots,X_{n}>VaR_{q_{n}}(X_{n})],
\end{align}
 where
 $\boldsymbol{q}=(q_{1},\cdots,q_{n})\in(0,~1)^{n},$
 $VaR_{\boldsymbol{q}}(\mathbf{X})=(VaR_{q_{1}}(X_{1}),~VaR_{q_{2}}(X_{2}),\cdots,VaR_{q_{n}}(X_{n}))^{T},$ and $VaR_{q_{k}}(X_{k})$,~$k\in\{1,~2,\cdots,n\}$, denotes the $q_{k}$-th quantile of $X_{k}$. Specially, Landsman et al. (2016) is a special case as $\boldsymbol{q}=(q,\cdots,q)\in(0,~1)^{n}.$
Landsman et al. (2018) further proposed a novel form of multivariate tail covariance (MTCov):
\begin{align}\label{(20)}
\nonumber\mathrm{MTCov}_{\boldsymbol{q}}(\mathbf{X})&=\mathrm{E}\left[(\mathbf{X}-\mathrm{MTCE}_{\boldsymbol{q}}(\mathbf{X}))(\mathbf{X}-\mathrm{MTCE}_{\boldsymbol{q}}(\mathbf{X}))^{T}|\mathbf{X}>VaR_{\boldsymbol{q}}(\mathbf{X})\right]\\
&=\inf_{\boldsymbol{c}\in\mathbb{R}^{n}}\mathrm{E}\left[(\mathbf{X}-\boldsymbol{c})(\mathbf{X}-\boldsymbol{c})^{T}|\mathbf{X}>VaR_{\boldsymbol{q}}(\mathbf{X})\right].
\end{align}

From Proposition 1, we give following corollary.\\
$\mathbf{Corollary~5}$ Let $\mathbf{Y}\sim GSE_{n}(\boldsymbol{\mu},~\mathbf{\Sigma},~g_{n},~H)$ $(n\geq2)$ be as in $(\ref{(1)})$. Suppose that it satisfies conditions
\begin{align}\label{(21)}
\lim_{z_{k}\rightarrow +\infty}z_{k}H(\boldsymbol{z})\overline{G}_{n}\left(\frac{1}{2}\boldsymbol{z}^{T}\boldsymbol{z}\right)=0,~k\in\{1,2,\cdots,n\}
\end{align}
and
\begin{align}\label{(22)}
\lim_{z_{k}\rightarrow +\infty}H(\boldsymbol{z})\overline{\mathcal{G}}_{n}\left(\frac{1}{2}\boldsymbol{z}^{T}\boldsymbol{z}\right)=0,~k\in\{1,2,\cdots,n\}.
\end{align}
Further, assume $\partial_{i}H$ and $\partial_{ij}H$ exist for $i,j\in\{1,2,\cdots,n\}$.
Then
\begin{align}\label{(23)}
\mathrm{(\widetilde{I})}~\mathrm{MTE}_{\boldsymbol{q}}(\mathbf{Y})&=\boldsymbol{\mu}+\frac{\mathbf{\Sigma}^{\frac{1}{2}}\boldsymbol{\delta}}{\overline{F}_{\mathbf{Z}}(\boldsymbol{\xi_{q}})},
\end{align}
\begin{align}\label{(24)}
\mathrm{(\widetilde{II})}~&\mathrm{MDTCov}_{(\boldsymbol{a},\boldsymbol{b})}(\mathbf{Y})=\mathbf{\Sigma}^{\frac{1}{2}}\left[ \frac{\boldsymbol{\Omega}}{\overline{F}_{\mathbf{Z}}(\boldsymbol{\xi_{q}})}-\frac{\boldsymbol{\delta}\boldsymbol{\delta}^{T}}{\overline{F}_{\mathbf{Z}}^{2}(\boldsymbol{\xi_{q}})} \right]  \mathbf{\Sigma}^{\frac{1}{2}},
\end{align}
where $\mathbf{\Omega}=(\Omega_{ij})_{i,j=1}^{n}$ is an $n\times n$ symmetric matrix, and $\boldsymbol{\delta}=(\delta_{1},\delta_{2},\cdots,\delta_{n})^{T}$ is an $n\times1$ vector.
Here \begin{align*}
\nonumber\Omega_{ij}=& 2\bigg\{\frac{c_{n}}{c_{n-2,\xi_{\boldsymbol{q}i},\xi_{\boldsymbol{q}j}}^{\ast\ast}}\overline{\mathrm{E}}^{\boldsymbol{\xi}_{\boldsymbol{q},-ij}}_{\boldsymbol{M}_{-\xi_{\boldsymbol{q}i},\xi_{\boldsymbol{q}j}}^{\ast\ast}}[H(\boldsymbol{M}_{\xi_{\boldsymbol{q}i},\xi_{\boldsymbol{q}j}}^{\ast\ast})]+\frac{c_{n}}{c_{n-1,\xi_{\boldsymbol{q}i}}^{\ast\ast}}\overline{\mathrm{E}}^{\boldsymbol{\xi}_{\boldsymbol{q},-i}}_{\boldsymbol{M}_{-\xi_{\boldsymbol{q}i}}^{\ast\ast}}[\partial_{j}H(\boldsymbol{M}_{\xi_{\boldsymbol{q}i}}^{\ast\ast})]\\
&+\frac{c_{n}}{c_{n-1,\xi_{\boldsymbol{q}j}}^{\ast\ast}}\overline{\mathrm{E}}^{\boldsymbol{\xi}_{\boldsymbol{q},-j}}_{\boldsymbol{M}_{-\xi_{\boldsymbol{q}j}}^{\ast\ast}}[\partial_{i}H(\boldsymbol{M}_{\xi_{\boldsymbol{q}j}}^{\ast\ast})]+\frac{c_{n}}{c_{n}^{\ast\ast}}\overline{\mathrm{E}}^{\boldsymbol{\xi}_{\boldsymbol{q}}}_{\mathbf{M}^{\ast\ast}}[\partial_{ij}H(\boldsymbol{M}^{\ast\ast})]\bigg\},~i\neq j,
\end{align*}
\begin{align*}
\nonumber\Omega_{ii}=& 2\bigg\{\frac{c_{n}\xi_{\boldsymbol{q},i}}{c_{n-1,\xi_{\boldsymbol{q}i}}^{\ast}}\overline{\mathrm{E}}^{\boldsymbol{\xi}_{\boldsymbol{q},-i}}_{\boldsymbol{M}_{-\xi_{\boldsymbol{q}i}}^{\ast}}[H(\boldsymbol{M}_{\xi_{\boldsymbol{q}i}}^{\ast})]+\frac{c_{n}}{c_{n-1,\xi_{\boldsymbol{q}i}}^{\ast\ast}}\overline{\mathrm{E}}^{\boldsymbol{\xi}_{\boldsymbol{q},-i}}_{\boldsymbol{M}_{-\xi_{\boldsymbol{q}i}}^{\ast\ast}}[\partial_{i}H(\boldsymbol{M}_{\xi_{\boldsymbol{q}i}}^{\ast\ast})]+\frac{c_{n}}{c_{n}^{\ast\ast}}\overline{\mathrm{E}}^{\boldsymbol{\xi}_{\boldsymbol{q}}}_{\mathbf{M}^{\ast\ast}}[\partial_{ii}H(\boldsymbol{M}^{\ast\ast})]\bigg\}
+\frac{c_{n}}{c_{n}^{\ast}}\overline{F}_{\mathbf{Z^{\ast}}}(\boldsymbol{\xi}_{\boldsymbol{q}})
,
\end{align*}
\begin{align*}
\delta_{k}=&2\bigg\{\frac{c_{n}}{c_{n-1,\xi_{\boldsymbol{q}k}}^{\ast}}\overline{\mathrm{E}}^{\boldsymbol{\xi}_{\boldsymbol{q},-k}}_{\boldsymbol{M}_{-\xi_{\boldsymbol{q}k}}^{\ast}}[H(\boldsymbol{M}_{\xi_{\boldsymbol{q}k}}^{\ast})]+\frac{c_{n}}{c_{n}^{\ast}}\overline{\mathrm{E}}^{\boldsymbol{\xi}_{\boldsymbol{q}}}_{\mathbf{M}^{\ast}}[\partial_{k}H(\boldsymbol{M}^{\ast})]\bigg\},
\end{align*}
 $~i,j,k\in\{1,~2,\cdots,n\},$ $\boldsymbol{\xi_{q}}=\left(\xi_{\boldsymbol{q},1},~\xi_{\boldsymbol{q},2},\cdots,\xi_{\boldsymbol{q},n}\right)^{T}=\mathbf{\Sigma}^{-\frac{1}{2}}(VaR_{\boldsymbol{q}}(\boldsymbol{Y})-\boldsymbol{\mu}),$ $\mathbf{Z}$, $\mathbf{M}^{\ast}$,$\mathbf{Z}^{\ast}$, $\mathbf{M}^{\ast\ast}$,  $\mathbf{M^{\ast}}_{-\xi_{\boldsymbol{s}k}}$, $\mathbf{M^{\ast\ast}}_{-\xi_{\boldsymbol{s}k}}$ and $\mathbf{M^{\ast\ast}}_{-\xi_{\boldsymbol{s}k},\xi_{\boldsymbol{t}l}}$
   are given in  Theorem 1. In addition, $\overline{F}_{\mathbf{W}}(\cdot)$ is tail function of $\mathbf{W}$.\\
\noindent $\mathbf{Proof}$. Letting $\boldsymbol{a}=VaR_{\boldsymbol{q}}(\boldsymbol{Y})$ and $\boldsymbol{b\rightarrow+\infty}$ in Proposition 1, and combining with conditions (\ref{(21)}) and (\ref{(22)}), we can instantly obtain above results.

Note that the results $\mathrm{(\widetilde{I})}$ and $\mathrm{(\widetilde{II})}$ of Theorem 2 are coincide with Theorem 3.1 and Theorem 1 in Zuo and Yin (2020,~2021b), respectively.
\section{Numerical illustration}
In this section, We examine multivariate doubly truncated expectation, multivariate doubly truncated covariance and MTCE for skew-normal (SN) distributions.\\
 Now, we compute multivariate doubly truncated expectation, multivariate doubly truncated covariance for SN distributions.
  Let
$\boldsymbol{a}=\left(2,2\right)^{T},~\boldsymbol{b}=\left(6,7\right)^{T}$, $\mathbf{P}=(P_{1},P_{2})^{T}\sim SN_{2}(\boldsymbol{\mu},\mathbf{\Sigma},\boldsymbol{\gamma})$ with parameters
\begin{align*}
\boldsymbol{\mu}=\left(\begin{array}{ccccccccccc}
3\\
4
\end{array}
\right),
\mathbf{\Sigma}=\left(\begin{array}{ccccccccccc}
0.9&0.5\\
0.5&0.7
\end{array}
\right),~
\boldsymbol{\gamma}=\left(\begin{array}{ccccccccccc}
1\\
2
\end{array}
\right).
\end{align*}
From Proposition 1 and Example 1, firstly, we give $\boldsymbol{\delta}$ and $\boldsymbol{\Omega}$:
\begin{align*}
\boldsymbol{\delta}=\left(\begin{array}{ccccccccccc}
0.4172\\
0.4875
\end{array}
\right),
\boldsymbol{\Omega}=\left(\begin{array}{ccccccccccc}
0.7886&0.0993\\
0.0993&0.6688
\end{array}
\right).
\end{align*}
Then the multivariate doubly truncated expectation and multivariate doubly truncated covariance matrix of $\mathbf{P}$ for $(\boldsymbol{a},~\boldsymbol{b})$ can be obtained as
\begin{align*}
\mathrm{MDTE}_{(\boldsymbol{a},\boldsymbol{b})}(\mathbf{P})=\left(\begin{array}{ccccccccccc}
3.4459\\
4.6047
\end{array}
\right),
\mathrm{MDTCov}_{(\boldsymbol{a},\boldsymbol{b})}(\mathbf{P})=\left(\begin{array}{ccccccccccc}
0.6008&0.2436\\
0.2436&0.2760
\end{array}
\right).
\end{align*}
Next,  we compute MTCE for skew-normal (SN) distributions. Let $$\boldsymbol{q_{1}}=(0.90,0.90,0.90,0.90,0.90)^{T},~\boldsymbol{q_{2}}=(0.90,0.80,0.70,0.60,0.50)^{T},$$ $\mathbf{Q}=(Q_{1},Q_{2},Q_{3},Q_{4},Q_{5})^{T}\sim SN_{5}(\boldsymbol{\mu},\mathbf{\Sigma},\boldsymbol{\gamma})$ with parameters
\begin{align*}
\boldsymbol{\mu}=\left(\begin{array}{ccccccccccccccccc}
0.4\\
 0.1\\
  1.4\\
   1.1\\
    1.0
\end{array}
\right)
,
\mathbf{\Sigma}=\left(\begin{array}{ccccccccccccccccccc}
0.84 &-0.27 &0.13 &-0.09 &-0.08\\
-0.27 &1.32 &0.11 &-0.07 &-0.15\\
0.13 &0.11 &0.82 &0.02 &-0.02\\
-0.09 &-0.07 &0.02 &0.64 &0.22\\
-0.08 &-0.15 &-0.02 &0.22 &0.67
\end{array}
\right),~
\boldsymbol{\gamma}=\left(\begin{array}{ccccccccccccccc}
0.1\\
 0.2\\
 -0.3\\
 0.1\\
  -0.2
\end{array}
\right).
\end{align*}
From Corollary~5 and Example 1, VaRs and MTCEs of $\mathbf{Q}$ for $\boldsymbol{q_{1}}$ and $\boldsymbol{q_{2}}$ are shown in Tables 1 and 2, respectively:
\begin{table}[!htbp]
\centering
\begin{tabular}{|c||c|c|c|c|c|c|}
  \hline
  {$i$}&$1$& $2$&$3$&4&5\\
  \hline
  $VaR_{\boldsymbol{q}_{1}}(\mathbf{Q})_{i}$ &1.6435 & 1.7330& 2.3206&2.1851&1.9080\\
  \hline
  $VaR_{\boldsymbol{q}_{2}}(\mathbf{Q})_{i}$&1.6435 & 1.2370&  1.6554&1.3691&0.8722\\
  \hline
\end{tabular}\\
Table 1: The VaRs of $\mathbf{Q}$ for $\boldsymbol{q_{1}}=(0.90, 0.90, 0.90,0.90,0.90)^T$ and $\boldsymbol{q_{2}}=(0.90, 0.80, 0.70,0.60,0.50)^T$.
\end{table}
\begin{table}[!htbp]
\centering
\begin{tabular}{|c||c|c|c|c|c|c|}
  \hline
  {$i$}&$1$& $2$&$3$&4&5\\
  \hline
  $\mathrm{MTCE}_{\boldsymbol{q_{1}}}(\mathbf{Q})_{i}$ &7.0066 & 10.6656& -3.0526&5.1563&-2.2325\\
  \hline
  $\mathrm{MTCE}_{\boldsymbol{q_{2}}}(\mathbf{Q})_{i}$&3.8282 & 3.3809&  3.0977&2.5537&1.6741\\
  \hline
\end{tabular}\\
Table 2: The MTCEs of $\mathbf{Q}$ for $\boldsymbol{q_{1}}=(0.90, 0.90, 0.90,0.90,0.90)^T$ and $\boldsymbol{q_{2}}=(0.90, 0.80, 0.70,0.60,0.50)^T$.
\end{table}
 As we see in Table 2 and Figure 1, there is a clear difference between the $\mathrm{MTCE}_{\boldsymbol{q_{1}}}(\mathbf{Q})$ and $\mathrm{MTCE}_{\boldsymbol{q_{2}}}(\mathbf{Q})$. The MTCE measures of Risks 1, 2 and 4 for $\boldsymbol{q_{1}}$  are greater than the MTCE measure of corresponding Risks for $\boldsymbol{q_{2}}$. However, Risks 3 and 5 are opposite.
 \section{Concluding remarks}
 In this paper, we have studied multivariate doubly truncated first two moments
of GSE distributions, which provides further
generalisation of the moments for doubly truncated multivariate normal mean-variance mixture distributions (Roozegar et al., 2020). Several important cases, for examples, GSN, GSSt, GSLo and GSLa distributions, are given. We also have presented multivariate doubly truncated expectation (MDTE) and covariance (MDTCov) for GSE distributions providing further
generalisation of the MDTE and MDTCov for elliptical distributions, considered in Zuo and Yin (2021c).
Note that when $\boldsymbol{a\rightarrow-\infty}$, those doubly truncated moments will degenerate into lower truncated moments; When $\boldsymbol{b\rightarrow+\infty}$, those doubly truncated moments will degenerate into upper truncated moments; When $\boldsymbol{a\rightarrow-\infty}$ and $\boldsymbol{b\rightarrow+\infty}$, those doubly truncated moments will degenerate into usual moments. As applications of our results, the MTCE and MTCov for generalized skew-elliptical distributions are given. Aim to examine established results, we have presented the numerical illustrations.
 Moreover, Ignatieva and Landsman (2021) computed tail conditional expectation (TCE) for generalised hyper-elliptical distributions. In Zuo and Yin (2021a), the authors derived formula of MTCE for location-scale mixture of elliptical distributions. It will, of course, be
of interest to generalize the results (established here) to those mixture distributions and we hope to report the findings in a future paper.
\section*{Acknowledgments}
\noindent  The research was supported by the National Natural Science Foundation of China (No. 12071251)
\section*{Conflicts of Interest}
\noindent The authors declare that they have no conflicts of interest.
\section*{References}
\bibliographystyle{model1-num-names}







\end{document}